\documentclass[11pt, a4paper]{amsart}

\usepackage{amsfonts,amsmath,amssymb, amscd,fullpage}
\usepackage{stmaryrd}
\usepackage[all]{xy}

\newtheorem{theorem}{Theorem}[section]
\newtheorem{lemma}[theorem]{Lemma}
\newtheorem{definition}[theorem]{Definition}
\newtheorem{proposition}[theorem]{Proposition}
\newtheorem{corollary}[theorem]{Corollary}
\newtheorem{example}[theorem]{Example}
\newtheorem{remark}[theorem]{Remark}

\newtheorem{notation}[theorem]{Notation}

\theoremstyle{definition}

\newcommand\pf{\begin{proof}}
\newcommand\epf{\end{proof}}

\newcommand\B{\mathcal{B}}

\newcommand\coc{\mathcal C}
\newcommand\yd{\mathcal{YD}}

\DeclareMathOperator{\Hom}{Hom}

\DeclareMathOperator{\GL}{GL}
\DeclareMathOperator{\SL}{SL}

\numberwithin{equation}{section}

\title{Hopf-Galois objects and cogroupoids}

\author{Julien Bichon}
\address{
Laboratoire de Math\'ematiques,
Universit\'e Blaise Pascal,
Complexe universitaire des C\'ezeaux,
63177~Aubi\`ere Cedex, France}
\email{Julien.Bichon@math.univ-bpclermont.fr}


\begin{document}

\maketitle

\begin{abstract}
 We survey some aspects of the theory of Hopf-Galois objects
that may studied advantageously by using the language of cogroupoids.
These are the notes for a series of lectures given at Cordob\'a University, may 2010.
The lectures are part of the course ``Hopf-Galois theory'' by Sonia Natale.
\end{abstract}

\tableofcontents

\section*{Introduction}

These are the notes for a series of lectures given at Cordob\'a University, may 2010.
The lectures are part of the course ``Hopf-Galois theory'' by Sonia Natale.
We assume therefore some knowledge on Hopf-Galois theory, although some basic facts will be recalled for other readers.

We shall focus on the following two concrete problems. Let $H$ be a Hopf algebra (over a fixed
base field $k$).

\begin{enumerate}
 \item Given an $H$-comodule algebra $A$ that we suspect to be an $H$-Galois object, how can we prove ``nicely'' that $A$ is an $H$-Galois object?
 \item How to classify the $H$-Galois objects?  
\end{enumerate}

Of course these two questions are strongly linked to each other. 
The answers (or tentative answers) we propose rely on Ulbrich's \cite{u89}
equivalence of categories
between the categories of $H$-Galois objects and fibre functors over the category
of $H$-comodules.

To answer the first question, we propose here to use the language of cogroupoids.
This forces us to introduce more notations and concepts, so 
what we have to do first is to explain our motivation. 
Recall that an $H$-Galois object is an $H$-comodule algebra for which a certain
``canonical'' linear map is bijective. The definition is clear and concise, so why should we need
to make it more complicated? Of course the basic answer, which probably
provides enough motivation, is
that we want to be able to check that the canonical map 
is indeed bijective. Another more conceptual motivation comes
with the following parallel situation: given a monoid $G$, then $G$
is a group if and only if the following ``canonical map''
\begin{align*}
 G \times G &\longrightarrow G \times G\\
(x,y) &\longmapsto (xy,y)
\end{align*}
 is bijective. Of course this gives a short definition of groups that
to not use the axioms of inverses, but for many obvious reasons
we prefer to use the (slightly more involved) axioms of inverses.
This is exactly the same philosophy that leads to the use of cogroupoids
in Hopf-Galois theory: we will have (much) more axioms but on the
other hand they should be more natural and easier to deal with.
We present various examples of cogroupoids (Section 3). We hope that
they will convince the reader that it is not more difficult (and in some sense easier)
to work with cogroupoids rather than with Galois objects.

One of the main motivations for studying Hopf-(bi)Galois objects is an important
result by Schauenburg \cite{sc1} stating that the comodule categories
over two Hopf algebras $H$ and $L$ are monoidally equivalent
if and only if there exists an $H$-$L$-bi-Galois object. 
The knowledge of the full
cogroupoid structure (rather than ``only'' the bi-Galois object)
might be useful to exactly determine the image of an object
by the monoidal equivalence. It is also the aim
of the notes to present several applications of this:
construction of new explicit resolutions from old ones in homological algebra, invariant theory,
monoidal equivalences between categories of Yetter-Drinfeld 
modules with applications to bialgebra cohomology and Brauer groups.
The use of the cogroupoid structure is probably not necessary everywhere,
but we believe that it can help!

The second question that we are concerned with is the classification problem
of the Hopf-Galois objects of a given Hopf algebra. Once again this is motivated
by monoidal categories of comodules and Ulbrich's theorem \cite{u89}
already mentioned.
We shall see that in fact Ulbrich's theorem is, in some situations,
a very convenient tool to classify Hopf-Galois objects.
\smallskip

These notes are organized as follows.
In Section 1 we recall some basic definitions
and some important results on Hopf-Galois objects.
In most cases we  do not give proofs but at  
 some occasions we give parts of the proofs, when these are useful
for the rest of the paper. In Section 2 we introduce cogroupoids and prove some basic results. We show that a connected cogroupoid 
induces a bi-Galois object and hence a monoidal equivalence between
the comodule categories over two Hopf algebras. Conversely it is shown
that any Hopf-Galois object 
and  any monoidal equivalence between comodule categories always arises from a connected cogroupoid. We also show that a connected cogroupoid  gives rise to a weak Hopf algebra. In Section 3 we present various examples of cogroupoids. 
In Section 4 we use the fibre functor method (and the constructions of Section 3) to classify  the Galois objects over the Hopf algebras of bilinear forms and the universal cosovereign Hopf algebras.
Section 5 is devoted to applications 
of Hopf-Galois objects and monoidal equivalences to comodule algebras: we describe a model comodule algebra for the Hopf algebra of a bilinear form and
we show how to use Hopf-Galois objects to get new results from old ones in invariant theory. 
In Section 6 we show how monoidal equivalences between comodule categories
extend to categories of Yetter-Drinfeld modules and we give an application
 to Brauer groups of Hopf algebras.
Section 7 is devoted to applications in homological algebra: we show that the Hochshild (co)homology of a Hopf-Galois object is determined by the Hochshild (co)homology of the corresponding Hopf algebras, we show how to transport equivariant resolutions, and we give an application of the result on Yetter-Drindeld modules to bialgebra cohomology.

\smallskip

These notes are mostly a compilation or reformulation of well-known results.
One of the only  new results is the construction (Proposition \ref{coweak})
of a weak Hopf algebra from a cogroupoid (hence from a Hopf-Galois object).
This (unpublished) result was obtained in collaboration with Grunspan in 2003-2004.
A similar result has been obtained independently by De Commer \cite{dec1,dec2}.
The notions of cocategory and cogroupoid were presented by the author
in several talks between 2004 and 2008 (together with the result with Grunspan).
It seems to be their first formal appearance in printed form, but of course
these notions are so natural that they might have been defined
by any mathematician who would have needed them! 

\smallskip

There are many aspects of Hopf-Galois theory (in particular Galois correspondences and Hopf-Galois extensions with non-trivial coinvariants)
that are ignored in these notes.
The reader might consult the excellent survey papers
\cite{monsurv,scsurv} for these topics.

\medskip

\noindent
\textbf{Notations and conventions.}
Throughout these notes we work over a fixed base field denoted $k$.
We assume that the reader has some knowledge of Hopf algebra theory
and Hopf-Galois theory 
(as in \cite{mon}, for Hopf-Galois theory we recall all the necessary
definitions) and on monoidal category theory \cite{kas}.
We use the standard notations, and in particular Sweedler's notation $\Delta(x) = x_{(1)}\otimes x_{(2)}$. The $k$-linear monoidal category
of right $H$-comodules (resp. finite-dimensional right $H$-comodules)
over a Hopf algebra $H$ is denoted ${\rm Comod}(H)$
(resp. ${\rm Comod}_f(H)$), and the set of $H$-comodule morphisms
($H$-colinear maps)
between  $H$-comodules $V$, $W$ is denoted $\Hom_H(V,W)$.

\medskip

\noindent
\textbf{Acknowledgements.} It is a pleasure to thank Sonia Natale for the invitation to give
these lectures and for her kind hospitality at Cordob\'a. The visit at Cordob\'a University 
was supported by the program Premer-Prefalc and Conicet (PIP CONICET  112-200801-00566t).

\section{Background on Hopf-Galois objects}

In this section we collect the definitions, constuctions and results needed in the paper.

\subsection{Hopf-Galois objects}

\begin{definition}
 Let $H$ be a Hopf algebra. A \textbf{left $H$-Galois object} is a left $H$-comodule algebra 
 $A \not = (0)$ such that
 if $\alpha : A \longrightarrow H \otimes A$ denotes the coaction of $H$ on $A$, the linear map
 \begin{equation*}
\begin{CD}
\kappa_l : A \otimes A @>\alpha \otimes 1_A>>
H \otimes A \otimes A @>1_H \otimes m>> H \otimes A
\end{CD}
\end{equation*}
is an isomorphism.  
A \textbf{right $H$-Galois object} 
is a right $H$-comodule algebra $A\not=(0)$ such that if $\beta : A \longrightarrow A \otimes H$ denotes the coaction of $H$ on $A$, the linear map
$\kappa_r$ defined by the composition
\begin{equation*}
\begin{CD}
\kappa_r : A \otimes A @>1_A \otimes \beta>>
A \otimes A \otimes H @>m \otimes 1_H>> A \otimes H
\end{CD}
\end{equation*}
is an isomorphism.

If $L$ is another Hopf algebra, an \textbf{$H$-$L$-bi-Galois object} is an $H$-$L$-bicomodule algebra
which is both a left $H$-Galois object and a right $L$-Galois object.
\end{definition}

The maps $\kappa_l$ and $\kappa_r$ are often called the \textbf{canonical maps}.

\begin{example}{\rm 
 The Hopf algebra $H$, endowed with its comultiplication $\Delta$, is itself an $H$-$H$-bi-Galois object.

Note that if $H$ is a bialgebra, the maps $\kappa_l$ and $\kappa_l$
are well defined and  $H$ is a Hopf algebra if and only if
$\kappa_l$ is bijective if and only if $\kappa_r$ is bijective. 
This gives a definition of Hopf algebras without using the antipode axioms.
However, the definition with the antipode is much more preferable.
}
\end{example}

\begin{example}
 {\rm 
Let $H$ be a Hopf algebra. Recall (see e.g. \cite{doi})
that a \textbf{2-cocycle} on $H$ is a convolution invertible linear map
$\sigma : H \otimes H \longrightarrow k$ satisfying
$$\sigma(x_{(1)}, y_{(1)}) \sigma(x_{(2)}y_{(2)},z) =
\sigma(y_{(1)},z_{(1)}) \sigma(x,y_{(2)} z_{(2)})$$
and $\sigma(x,1) = \sigma(1,x) = \varepsilon(x)$, for all $x,y,z \in H$.
The set of 2-cocycles on $H$ is denoted $Z^2(H)$. When $H =k[G]$ is a group
algebra, it is easy to check that we have an identification $Z^2(k[G]) \simeq
Z^2(G,k^*)$. Note however that in general there is no natural group structure
on $Z^2(H)$.

The convolution inverse of $\sigma$, denoted $\sigma^{-1}$, satisfies 
$$\sigma^{-1}(x_{(1)}y_{(1)},z) \sigma^{-1}(x_{(2)},y_{(2)}) =
\sigma^{-1}(x,y_{(1)}z_{(1)}) \sigma^{-1}(y_{(2)}, z_{(2)})$$
and 
$\sigma^{-1}(x,1) = \sigma^{-1}(1,x) = \varepsilon(x)$, 
for all $x,y,z \in H$.

The algebra 
$_{\sigma} \! H$ is defined as follows. 
As a vector space $_{\sigma} \! H = H$ and the product
of $_{\sigma}H$ is defined to be
$$x \cdot y = \sigma(x_{(1)}, y_{(1)}) x_{(2)} y_{(2)}, 
\quad x,y \in H.$$
That $_{\sigma} \! H$ is an associative algebra with $1$ as unit
follows from the $2$-cocycle condition. Moreover 
$_{\sigma}H$ is a right $H$-comodule algebra with $\Delta : {_{\sigma}\!H} \longrightarrow {_{\sigma}\!H} \otimes H$ as a coaction, and is a right $H$-Galois object. 

Similarly we  have the algebra $H_{\sigma^{-1}}$. As a vector space 
$H_{\sigma^{-1}} = H$ and the product of 
$H_{\sigma^{-1}}$ is defined to be
$$x \cdot y= \sigma^{-1}(x_{(2)}, y_{(2)}) x_{(1)} y_{(1)}, 
\quad x,y \in H.$$
Then $H_{\sigma^{-1}}$ is a left $H$-comodule algebra with coaction
$\Delta : H_{\sigma^{-1}} \longrightarrow H \otimes H_{\sigma^{-1}}$ and
$H_{\sigma^{-1}}$ is a left $H$-Galois object.
We will (re)prove these facts in Subsection 3.3.
}
\end{example}

\begin{example}\label{opp}{\rm
 Let $H$ be a Hopf algebra and let $A$ be a left $H$-Galois object
with coaction $\alpha : A \longrightarrow H \otimes A$, $\alpha(a) =a_{(-1)} \otimes a_{(0)}$. Then the linear map $\beta : A^{\rm op} \longrightarrow A^{\rm op} \otimes H$, $\beta(a) = a_{(0)} \otimes S(a_{(-1)})$ endows $A^{\rm op}$ with a right
$H$-comodule algebra structure, and $A^{\rm op}$ is right $H$-Galois if 
the antipode $S$ is bijective (exercise).

Hence if the antipode of $H$ is bijective, there is no essential difference
between left anf right $H$-Galois objects. 
}
\end{example}

\begin{definition}
Let $H$ be a Hopf algebra. The \textbf{category of left $H$-Galois objects} (resp. \textbf{right $H$-Galois objects}), denoted $\underline{{\rm Gal}}^l(H)$
(resp.  $\underline{{\rm Gal}}^r(H)$), is the category whose objects
are  left $H$-Galois objects (resp. right $H$-Galois objects) and whose morphisms
are $H$-colinear algebra maps (i.e. $H$-comodule algebra maps).

The set of isomorphism classes of left $H$-Galois objects (resp. right $H$-Galois objects) is denoted ${\rm Gal}^l(H)$
(resp.  ${\rm Gal}^r(H)$).
\end{definition}

It follows from Example \ref{opp} that if $H$ has bijective antipode, the 
categories $\underline{{\rm Gal}}^l(H)$ and  $\underline{{\rm Gal}}^r(H)$
are isomorphic, and in this case we simply put ${\rm Gal}(H)=
{\rm Gal}^l(H)={\rm Gal}^r(H)$.

\medskip

The following result (Remark 3.11 in \cite{schnei}), 
means that the categories $\underline{{\rm Gal}}^l(H)$ and  $\underline{{\rm Gal}}^r(H)$ are groupoids (i.e. every morphism is an isomorphism).
It is very useful for classification results.

\begin{proposition}\label{galgroupo}
 Let $H$ be a Hopf algebra and let $A$, $B$ some (left or right) $H$-Galois objects.
Any $H$-colinear algebra map $f : A  \longrightarrow B$ is an isomorphism.
\end{proposition}

\begin{proof}
 Assume for example that $A$ and $B$ are left $H$-Galois and denote
 by $\kappa_l^A$ and $\kappa_l^B$ the respective canonical maps.
 Endow $B$ with the left $A$-module structure induced by $f$.
 The diagram
  \begin{equation*}
\begin{CD}
A \otimes B  @>\kappa_l^B\circ (f \otimes 1_B)>>  H \otimes B\\
@V \simeq VV 
@V\simeq VV\\
(A \otimes A)\otimes_A B  @>\kappa_l^A \otimes_A1_B>> 
(H\otimes A)  \otimes_A B
\end{CD}
\end{equation*}
commutes and hence $f \otimes 1_B$ is an isomorphism, and so is $f$.
\end{proof}

We get a simple criterion to test if a Galois object is trivial.

\begin{proposition}
 Let $H$ be a Hopf algebra and let $A$ be a left or right $H$-Galois object.
 Then $A \cong H$ as $H$-comodule algebras if and only if there exists an algebra 
 map $\phi : A \longrightarrow k$.
\end{proposition}

\begin{proof}
 It is clear, using the counit, that if $A \cong H$ as comodule algebras, then there exists
 an algebra map $A \longrightarrow k$. Conversely, assume for example that
 $A$ is left $H$-galois and that there exists an algebra map $\phi : A \longrightarrow k$.
 Then the map $A \longrightarrow H$, $a \longmapsto \phi(a_{(0)}) a_{(-1)}$
 is a left $H$-colinear algebra map, and is an isomorphism by the previous proposition. 
\end{proof}

\subsection{Cleft Hopf-Galois objects} 

In this subsection we briefly focus on an important class of Hopf-Galois objects, called \textbf{cleft}.

We have already seen (Example 1.3) how to associate a Galois object to a 2-cocyle.
This construction is axiomatized by the
following result, which summarizes work of Doi-Takeuchi \cite{doitak}
and Blattner-Montgomery \cite{blamon}.
The proof can be found in \cite{mon} or \cite{scsurv}.

\begin{theorem}
 Let $H$ be a Hopf algebra and let $A$ be a left $H$-Galois object.
The following assertions are equivalent:
\begin{enumerate}
 \item There exists $\sigma \in Z^2(H)$ such that $A \cong H_{\sigma^{-1}}$ as left comodule algebras.
\item $A \cong H$ as left $H$-comodules.
\item There exists a convolution invertible $H$-colinear map $\phi : H \longrightarrow A$
\end{enumerate}
A left $H$-Galois object is said to be \textbf{cleft} if it satisfies the above equivalent conditions. 
\end{theorem}

Of course there is a similar result for right $H$-Galois objects.

There are nice classes of Hopf algebras for which cleftness is automatic.
Recall that a Hopf algebra is said to be pointed if all its simple
comodules are $1$-dimensional.

\begin{theorem}\label{finipoin}
 Let $H$ be a finite-dimensional or pointed Hopf algebra.
Then any $H$-Galois object is cleft.
\end{theorem}

 For the pointed case we refer the reader to Remark 10 in \cite{gun} 
and for the finite-dimensional case we refer the reader
to \cite{kc}.
See however the last remark in the next subsection, where we give a proof by using
fibre functors.

\subsection{Monoidal equivalences and Schauenburg's Theorem}

The following result from \cite{sc1} is one of the main motivations for the study of (bi-)Galois objects.

\begin{theorem}[Schauenburg]\label{schau}
 Let $H$ and $L$ be some Hopf algebras. The following assertions are equivalent.
 \begin{enumerate}
  \item There exists a $k$-linear equivalence of monoidal categories
$${\rm Comod}(H) \cong^\otimes {\rm Comod}(L)$$
\item There exists an $H$-$L$-bi-Galois object.
 \end{enumerate}
\end{theorem}

So we have a very powerful tool to construct monoidal equivalences
between categories of comodules over Hopf algebras!
In the next subsection we present various very important
constructions related to Schauenburg's Theorem.

\subsection{Cotensor product and various constructions of functors}

In this subsection we recall several important constructions of functors
associated with Hopf-Galois objects.
The basic construction is the cotensor product.

\begin{definition}
Let $C$ be a coalgebra, let $V$ be a right $C$-comodule and let 
$W$ be a left $C$-comodule. The \textbf{cotensor product of $V$ and $W$}, denoted
$V\square_CW$, 
is defined to be the equalizer
$$0 \longrightarrow V \square_CW \longrightarrow V \otimes W \rightrightarrows V \otimes C \otimes W$$
i.e. the kernel of the map $\alpha_V \otimes 1_W - 1_V\otimes \alpha_W$, where
$\alpha_V$ and $\alpha_W$ are the respective coactions on $V$ and $W$.
\end{definition}

The first thing that the cotensor product allows us to do 
is to construct functors on comodule categories, as follows.

\begin{proposition}\label{basi}
 Let $C$ be a coalgebra. Any left $C$-comodule $W$
defines  a $k$-linear functor
\begin{align*}
\Omega^W : {\rm Comod}(C) & \longrightarrow {\rm Vect}(k) \\
 V & \longmapsto V \square_C W
\end{align*}
\end{proposition}

\begin{proof}
 If $f : V \longrightarrow V'$ is $C$-colinear map,
then it is easy to check that $f \otimes 1(V \square_C W) \subset V' \square_C W$
and hence we get the announced functor.
\end{proof}

\begin{remark}\label{trifun}
If $W \cong C$ as left $C$-comodules, then $\Omega^W$ is isomorphic
to the forgetful functor.
\end{remark}

\begin{proof}
Let  $f : C \longrightarrow W$ be a $C$-colinear isomorphism. If $V$
is a right $C$-comodule, we have a vector space isomorphism
\begin{align*}
 V &\longrightarrow V\square_C W =\Omega^W(V)\\
v &\longmapsto v_{(0)} \otimes f(v_{(1)})
\end{align*}
 (the inverse is given by $v \otimes a \mapsto \varepsilon(f^{-1}(a))v$).
The isomorphism is clearly functorial and we have the result. 
\end{proof}

The monoidal analogue of  Proposition \ref{basi}
is Ulbrich's theorem \cite{u89}:

\begin{theorem}[Ulbrich]\label{ulbrich}
Let $H$ be a Hopf algebra and let $A$ be a left $H$-Galois object.
The functor 
\begin{align*}
 \Omega^A : {\rm Comod}_f(H) & \longrightarrow {\rm Vect}_f(k) \\
 V & \longmapsto V \square_H A
\end{align*}
is a fibre functor: it is $k$-linear, exact, faithful and monoidal.
Conversely, any fibre functor arises in this way from a unique (up to isomorphism)
left $H$-Galois object. 
\end{theorem}

\begin{proof}[Partial proof]
 Consider the previous functor $\Omega^A : {\rm Comod}(H) \longrightarrow {\rm Vect}(k)$. We endow $\Omega^A$ with a monoidal structure as follows (the construction is from \cite{u87}).
First define $\widetilde\Omega_0^A : k \longrightarrow \Omega^A(k) = k \square_H A = A^{{\rm co} H}$, $1 \longmapsto 1$. This is an isomorphism since for
$a \in A^{{\rm co} H}$, we have $\kappa_l(a\otimes 1) =1\otimes a = \kappa_l(1\otimes a)$, hence   we have $a\otimes 1=1 \otimes a$ by the injectivity of $\kappa_l$
and $a \in k1$. 

Now let $V,W \in {\rm Comod}(H)$.
It is straightforward to check
that if $\sum_i v_i \otimes a_i \in V \square_{H}A$
and $\sum_j w_j \otimes b_j \in W \square_{H} A$, then
$$\sum_{i,j} v_i \otimes w_j\otimes a_ib_j \in (V \otimes W) \square_{H} A$$
 Thus we have a map
\begin{align*}
 (V \square_{H} A) \otimes
(W \square_{H} A) 
&\longrightarrow (V\otimes W) \square_{H} A \\
(\sum_i v_i \otimes a_i) \otimes
(\sum_j w_j \otimes b_j) &\longmapsto
\sum_{i,j} v_i \otimes w_j\otimes a_ib_j
\end{align*}
that we denote
$$ \widetilde{\Omega}^A_{V,W} : \Omega^A(V) \otimes \Omega^A(W) \longrightarrow \Omega^A(V \otimes W)$$
It is clear that this map is functorial, and we have to check
that $\Omega^A = (\Omega^A,  \widetilde{\Omega}^A_{\bullet,\bullet}, \widetilde \Omega_0)$ is a monoidal functor. It is easy to seen that the associativity (coherence) constraints
of a monoidal functor are satisfied, and what is really non trivial
is to check that $\widetilde{\Omega}^A_{V,W}$ is an isomorphism for
any $V,W \in {\rm Comod}(H)$. We use the following Lemma.

\begin{lemma}
\
 \begin{enumerate}
\item For any $V \in {\rm Comod}(H)$, the linear map
$1 \otimes \varepsilon \otimes 1:(V \otimes H) \square_H A \longrightarrow V \otimes A$
is an isomorphism.
\item For any $V \in {\rm Comod}(H)$, the linear map
$m_V:=1 \otimes m :(V \square_H A)\otimes A \longrightarrow V \otimes A$
is an isomorphism.
 \end{enumerate}
\end{lemma}
 \begin{proof}[Proof of the  Lemma]
(1) The reader will easily check that the inverse $V \otimes A \longrightarrow (V \otimes H)\square_H A$ is
given by $v \otimes a \longmapsto v_{(0)} \otimes S(v_{(1)}) a_{(-1)} \otimes a_{(0)}$.


(2) It is not difficult to check that $1 \otimes \kappa_l : V \otimes A \otimes A \longrightarrow V \otimes H \otimes A$ 
induces a map $1 \otimes \kappa_l : (V \square_H A) \otimes A \longrightarrow (V \square_H H) \otimes A$ and that the following diagram commutes
$$\xymatrix{
(V \square_H A) \otimes A  \ar[dr]_{m_V}  \ar[rr]^{1 \otimes \kappa_l} &  & (V \square_H H) \otimes A \ar[dl]_{1 \otimes \varepsilon \otimes 1}\\
& V \otimes A &
}$$
with the vertical arrow on the right bijective by (1), so $m_V$ is an isomorphism. 
\end{proof}

We are now ready to prove that $\widetilde{\Omega}^A_{V,W}$ is an isomorphism.
The following diagram commutes
{\footnotesize
$$\xymatrix{
  (V \square_{H} A) \otimes
(W \square_{H} A) \otimes A \ar[d]^{1 \otimes m_W} \ar[rrr]^{\widetilde{\Omega}_{V,W}^A\otimes 1} &
 & & ((V\otimes W) \square_{H} A) \otimes A \ar[d]^{m_{V \otimes W}} \\
(V \square_{H} A) \otimes W \otimes A \ar[r]^{1 \otimes \tau} & (V \square_{H} A) \otimes A \otimes W
\ar[r]^-{m_V \otimes 1} & V \otimes A \otimes W \ar[r]^{1 \otimes \tau} & V \otimes W \otimes A
}$$}
where the $\tau$'s denote the canonical flips $x \otimes y \mapsto y \otimes x$.
Since by the Lemma all the vertical and horizontal down morphisms are isomorphisms
we conclude that $\widetilde{\Omega}_{V,W}^A\otimes 1$ is an isomorphism,
and so is 
$\widetilde{\Omega}_{V,W}^A$: we have our monoidal functor
 $\Omega^A =(\Omega^A, \widetilde \Omega_{\bullet, \bullet}^A, \widetilde \Omega_0^A): {\rm Comod}(H) \longrightarrow {\rm Vect}(k)$.
That $\Omega^A : {\rm Comod}(H) \longrightarrow {\rm Vect}(k)$ restricts
to a monoidal functor $\Omega^A :{\rm Comod}_f(H) \longrightarrow {\rm Vect}_f(k)$
is a consequence of the forthcoming Proposition \ref{fini}.

Faithfulness of $\Omega^A$ is obvious while exactness is easy using
(3) in the Lemma (see \cite{u87}), hence $\Omega^A$ is a fibre functor. We do not prove the converse (see \cite{u87,u89}).
\end{proof}

The following result has been used in the partial proof of Ulbrich's theorem
and will also be useful elsewhere.

\begin{proposition}\label{fini}
 Let $H$ be a Hopf algebra and let $F : {\rm Comod}(H) \longrightarrow {\rm Vect}(k)$
be a monoidal functor. If $V$ is a finite-dimensional $H$-comodule, then $F(V)$
is a finite-dimensional vector space. Moreover we have $\dim (V)=1 \Rightarrow\dim(F(V))=1$, and if $F$ is a fibre functor then $\dim(F(V))=1 \Rightarrow \dim( V)=1$.
\end{proposition}

\begin{proof}
The key point is that a vector space $V$ is finite-dimensional
if and only if there exists a vector space $W$ and linear maps
$e : W \otimes V \rightarrow k$ and $d : k \rightarrow V \otimes W$ such that
$$1_V = (1_V \otimes e) \circ (d\otimes 1_V), \ 
1_W = (e \otimes 1_W) \circ (1_W\otimes d)$$
If $V \in {\rm Comod}_f(H)$, the dual comodule $V^*$
satisfies to the above requirement with $e$ and $d$ $H$-colinear, where
$e$ is the evaluation map: in other words $V$ has a left dual in ${\rm Comod}(H)$,
see \cite{js} or \cite{kas}. Applying the monoidal functor
$F$, we easily see that $F(V)$ satisfies to the above condition and hence
that $F(V)$ is finite-dimensional. 
If $\dim(V)=1$, then $e: V^* \otimes V \longrightarrow k$ is an isomorphism and
so is the induced composition $F(V^*) \otimes F(V) \cong F(V^*\otimes V) \cong F(k) \cong k$, hence $\dim(F(V))=1$.  If $\dim (F(V))=1$, then the previous composition is an isomorphism (because the monoidal functor transforms left duals into left duals),
and hence $F(e)$ is an isomorphism. If $F$ is a fibre functor then it is exact faithful so $e$ is an isomorphism, which shows that $\dim (V)=1$. 
\end{proof}

Of course there is a left-right version of Ulbrich's Theorem.
Note that Ulbrich's theorem in \cite{u89} is a stronger statement:
it states an equivalence of categories between left $H$-Galois objects
and fibre functors on ${\rm Comod}_f(H)$.

In the  setting of fibre functors there is a nice characterization of cleftness in terms of fibre functors, essentially due to Etingof-Gelaki \cite{eg}.

\begin{theorem}\label{cleftfibre}
 Let $H$ be a Hopf algebra and let $A$ be a left $H$-Galois object.
The following assertions are equivalent.
\begin{enumerate}
 \item $A$ is a cleft $H$-Galois object.
\item The fibre functor $\Omega^A : {\rm Comod}_f(H) \longrightarrow {\rm Vect}_f(k)$ preserves the dimensions of the underlying vector spaces. 
\end{enumerate}
\end{theorem}

\begin{proof}[Partial proof]
 We give the proof of $(1) \Rightarrow (2)$ (the easy part) and we refer
to \cite{eg} for the proof of $(2) \Rightarrow (1)$. If $A$ is cleft
there exists an $H$-colinear isomorphism $H \cong A$ and by 
Remark \ref{trifun}
  $\Omega^A$ is isomorphic
to the forgetful functor, which proves the result.
\end{proof}

\begin{remark}{\rm
 We can use the fibre functor interpretation to give a proof of
Theorem \ref{finipoin}. Let $A$ be a left $H$-Galois object.

Assume first that $H$ is pointed. Any simple $H$-comodule is one dimensional, so  $F_A$ preserves the dimension of any simple
$H$-comodule by Proposition \ref{fini}. An induction now shows that $F_A$ preserves the dimension of any
finite-dimensional $H$-comodule, and hence $A$ is cleft by Theorem \ref{cleftfibre}. 

Assume now that $H$ is finite-dimensional. We have a linear isomorphism $A \cong H\square_H A=\Omega^A(H)$, $a \mapsto a_{(-1)}\otimes a_{(0)}$,
so $A$ is finite-dimensional since $\Omega^A$ has its values in ${\rm Vect}_f(k)$.
Denote by $A_0$ the trivial $H$-comodule whose underlying vector space is $A$.
The canonical map
$\kappa_l : A \otimes A_0 \longrightarrow H \otimes A_0$ is a left $H$-comodule
isomorphism. Hence we have a left $H$-comodule isomorphism
$A^{\dim(A)} \cong H^{\dim(A)}$ and so the Krull-Schmidt Theorem shows
that $A\cong H$ as left $H$-comodules.
}
\end{remark}

We now come to the construction of functors between categories of comodules.
The basic result is the following one.

\begin{proposition}\label{basi2}
 Let $C$ and $D$ be some coalgebras.
Any $C$-$D$-bicomodule $X$ defines a $k$-linear functor
\begin{align*}
F^X : {\rm Comod}(C) & \longrightarrow {\rm Comod}(D) \\
 V & \longmapsto V \square_C X
\end{align*}
\end{proposition}

\begin{proof}
It is clear that $V \otimes X$ is a right $D$-comodule for the coaction
$1_V \otimes \beta_X$, where $\beta_X : X \longrightarrow X \otimes D$
is the right $D$-coaction on $X$. We have to check that
$V \square_C X$ is a $D$-subcomodule of $V \otimes X$.
So let $\sum_i v_i \otimes x_i \in V \square_C X$. We have
$$\sum_i v_{i(0)} \otimes v_{i(1)} \otimes x_i
= \sum_i v_{i} \otimes x_{i(-1)} \otimes x_{i(0)}$$
and hence
$$\sum_i v_{i(0)} \otimes v_{i(1)} \otimes x_{i(0)} \otimes x_{i(1)}
= \sum_i v_{i} \otimes x_{i(-1)} \otimes x_{i(0)}\otimes x_{i(1)}$$
which shows that $\sum_i v_{i}  \otimes x_{i(0)} \otimes x_{i(1)} \in
(V\square_C X) \otimes D$, and hence $V \square_C X$ is a $D$-subcomodule of $V \otimes X$. It is clear that if $f : V \longrightarrow W$ is a $C$-comodule
map, then $f\otimes 1_X$ induces a $D$-comodule map $V \square_C X\longrightarrow W \square_C X$, and we have our functor.
\end{proof}

Takeuchi \cite{tak} has given the precise conditions for which
the above functor is an equivalence: the bicomodule has to be part
of some more structured data, known now as a Morita-Takeuchi context, 
and he has shown that
any $k$-linear equivalence between categories of comodules arises in that way.
The axioms of cocategories discussed in the next section are quite
close from those of Morita-Takeuchi contexts. 

\medskip

We finish the subsection by  a monoidal version of Proposition \ref{basi2}. 

\begin{proposition}\label{basimonoidal}
 Let $H$ and $L$ be some Hopf algebras. Let $A$ be an $H$-$L$-bicomodule
algebra such that $A$ is left $H$-Galois.
Then $A$ defines a $k$-linear monoidal functor
\begin{align*}
F^A : {\rm Comod}(H) & \longrightarrow {\rm Comod}(L) \\
 V & \longmapsto V \square_H A
\end{align*}
\end{proposition}

\begin{proof}
The functor is provided by Proposition \ref{basi2}.
By the (partial) proof of Ulbrich's Theorem we have
for any $V,W \in {\rm Comod}(H)$ linear isomorphisms
\begin{align*}
 (V \square_{H} A) \otimes
(W \square_{H} A) 
&\longrightarrow (V\otimes W) \square_{H} A \\
(\sum_i v_i \otimes a_i) \otimes
(\sum_j w_j \otimes b_j) &\longmapsto
\sum_{i,j} v_i \otimes w_j\otimes a_ib_j
\end{align*}
which are easily seen to be $L$-colinear.
In this way $F^A$ is a monoidal functor, as announced.
\end{proof}

The proof of the $(2)\Rightarrow (1) $ part of Schauenburg's theorem
(from a bi-Galois object to a monoidal equivalence) uses the construction of 
Proposition \ref{basimonoidal}.

\section{Cocategories and cogroupoids}

In this section we  put Hopf bi-Galois objects into a more structured framework.
The idea is that although the axiomatic becomes more complicated at first sight, it should 
be more natural and easier to deal with.

\subsection{Basic definitions}
The first step is to propose a notion that is dual to the one of category, as follows. 

\begin{definition}
A \textbf{cocategory} (or $k$-cocategory) $\coc$ consists of:

\noindent
$\bullet$ a set of objects ${\rm ob}(\coc)$.

\noindent
$\bullet$ For any $X,Y \in {\rm ob}(\coc)$, a $k$-algebra 
$\coc(X,Y)$. 

\noindent
$\bullet$ For any $X,Y,Z \in {\rm ob}(\coc)$, algebra morphisms
$$\Delta_{X,Y}^Z : \coc(X,Y) \longrightarrow \coc(X,Z) \otimes \coc(Z,Y)
\quad {\rm and} \quad \varepsilon_X : \coc(X,X) \longrightarrow k$$
such that for any $X,Y,Z,T \in {\rm ob}(\coc)$,
the following diagrams commute:
\begin{equation*}
\begin{CD}
\coc(X,Y) @>\Delta_{X,Y}^Z>> \coc(X,Z) \otimes \coc(Z,Y) \\
@V\Delta_{X,Y}^TVV
@V\Delta_{X,Z}^T \otimes 1VV\\
\coc(X,T) \otimes \coc(T,Y)  @>1 \otimes \Delta_{T,Y}^Z>> 
\coc(X,T)\otimes \coc(T,Z)  \otimes \coc(Z,Y)
\end{CD}
\end{equation*}
$$\xymatrix{
\coc(X,Y) \ar[d]^{\Delta_{X,Y}^Y} \ar@{=}[rrd] \\
 \coc(X,Y) \otimes \coc(Y,Y) \ar[rr]^{1 \otimes \varepsilon_Y} && \coc(X,Y)}
\quad
\xymatrix{
\coc(X,Y) \ar[d]^{\Delta_{X,Y}^X} \ar@{=}[rrd] \\
 \coc(X,X) \otimes \coc(X,Y) \ar[rr]^{\varepsilon_X \otimes 1} && \coc(X,Y)}
$$
\end{definition}

The following results are immediate consequences of the axioms.

\begin{proposition}\label{consequence}
 Let $\coc$ be a cocategory and let $X,Y, Z \in {\rm ob}(\coc)$.
 \begin{enumerate}
  \item $\coc(X,X)= (\coc(X,X), \Delta_{X,X}^X, \varepsilon_X)$ is a bialgebra.
  \item $\coc(X,Y)$ is a $\coc(X,X)$-$\coc(Y,Y)$-bicomodule algebra, via $\Delta_{X,Y}^X$ and 
  $\Delta^{Y}_{X,Y}$ respectively.
\item $\Delta_{X,Y}^Z : \coc(X,Y) \longrightarrow \coc(X,Z) \otimes \coc(Z,Y)$
is a $\coc(X,X)$-$\coc(Y,Y)$-bicomodule algebra
morphism and we have $\Delta_{X,Y}^Z(\coc(X,Y)) \subset \coc(X,Z) \square_{\coc(Z,Z)} \coc(Z,Y)$. 
 \end{enumerate}
\end{proposition}

Thus a cocategory with one object is just a bialgebra.

\begin{definition}
A cocategory $\coc$ is said to be \textbf{connected} if
$\coc(X,Y)$ is a non zero algebra  
for any $X,Y \in {\rm ob}(\coc)$.
\end{definition}

A groupoid is a category whose morphisms all are isomorphisms.
A dual notion is the following one. 

\begin{definition}
A \textbf{cogroupoid} (or a $k$-cogroupoid) $\coc$ consists of a cocategory $\coc$ together
with, for any $X,Y \in {\rm ob}(\coc)$, linear maps
$$S_{X,Y} : \coc(X,Y) \longrightarrow \coc(Y,X)$$
such that for any $X, Y \in {\rm ob}(\coc)$, the following diagrams commute: 
{\footnotesize
$$\xymatrix{
\coc(X,X) \ar[r]^{\varepsilon_X} \ar[d]^{\Delta_{X,X}^Y} &
k \ar[r]^{u}  & \coc(X,Y) \\
\coc(X,Y) \otimes \coc(Y,X) \ar[rr]^{1 \otimes S_{Y,X}} && \coc(X,Y) \otimes \coc(X,Y) \ar[u]^{m}} \quad \quad
\xymatrix{
\coc(X,X) \ar[r]^{\varepsilon_X} \ar[d]^{\Delta_{X,X}^Y} &
k \ar[r]^{u}  & \coc(Y,X) \\
\coc(X,Y) \otimes \coc(Y,X) \ar[rr]^{S_{X,Y} \otimes 1} && \coc(Y,X) \otimes \coc(Y,X) \ar[u]^{m}}
$$} where $m$ denotes the multiplication of $\coc(X,Y)$ and $u$ is the unit map.
\end{definition}

\begin{remark} [on terminology]
 {\rm A connected cogroupoid with two objects is exactly what Grunspan
called a total Hopf-Galois system in \cite{gru1} (some axioms are redundant in \cite{gru1}), which was a symmetrisation
of the notion of Hopf-Galois system from  \cite{bi2}.
It seems that the more compact present formulation is much more convenient and pleasant to deal with.  
}
\end{remark}

\begin{definition}
   A \textbf{full subcocategory} (resp. \textbf{full subcogroupoid}) of a 
cocategory (resp. cogroupoid) $\coc$ is a cocategory (resp. cogroupoid)
$\mathcal D$
with ${\rm ob}(\mathcal D) \subset {\rm ob}(\coc)$, with $\mathcal D(X,Y) = \mathcal C(X,Y)$, $\forall X,Y \in {\rm ob}(\mathcal D)$, and whose structural
morphisms are those induced by the ones of $\coc$.
\end{definition}

\begin{notation}
We now introduce Sweedler's notation for cocategories and cogroupoids.
Let $\coc$ be a cocategory. For $a^{X,Y} \in \coc(X,Y)$, we write 
$$\Delta_{X,Y}^Z(a^{X,Y})= a_{(1)}^{X,Z} \otimes a_{(2)}^{Z,Y}$$
The cocategory axioms now read 
$$(\Delta_{X,Z}^T \otimes 1)\circ \Delta_{X,Y}^Z(a^{X,Y}) = 
a_{(1)}^{X,T} \otimes a_{(2)}^{T,Z} \otimes a_{(3)}^{Z,Y}= 
(1 \otimes \Delta_{T,Y}^Z)\circ \Delta_{X,Y}^T(a^{X,Y})
$$
$$\varepsilon_X(a_{(1)}^{X,X}) a_{(2)}^{X,Y} = a^{X,Y} =
\varepsilon_Y(a_{(2)}^{Y,Y})  a_{(1)}^{X,Y}$$
and the additional cogroupoid axioms are
$$S_{X,Y}(a_{(1)}^{X,Y}) a_{(2)}^{Y,X} = \varepsilon_X(a^{X,X})1 =
a_{(1)}^{X,Y} S_{Y,X}(a_{(2)}^{Y,X})
$$

\end{notation}

\subsection{Back to Hopf-galois objects}

\begin{proposition}\label{back}
 Let $\coc$ be a connected cogroupoid. Then for any $X,Y \in {\rm ob}(\coc)$, 
 the algebra $\coc(X,Y)$ is a $\coc(X,X)-\coc(Y,Y)$ bi-Galois object.
\end{proposition}

\begin{proof}
 We give a proof by using Sweedler's notation (a proof with morphisms can be found  in \cite{bi2}).
 Let $\eta_l : \coc(X,X) \otimes \coc(X,Y) \longrightarrow \coc(X,Y) \otimes \coc(X,Y)$
be defined by 
$$\eta_l = (1 \otimes m) \circ (1 \otimes S_{Y,X} \otimes 1)
\circ (\Delta_{X,X}^Y \otimes 1)$$
i.e. $\eta_l(a^{X,X} \otimes b^{X,Y}) = 
a_{(1)}^{X,Y} \otimes S_{Y,X}(a_{(2)}^{Y,X})b^{X,Y}$. 
We show that $\eta_l$ is an inverse for $\kappa_l$. Let $a^{X,Y},b^{X,Y} \in \coc(X,Y)$.
We have
\begin{align*}
\eta_l \circ \kappa_l(a^{X,Y} \otimes b^{X,Y}) & =
\eta_l(a_{(1)}^{X,X} \otimes a_{(2)}^{X,Y}b^{X,Y})  \\
& =a_{(1)}^{X,Y} \otimes S_{Y,X}(a_{(2)}^{Y,X}) a_{(3)}^{X,Y}b^{X,Y}
= a_{(1)}^{X,Y} \otimes \varepsilon_Y(a_{(2)}^{Y,Y})b^{X,Y} =
a^{X,Y} \otimes b^{X,Y}
\end{align*}
Hence $\eta_l \circ \kappa_l$ is the identity map.
Similarly one checks that $\kappa_l \circ \eta_r= {\rm id}$ and hence 
$\coc(X,Y)$ is left $\coc(X,X)$-Galois.
 Similarly,
we define  $\eta_r : \coc(X,Y) \otimes \coc(Y,Y) \longrightarrow \coc(X,Y) \otimes \coc(X,Y)$
by 
$$\eta_r = (m \otimes 1) \circ (1 \otimes S_{Y,X} \otimes 1)
\circ (1 \otimes \Delta_{Y,Y}^X),$$
and check that $\eta_r$ is an inverse for $\kappa_r$. 
Hence $\coc(X,Y)$ is right $\coc(Y,Y)$-Galois.
\end{proof}

\begin{corollary}\label{cogroumonoidal}
 Let $\coc$ be a connected cogroupoid. Then for any $X,Y \in {\rm ob}(\coc)$, 
 there exists a $k$-linear equivalence of monoidal categories
$${\rm Comod}(\coc(X,X)) \cong^\otimes {\rm Comod}(\coc(Y,Y))$$
\end{corollary}

\begin{proof}
 Just combine the previous result with Schauenburg's Theorem. 
\end{proof}

We now state two results that mean that the theory of Hopf-Galois objects
is actually equivalent to the theory of connected cogroupoids.

\begin{theorem}\label{schacog}
  Let $H$ and $L$ be some Hopf algebras. The following assertions are equivalent.
 \begin{enumerate}
  \item There exists a $k$-linear equivalence of monoidal categories
$${\rm Comod}(H) \cong^\otimes {\rm Comod}(L)$$
\item There exists a connected cogroupoid $\coc$ with two objects $X,Y$
such that $H= \coc(X,X)$ and $L=\coc(Y,Y)$.
 \end{enumerate}
\end{theorem}

\begin{theorem}\label{galoiscog}
 Let $H$ be a Hopf algebra and let $A$ be a left $H$-Galois object.
Then there exists a connected cogroupoid $\coc$ with two objects $X,Y$
such that $H= \coc(X,X)$ and $A=\coc(X,Y)$.
\end{theorem}

The proof of $(1) \Rightarrow (2)$ in Theorem \ref{schacog}
is done in \cite{bi2} by 
using Tannaka-Krein reconstruction techniques (\cite{js,sctan}). 
We give a sketch of the proof
in  Subsection 2.4.

The proof of Theorem \ref{galoiscog} is also given in Subsection 2.4.
Another proof of this result is provided by Grunspan \cite{gru2}, using
quantum torsors and results by Schauenburg (see \cite{sc1} and  Subsection 2.8 of \cite{scsurv}).
This proof has the merit to work for Hopf algebras over rings (of course with flatness assumptions).

The proofs of these two results in Subsection 2.4 will never be used in the rest of the paper, so the reader might skip them first, and get back
if he needs to.

The $(2) \Rightarrow (1)$ part in Theorem \ref{schacog}
is the previous corollary. 
We shall give a proof of it without using Schauenburg's theorem in the next subsection.
More precisely, the  
explicit form of the monoidal equivalence constructed by starting from
a connected cogroupoid is as follows. Having this explicit form
in hand is useful in some applications.

\begin{theorem}\label{explicit}
 Let $\coc$ be a connected cogroupoid. Then for any $X,Y \in {\rm ob}(\coc)$
we have  $k$-linear equivalences of monoidal categories that are inverse of each other
\begin{align*}
{\rm Comod}(\coc(X,X)) & \cong^\otimes {\rm Comod}(\coc(Y,Y)) 
\quad & {\rm Comod}(\coc(Y,Y)) & \cong^\otimes {\rm Comod}(\coc(X,X))\\
V &\longmapsto V \square_{\coc(X,X)} \coc(X,Y) \quad 
& V &\longmapsto V \square_{\coc(Y,Y)} \coc(Y,X)
\end{align*}
\end{theorem}

We give a proof in the next subsection.

\subsection{Some basic properties of cogroupoids}
This section is devoted to state and prove some basic properties
of cogroupoids, and to give a proof of Theorem \ref{explicit}.

We begin by examinating some properties of the ``antipodes''.
Part (1) of the following result is proved very indirectly
in \cite{bi2}, with a direct proof given in \cite{hob}.
In view of the next  Proposition \ref{coweak}, theses properties
also are consequences of results on weak Hopf algebras \cite{bonisz}.

\begin{proposition}\label{anti}
 Let $\coc$ be a cogroupoid and let $X, Y \in {\rm ob}(\coc)$.
\begin{enumerate}
 \item $S_{Y,X} : \coc (Y,X) \longrightarrow \coc(X,Y)^{\rm op}$ is an 
algebra morphism.
\item For any $Z \in {\rm ob}(\coc)$ and $a^{Y,X} \in \coc(Y,X)$, we have
$$\Delta_{X,Y}^Z(S_{Y,X}(a^{Y,X})) = S_{Z,X}(a_{(2)}^{Z,X}) \otimes S_{Y,Z}(a_{(1)}^{Y,Z})$$
\end{enumerate}
\end{proposition}

\begin{proof}
 Let $a^{Y,X}, b^{Y,X} \in \coc(Y,X)$. We have
\begin{align*}
 S_{Y,X}(a^{Y,X} b^{Y,X}) &= S_{Y,X}(a_{(1)}^{Y,X}b_{(1)}^{Y,X})\varepsilon_X(a_{(2)}^{X,X})\varepsilon_X( b_{(2)}^{X,X}) \\
& = S_{Y,X}(a_{(1)}^{Y,X}b_{(1)}^{Y,X})a_{(2)}^{X,Y}b_{(2)}^{X,Y}
S_{Y,X}(b_{(3)}^{Y,X})S_{Y,X}(a_{(3)}^{Y,X}) \\
& = \varepsilon_Y(a_{(1)}^{Y,Y}b_{(1)}^{Y,Y})
S_{Y,X}(b_{(2)}^{Y,X})S_{Y,X}(a_{(2)}^{Y,X}) \\
& = S_{Y,X}(b^{Y,X})S_{Y,X}(a^{Y,X})
\end{align*}
and this proves (1). We also have
\begin{align*}
 S_{Z,X}(a_{(2)}^{Z,X}) \otimes S_{Y,Z}(a_{(1)}^{Y,Z})  & = 
S_{Z,X}(a_{(2)}^{Z,X}\varepsilon_X(a_{(3)}^{X,X})) \otimes S_{Y,Z}(a_{(1)}^{Y,Z}) \\
&= S_{Z,X}(a_{(2)}^{Z,X}) \otimes S_{Y,Z}(a_{(1)}^{Y,Z}) \cdot \left(\varepsilon_X(a_{(3)}^{X,X}) 1 \otimes 1\right) \\
&= S_{Z,X}(a_{(2)}^{Z,X}) \otimes S_{Y,Z}(a_{(1)}^{Y,Z}) \cdot \left(\Delta_{X,Y}^Z(a_{(3)}^{X,Y} S_{Y,X}(a_{(4)}^{Y,X}))\right) \\
&= S_{Z,X}(a_{(2)}^{Z,X}) \otimes S_{Y,Z}(a_{(1)}^{Y,Z}) \cdot 
\left(a_{(3)}^{X,Z} \otimes a_{(4)}^{Z,Y}\right) \cdot \Delta_{X,Y}^Z(S_{Y,X}(a_{(5)}^{Y,X})) \\
& = \varepsilon_Z(a_{(2)}^{Z,Z})1 \otimes S_{Y,Z}(a_{(1)}^{Y,Z})a_{(3)}^{Z,Y}
\cdot \Delta_{X,Y}^Z(S_{Y,X}(a_{(4)}^{Y,X})) \\
& = 1 \otimes S_{Y,Z}(a_{(1)}^{Y,Z})a_{(2)}^{Z,Y}
\cdot \Delta_{X,Y}^Z(S_{Y,X}(a_{(3)}^{Y,X})) \\
& =\Delta_{X,Y}^Z(S_{Y,X}(a^{Y,X}))
\end{align*}
and this proves (2).
\end{proof}

The following result is useful to prove connectedness properties of cogroupoids,
and also for the proof of Theorem \ref{explicit} 

\begin{lemma}\label{useful}
 Let $\coc$ be a cogroupoid and let $X,Y,Z \in {\rm ob}(\coc)$.
Assume that $\coc(Z,Y) \not =(0)$ or $\coc(X,Z) \not =(0)$.
Then $\Delta_{X,Y}^Z : \coc(X,Y)\longrightarrow\coc(X,Z) \otimes \coc(Z,Y)$
is split injective, and induces a $\coc(X,X)-\coc(Y,Y)$-bicomodule algebra
isomorphism
$$
\coc(X,Y)\cong\coc(X,Z) \square_{\coc(Z,Z)} \coc(Z,Y)$$
\end{lemma}

\begin{proof}
 Assume first that $\coc(Z,Y) \not = (0)$, and let
$\psi : \coc(Z,Y) \longrightarrow k$ be a linear map such that
$\psi(1)=1$. Define $f : \coc(X,Z) \otimes \coc(Z,Y) \longrightarrow  \coc(X,Y)$
by 
$$f(a^{X,Z}\otimes b^{Z,Y}) = \psi\left(S_{Y,Z}(a^{Y,Z}_{(2)})b^{Z,Y}\right)a^{X,Y}_{(1)}$$
Then 
$$f \circ \Delta_{X,Y}^Z(a^{X,Y})
=f(a_{(1)}^{X,Z}\otimes a_{(2)}^{Z,Y})
= \psi\left(S_{Y,Z}(a^{Y,Z}_{(2)})a_{(3)}^{Z,Y}\right)a^{X,Y}_{(1)}
=\varepsilon_Y(a_{(2)}^{Y,Y})a^{X,Y}_{(1)}=a^{X,Y}$$
This proves that $\Delta_{X,Y}^Z$ is split injective. 
We know from Proposition \ref{consequence} that 
$\Delta_{X,Y}^Z$ is a $\coc(X,X)$-$\coc(Y,Y)$-bicomodule algebra
morphism and that
$\Delta_{X,Y}^Z(\coc(X,Y)) \subset \coc(X,Z) \square_{\coc(Z,Z)} \coc(Z,Y)$. 
Let $\sum_i a_i^{X,Z} \otimes b_i^{Z,Y}\in \coc(X,Z) \square_{\coc(Z,Z)} \coc(Z,Y)$.
We have
\begin{align*}
 \Delta_{X,Y}^Z\circ f &(\sum_i a_i^{X,Z} \otimes b_i^{Z,Y})=
\Delta_{X,Y}^Z\left(\sum_i\psi\left(S_{Y,Z}(a^{Y,Z}_{i(2)})b_i^{Z,Y}\right)a^{X,Y}_{i(1)}\right) \\
&= \sum_i\psi\left(S_{Y,Z}(a^{Y,Z}_{i(3)})b_i^{Z,Y}\right)a^{X,Z}_{i(1)} \otimes a_{i(2)}^{Z,X} = \sum_i\psi\left(S_{Y,Z}(b^{Y,Z}_{i(2)})b_{i(3)}^{Z,Y}\right)a^{X,Z}_{i}\otimes 
b_{i(1)}^{Z,Y} \\
&  = \sum_i \varepsilon_Y(b_{i(2)}^{Y,Y})a_{i(1)}^{X,Z} \otimes b_{i(1)}^{Z,Y}
 = \sum_i a_i^{X,Z} \otimes b_i^{Z,Y}
\end{align*}
which proves the result.
Assume now that $\coc(X,Z) \not=(0)$,
and let
$\phi : \coc(X,Z) \longrightarrow k$ be a linear map such that
$\phi(1)=1$. Define $g : \coc(X,Z) \otimes \coc(Z,Y) \longrightarrow  \coc(X,Y)$
by 
$$g(a^{X,Z}\otimes b^{Z,Y}) = \phi\left(a^{X,Z}S_{Z,X}(b^{Z,X}_{(1)})\right)b^{X,Y}_{(2)}$$
Then 
$$g \circ \Delta_{X,Y}^Z(a^{X,Y})
=g(a_{(1)}^{X,Z}\otimes a_{(2)}^{Z,Y})
= \phi\left(a_{(1)}^{X,Z}S_{Z,X}(a^{Z,X}_{(2)})\right)a^{X,Y}_{(3)}
=\varepsilon_X(a_{(1)}^{X,X})a^{X,Y}_{(2)}=a^{X,Y}$$
and hence $\Delta_{X,Y}^Z$ is split injective. The rest of the proof
is then similar to the previous case, and is left to the reader.
\end{proof}

We get a useful criterion to show that a cogroupoid
is connected.

\begin{proposition}\label{connec}
 Let $\coc$ be a cogroupoid. The following assertions are equivalent.
\begin{enumerate}
 \item $\coc$ is connected.
\item There exists $X_0 \in {\rm ob}(\coc)$ such that $\forall Y \in {\rm ob}(\coc)$, $\coc(X_0,Y) \not = (0)$.
\item There exists $X_0 \in {\rm ob}(\coc)$ such that $\forall Y \in {\rm ob}(\coc)$, $\coc(Y,X_0) \not = (0)$.
\end{enumerate}
\end{proposition}

\begin{proof}
 Assume that (2) holds. Then for $X,Y \in {\rm ob}(\coc)$, the previous
lemma ensures that
$$\coc(X_0,X) \otimes \coc(X,Y) \simeq \coc(X_0,Y) \oplus W$$
for some vector space $W$. Hence $\coc(X,Y) \not = (0)$ and $\coc$ is connected.

Assume that (3) holds. Then for $X,Y \in {\rm ob}(\coc)$, the previous
lemma ensures that
$$\coc(X,Y) \otimes \coc(Y,X_0) \simeq \coc(X,X_0) \oplus W'$$
for some vector space $W'$. Hence $\coc(X,Y) \not = (0)$ and $\coc$ is connected.
\end{proof}

We now provide a self-contained proof of Theorem \ref{explicit}.

\begin{proof}[Proof of Theorem \ref{explicit}]
Let $\coc$ be a connected cogroupoid and let $X,Y \in {\rm ob}(\coc)$.
By Proposition \ref{consequence} and Proposition \ref{basi2}  
we have two $k$-linear functors
 $$F : {\rm Comod}(\coc(X,X))\longrightarrow {\rm Comod}(\coc(Y,Y)), 
  V \longmapsto V \square_{\coc(X,X)} \coc(X,Y)$$
  $$G : {\rm Comod}(\coc(Y,Y))\longrightarrow {\rm Comod}(\coc(X,X)), V \longmapsto V \square_{\coc(Y,Y)} \coc(Y,X)$$
Let us prove that these functors are inverse equivalences.
Let $V \in {\rm Comod}(\coc(X,X))$. We have a $\coc(X,X)$-colinear isomorphism
$\theta_V : V \cong G\circ F(V)$ defined by the composition:
{\footnotesize
\begin{equation*}
\xymatrix{
V \ar[r]^{\cong} &
V\square_{\coc(X,X)}\coc(X,X) \ar[d]^{1 \otimes
\Delta_{X,X}^Y} \\ 
& V \square_{\coc(X,X)}\left(\coc(X,Y) \square_{\coc(Y,Y)} \coc(Y,X)\right) 
 \ar[r]^-{\cong} & \left(V \square_{\coc(X,X)}(\coc(X,Y)\right) \square_{\coc(Y,Y)} \coc(Y,X) = G\circ F(V)
}
\end{equation*}}
where the first isomorphism is induced by the $\coc(X,X)$-coaction on $V$,
the second map is an isomorphism by Lemma \ref{useful} ($\Delta_{X,X}^Y$ is 
$\coc(X,X)$-colinear), and the third isomorphism is the associativity isomorphism
of cotensor products. This isomorphism is clearly natural in $V$, so
we have a functor isomorphism ${\rm id} \cong G\circ F$, and we also
have ${\rm id} \cong F\circ G$, so $F$ and $G$ are inverse equivalences.

We know from Proposition \ref{basimonoidal}
that $F$ and $G$ are monoidal functors, and it is easy to
check that the above functor isomorphisms are isomorphisms of monoidal functors.
\end{proof}

\begin{remark}{\rm 
It is not difficult (and probably interesting)
to give a proof of the fact that $F$ and $G$ are monoidal
functors without using Proposition \ref{basimonoidal} (and hence Ulbrich's Theorem)
directly from the cogroupoid axioms.

We have to check that for $V,W \in {\rm Comod}(\coc(X,X))$ the map
\begin{align*}
 \widetilde{F}_{V,W} :(V \square_{\coc(X,X)} \coc(X,Y)) \otimes
(W \square_{\coc(X,X)} \coc(X,Y)) 
&\longrightarrow (V\otimes W) \square_{\coc(X,X)} \coc(X,Y) \\
(\sum_i v_i \otimes a_i^{X,Y}) \otimes
(\sum_j w_j \otimes b_j^{X,Y}) &\longmapsto
\sum_{i,j} v_i \otimes w_j\otimes a_i^{X,Y}b_j^{X,Y}
\end{align*}
is an isomorphism.
It is immediate to check that the following diagram commutes
$$\xymatrix{
GF(V) \otimes GF(W)\ar[rr]^{\widetilde{G}_{F(V),F(W)}}  &&
G(F(V)\otimes F(W)) \ar[d]^{G(\widetilde{F}_{V,W})}  \\ 
V \otimes W \ar[u]^{\theta_V \otimes \theta_W} \ar[rr]_{\theta_{V \otimes W}}
&& GF(V \otimes W) } $$
(with the previous notation)
and hence $G(\widetilde{F}_{V,W}) \circ \widetilde{G}_{F(V),F(W)}$ is an isomorphism.
The same reasoning for $G$ shows that
$F(\widetilde{G}_{F(V),F(W)}) \circ \widetilde{F}_{GF(V),GF(W)}$ is an isomorphism.
We conclude that $\widetilde{G}_{F(V),F(W)}$ is an isomorphism and hence
so is $\tilde{F}_{V,W}$.}
\end{remark}

\begin{remark}{\rm
 Let us say that a cocategory $\coc$ is a Takeuchi cocategory
if for all $X,Y,Z \in {\rm ob}(\coc)$ the algebra map
$$\Delta_{X,Y}^Z : \coc(X,Y)\longrightarrow\coc(X,Z) \square_{\coc(Z,Z)} \coc(Z,Y)$$
is an isomorphism.
This terminology comes from the fact
that every pair of objects in a Takeuchi cocategory produces a set of equivalence
data in the sense of \cite{tak} (a Morita-Takeuchi equivalence). Lemma \ref{useful} ensures that a connected cogroupoid
is a Takeuchi cocategory and the proof that
the functors of Theorem \ref{explicit} are inverse equivalences
is just the classical proof that a Morita-Takeuchi equivalence produces
an equivalence of categories.
The proof of monoidality we gave in the last remark shows that
if $X,Y$ are objects of a Takeuchi cocategory $\coc$
the comodule categories over the bialgebras $\coc(X,X)$ and $\coc(Y,Y)$
are monoidally equivalent. 
}
\end{remark}


\subsection{From Hopf-Galois objects and monoidal equivalences
to cogroupoids}

In this subsection we give (sketches of) the proofs of Theorem \ref{galoiscog}
and of the implication $(1) \Rightarrow (2)$ in Theorem
\ref{schacog}. The techniques 
used here are never used  in the rest of the paper, so
the reader who would prefer to see examples and applications
 might skip the subsection first, and get back
if he wants or needs to.

We use Tannaka-Krein reconstruction techniques (see \cite{js,sctan}).
Let us begin with the following general situation.
Let $\mathcal C$ be a small category
and 
let $F,G : \mathcal C \longrightarrow {\rm Vect}_f(k)$ be some functors.  
Following \cite{js}, Section 3, we associate a vector space
$\mathsf{Hom}^{\vee}\!(F,G)$ to such a pair:
$$
\mathsf{Hom}^{\vee}\!(F,G) = \bigoplus_{X \in {\rm ob}(\mathcal C)} 
{\rm Hom}_k(F(X),G(X))/\mathcal N$$
where
$\mathcal N$ is the linear subspace of $\bigoplus_{X \in {\rm ob}(\mathcal C)} {\rm
Hom}_k(F(X),G(X))$ generated by the elements
$G(f) \circ u - u \circ
F(f)$, with $f \in {\rm Hom}_{\mathcal C}(X,Y)$ and $u \in {\rm
Hom}_k(F(Y),G(X))$.
The class of an element $u \in {\rm Hom}_k(F(X), G(X))$ is denoted
by $[X,u]$ in  $\mathsf{Hom}^{\vee}\!(F,G)$
(note that we have changed the order of the functors in \cite{js}: our
$\mathsf{Hom}^{\vee}\!(F,G)$ is $\mathsf{Hom}^{\vee}\!(G,F)$ in \cite{js}).
The vector space $\mathsf{Hom}^{\vee}\!(F,G)$  has the following universal property.

\begin{lemma}\label{univcoend}
The vector space  $\mathsf{Hom}^{\vee}\!(F,G)$ represents the functor 
$${\rm Vect}_f(k) \longrightarrow {\rm Vect}(k), \quad 
V \longmapsto {\rm Nat}(G, F \otimes V)$$
More precisely, there exists $\alpha_\bullet \in {\rm Nat}(G,F \otimes  \mathsf{Hom}^{\vee}\!(F,G))$ such that the map
\begin{align*}
\Hom_k( \mathsf{Hom}^{\vee}\!(F,G), V) &\longrightarrow {\rm Nat}(G,F \otimes V) \\
\phi &\longmapsto (1 \otimes \phi) \circ \alpha_\bullet
\end{align*}
is a bijection.
\end{lemma}

\begin{proof}
Let $X$ in ${\rm ob}(\coc)$ and let $e_1, \ldots , e_n$
be a basis of $F(X)$. Define $$\alpha_X : G(X) \longrightarrow F(X) \otimes  \mathsf{Hom}^{\vee}\!(F,G)$$ by $\alpha_X(x) = \sum_{i}e_i \otimes [X, e_i^* \otimes x]$. It is easily seen that $\alpha_X$ does not depend on the choice of a basis of $G(X)$, and that this procedure defines an element $\alpha_\bullet \in {\rm Nat}(G,F \otimes  \mathsf{Hom}^{\vee}\!(F,G))$. It is not difficult to check
that the map in the statement of the lemma is a bijection.
\end{proof}

The universal property of $\mathsf{Hom}^{\vee}\!(F,G)$ gives, for any functor
$K : \mathcal C \longrightarrow {\rm Vect}_f(k)$, a linear map 
$$\Delta_{F,G}^K : \mathsf{Hom}^{\vee}\!(F,G)
\longrightarrow \mathsf{Hom}^{\vee}\!(F,K) 
\otimes \mathsf{Hom}^{\vee}\!(K,G)$$
coassociative in the sense of cocategories. The map $\Delta_{F,G}^K$
may be described as follows. Let $X \in {\rm ob}(\mathcal C)$, let
$\phi \in F(X)^*$, let $x \in G(X)$ and let $e_1, \ldots , e_n$
be a basis of $K(X)$. Then
$$\Delta_{F,G}^K([X, \phi \otimes x])=
\sum_{i=1}^n [X,\phi \otimes e_i] \otimes [X, e_i^* \otimes x].$$   
As a particular case of the previous construction,
$\mathsf{End}^{\vee}\!(F) := \mathsf{Hom}^{\vee}\!(F,F)$ is a coalgebra, with counit
$\varepsilon_X$
induced by the trace: $\varepsilon_X([X,u]) = {\rm tr}(u)$ for $u \in \Hom(F(X),F(X))$.

The previous construction enables one to reconstruct a coalgebra from
its category of finite-dimensional comodules and the forgetful functor: this the Tannaka reconstruction theorem.

\begin{theorem}
 Let $C$ be a coalgebra and let
$\Omega^C : {\rm Comod}_f(C) \longrightarrow {\rm Vect}_f(k)$ be the forgetful functor. We have a coalgebra isomorphism $ \mathsf{End}^{\vee}\!(\Omega^C) \cong C$.
\end{theorem}

\begin{proof}
 Consider the natural transformation $\alpha^C :\Omega^C \longrightarrow \Omega^C \otimes C$ induced by the coactions of $C$ on its comodules. The universal property
of  $ \mathsf{End}^{\vee}$ yields a unique linear map $f : \mathsf{End}^{\vee}\!(\Omega^C) \longrightarrow C$ such that the following diagram commutes
$$\xymatrix{
\Omega^C \ar[rr]^{\alpha_\bullet} \ar[dr]_{\alpha^C} & & \Omega^C \otimes 
\mathsf{End}^{\vee}\!(\Omega^C) \ar[dl]^{1 \otimes f}\\
& \Omega^C \otimes C &
}
$$
It is easy to see that $f$ is a coalgebra morphism.
To prove that $f$ is an isomorphism, we proceed by following \cite{sctan}, Lemma 2.2.1 (for another proof see Section 6 in \cite{js}). Consider, 
for a vector space $V$, the linear map
\begin{align*}
\Phi : \Hom_k(C, V) &\longrightarrow {\rm Nat}(\Omega^C,\Omega^C \otimes V) \\
\phi &\longmapsto (1 \otimes \phi) \circ \alpha^C
\end{align*}
Let us prove that $\Phi$ is bijection. This will define a linear map $C \longrightarrow \mathsf{End}^{\vee}\!(\Omega^C)$ which necessarily will be an inverse of $f$. To construct the inverse of $\Phi$, the key remark
is that if $\varphi \in  {\rm Nat}(\Omega^C,\Omega^C \otimes V)$ and $N$, $M$ are two finite-dimensional subcomodules of $C$, then
$$((\varepsilon \otimes 1) \circ \varphi_{N})_{|N \cap M} = (\varepsilon \otimes 1) \circ \varphi_{N\cap M} =
((\varepsilon \otimes 1) \circ \varphi_{M})_{|N\cap M}$$
This follows from the naturality of $\varphi$. This enables to define a (linear) map
$$\Psi : {\rm Nat}(\Omega^C,\Omega^C \otimes V) \longrightarrow \Hom_k(C, V)$$ by
$\Psi(\varphi)(x) = (\varepsilon \otimes 1) \circ \varphi_{M}(x)$, $\forall x \in C$, where $M$ is any finite-dimensional subcomodule of $C$ containing $x$.

For $\phi \in {\rm Hom}_k(C,V)$, $x \in C$ and $M$ a finite-dimensional subcomodule of $C$ containing $x$, we have 
$$\phi(x) = \phi ((1 \otimes \varepsilon)  \circ \alpha_M^C(x)) =
(\varepsilon \otimes 1) \circ (1 \otimes \phi) \circ \alpha_M^C(x)=\Psi(\Phi(\phi))(x)$$
which proves that $\Psi \circ \Phi = {\rm id}$.

Let $\varphi \in {\rm Nat}(\Omega^C,\Omega^C \otimes V)$. Let $M$ be a finite-dimensional comodule and let $D$ be a finite-dimensional subcoalgebra of $C$ such that $\alpha^C(M) \subset M \otimes D$. Then $D$ is a subcomodule of $C$. Denote by $M_0 \otimes D$ the $C$-comodule whose coaction is given by $1 \otimes \Delta$.
It is clear that $\alpha_M^C : M \longrightarrow M_0 \otimes D$ is a $C$-comodule map, hence by naturality of $\varphi$ the following diagram commutes
$$\xymatrix{
M \ar[r]^{\varphi_M} \ar[d]^{\alpha_M^C} & M \otimes V \ar[d]^{\alpha_M^C \otimes 1} \\
M_0 \otimes D \ar[r]^{\varphi_{M_0\otimes D}} & M_0 \otimes D \otimes V 
}
$$ 
Any linear map $D \longrightarrow M_0 \otimes D$ is $C$-colinear, hence again the naturality of $\varphi$ shows that $\varphi_{M_0\otimes D}= 1 \otimes \varphi_D$.  
Hence we have
\begin{align*}
\Phi\circ \Psi(\varphi)_M&= (1 \otimes \Psi(\varphi)) \circ \alpha_M^C
=(1 \otimes \Psi(\varphi)_{|D}) \circ \alpha_M^C \\
&= (1 \otimes \varepsilon \otimes 1) \circ (1 \otimes \varphi_D) \circ \alpha_M^C
= (1 \otimes \varepsilon \otimes 1) \circ (\alpha_M^C \otimes 1) \circ \varphi_M
= \varphi_M
\end{align*}
and this proves that $\Phi \circ \Psi ={\rm id}$. We conclude that $\Phi$ is an isomorphism.
\end{proof}

\smallskip

Now assume that $\mathcal C$ is a monoidal category and that
$F,G$ are monoidal functors. Then 
$\mathsf{Hom}^{\vee}\!(F,G)$ inherits an algebra structure, whose
product may be described by the following formula:
$$[X,u].[Y,v] = [X \otimes Y, \widetilde G_{X,Y} \circ (u \otimes v)
\circ \widetilde F_{X,Y}^{-1}]$$
where the isomorphisms $\widetilde F_{X,Y} : F(X) \otimes F(Y) 
\longrightarrow F(X \otimes Y)$ and 
$\widetilde G_{X,Y} : G(X) \otimes G(Y) 
\longrightarrow G(X \otimes Y)$ are the constraints of the monoidal functors
$F$ and $G$. The unit element is $[I,\widetilde G_0 \circ \widetilde F_0^{-1}]$ where $I$ stands for the monoidal unit of $\coc$ (that we indeed have an associative algebra structure heavily depends
on the fact that $F$ and $G$ are monoidal functors).
It is not difficult to check that the maps
$\Delta_{F,G}^K$ and $\varepsilon_F$ are algebra maps. In particular $\mathsf{End}^{\vee}\!(F)$ is a bialgebra.

We summarize the above constructions as follows.

\begin{definition}
 Let $\coc$ be a monoidal category. The cocategory ${\rm Mon}_k(\coc)$ is the cocategory whose objects are the monoidal functors $\mathcal C \longrightarrow {\rm Vect}_f(k)$, with ${\rm Mon}_k(\coc)(F,G) =  \mathsf{Hom}^{\vee}\!(F,G)$ for 
$F,G \in {\rm Mon}_k(\coc)$, and with structural maps
$\Delta_{\bullet,\bullet}^\bullet$ and $\varepsilon_\bullet$ defined above.
\end{definition}

We shall need a monoidal version of Lemma \ref{univcoend}.
If $A$ is an algebra and $F,G : \coc \longrightarrow {\rm Vect}_f(k)$
are monoidal functors, denote by ${\rm Nat}_\otimes(G,F\otimes A)$
the set of elements $\theta \in {\rm Nat}(G,F\otimes A)$ such
that the following diagrams commute for any objects $X$,$Y$
$$\xymatrix{
G(X\otimes Y) \ar[r]^{\theta_{X\otimes Y}} & F(X\otimes Y)\otimes A \ar[r]^{\widetilde F_{X,Y}^{-1}} & F(X)\otimes F(Y) \otimes A \\
G(X) \otimes G(Y) \ar[u]^{\widetilde G_{X,Y}} \ar[r]^-{\theta_X \otimes \theta_Y} &
F(X)\otimes A \otimes F(Y) \otimes A \ar[r]^{1 \otimes \tau \otimes 1} & F(X) \otimes F(Y) \otimes A \otimes A
\ar[u]_{1 \otimes 1 \otimes m_A}
}$$
$$\xymatrix{
F(I) \ar[r]^{\theta_I}  & G(I)  \otimes A \\
k \ar[u]^{\widetilde F_0} \ar[r]^-{u} & A \cong k \otimes A \ar[u]_{\widetilde G_0\otimes 1}
}$$
 
\begin{lemma}\label{univcoendmono}
Consider the element $\alpha_\bullet \in {\rm Nat}(G,F \otimes  \mathsf{Hom}^{\vee}\!(F,G))$ 
of Lemma \ref{univcoend}. We have $\alpha_\bullet \in {\rm Nat}_{\otimes}(G,F \otimes  \mathsf{Hom}^{\vee}\!(F,G))$ 
and we have, for any algebra $A$, a map
\begin{align*}
\Hom_{k-{\rm alg}}( \mathsf{Hom}^{\vee}\!(F,G), A) &\longrightarrow {\rm Nat}_{\otimes}(G,F \otimes A) \\
\phi &\longmapsto (1 \otimes \phi) \circ \alpha_\bullet
\end{align*}
which is a bijection.
\end{lemma}


\begin{proof}
 The proof is straightforward.
\end{proof}

The bialgebra version of the Tannaka duality theorem is as follows.

\begin{theorem}\label{tanbial}
 Let $B$ be a bialgebra and let
$\Omega^B : {\rm Comod}_f(B) \longrightarrow {\rm Vect}_f(k)$ be the forgetful functor. We have a bialgebra isomorphism $\mathsf{End}^{\vee}\!(\Omega^B) \cong B$.
\end{theorem}

\begin{proof}
 The natural transformation $\alpha^B :\Omega^B \longrightarrow \Omega^B \otimes B$ induced by the coactions of $B$ on its comodules
is an element in
${\rm Nat}_{\otimes}(\Omega^B,\Omega^B \otimes B)$, by definition of the tensor product of $B$-comodules. Hence Lemma \ref{univcoendmono} yields a unique algebra map map $f : \mathsf{End}^{\vee}\!(\Omega^B) \longrightarrow B$ such that the following diagram commutes
$$\xymatrix{
\Omega^C \ar[rr]^{\alpha_\bullet} \ar[dr]_{\alpha^C} & & \Omega^C \otimes 
\mathsf{End}^{\vee}\!(\Omega^C) \ar[dl]^{1 \otimes f}\\
& \Omega^C \otimes C &
}
$$ We know from the proof of Lemma \ref{univcoendmono} that $f$ is a colagebra
isomorphism, and hence is a bialgebra isomorphism.
\end{proof}

Assume moreover that $\mathcal C$ is a rigid monoidal category. This
means that every object $X$ has a left dual (\cite{js, kas}), 
i.e. there exist a triplet $(X^*, e_X, d_X)$ where 
$X^* \in {\rm ob}(\mathcal C)$, while 
$e_X:  {X^* } \otimes X \longrightarrow I$ ($I$ is the monoidal
unit of $\mathcal C$) and 
$d_X : I \longrightarrow X \otimes  {X^* }$ are morphisms 
of $\mathcal C$ such that:
$$ 
(1_X \otimes e_X) \circ (d_X \otimes 1_X) = 1_X \quad 
{\rm and} \quad (e_X \otimes 1_{X^*})
 \circ (1_{X^*} \otimes d_X) = 1_{X^*}$$
The rigidity of $\mathcal C$ allows one to define a duality
endofunctor of $\mathcal C$, which will be used in the proof of
the following result, which generalizes \cite{u90}, using the same idea.

\begin{proposition}
 Let $\coc$ be a rigid monoidal category.
Then the cocategory ${\rm Mon}_k(\coc)$ has a cogroupoid structure.
\end{proposition}

\begin{proof}
We have to construct
the linear maps $S_{F,G}: \mathsf{Hom}^{\vee}\!(F,G) \longrightarrow
\mathsf{Hom}^{\vee}\!(G,F)$. 
We fix for every $X \in {\rm ob}(\mathcal C)$ a left dual $(X^*, e_X, d_X)$ (with $I^*=I$ for the monoidal unit). This defines a contravariant endofunctor
of $\coc$: $X \longmapsto X^*$, $f \longmapsto f^*$, where for 
$f : X \longrightarrow Y$, the morphism $f^* : Y^* \longrightarrow X^*$
is defined by the composition
$$\xymatrix{
Y^* \ar[r]^-{1 \otimes d_X} & Y^* \otimes (X \otimes X^*) \ar[r]^{1\otimes (f \otimes 1)} & Y^* \otimes (Y \otimes X^*) \ar[r]^{\sim} &(Y^* \otimes Y) \otimes X^*
\ar[r]^{e_Y \otimes 1} & I \otimes X^* \cong X^*
}$$
Let $X \in {\rm ob}(\mathcal C)$. The unicity
of a left dual in a monoidal category yields natural isomorphisms 
$$\lambda_X^F : F(X)^* \longrightarrow F(X^*) \quad {\rm and} \quad
\lambda_X^G : G(X)^* \longrightarrow G(X^*)$$
such that the following diagrams commute:
 $$\xymatrix{
F(X)^* \otimes F(X) \ar[r]^-{e_{F(X)}} \ar[d]^{\lambda_X^F \otimes 1_{F(X)}} &
I \ar[r]^{\widetilde F_0}  & F(I) \\
F(X^*) \otimes F(X) \ar[rr]^-{\widetilde F_{X^*,X}} && F(X^* \otimes X) \ar[u]^{F(e_X)}} \quad \quad
\xymatrix{
F(X) \otimes F(X)^*  \ar[d]^{1_{F(X)} \otimes \lambda_X^F} &
I  \ar[l]_-{d_{F(X)}} \ar[r]^{\widetilde F_0}  & F(I) \ar[d]^{F(d_X)} \\
F(X) \otimes F(X^*) \ar[rr]^-{\widetilde F_{X,X^*}} && 
F(X \otimes X^*)}
$$
(we have endowed ${\rm Vect}_f(k)$ with its standard duality).
Let $u \in {\rm Hom}_k(F(X), G(X))$. We put 
$$S_{F,G}([X,u]) = [X^*, \lambda_X^F \circ u^* \circ (\lambda_X^G)^{-1}]$$
It is easy to see that $S_{F,G}$ is a well defined linear map.
Now let $\phi \in F(X)^*$, let $x \in F(X)$ and let $e_1, \ldots, e_n$
be a basis of $G(X)$. Then we have
\begin{align*}
m & \circ (1 \otimes S_{G,F}) \circ \Delta_{F,F}^G ([X, \phi \otimes x]) 
= \sum_{i=1}^n [X, \phi \otimes e_i] [X^*, \lambda_X^G \circ
(e_i^* \otimes x) \circ (\lambda_X^F)^{-1}] = \\
= & \sum_{i=1}^n
[X \otimes X^*, \widetilde G_{X,X^*} \circ
(1_{G(X)} \otimes \lambda_X^G) \circ
((\phi \otimes e_i) \otimes (x \otimes e_i^*))
\circ (1_{F(X)} \otimes (\lambda_X^F)^{-1}) \circ
\widetilde F_{X,X^*}^{-1}] \\
= & [X \otimes X^*, \widetilde G_{X,X^*} \circ
(1_{G(X)} \otimes \lambda_X^G) \circ d_{G(X)} \circ (\phi \otimes x)
\circ (1_{F(X)} \otimes (\lambda_X^F)^{-1}) \circ
\widetilde F_{X,X^*}^{-1}] \\
= & [X \otimes X^*, G(d_X) \circ \widetilde G_0 \circ (\phi \otimes x) 
\circ (1_{F(X)} \otimes (\lambda_X^F)^{-1}) \circ
\widetilde F_{X,X^*}^{-1}] \\
= & [I, \widetilde G_0 \circ (\phi \otimes x)
\circ (1_{F(X)} \otimes (\lambda_X^F)^{-1}) \circ
\widetilde F_{X,X^*}^{-1} \circ F(d_X)] \\ 
= & [I, \widetilde G_0 \circ (\phi \otimes x) \circ d_{F(X)}
\circ \widetilde F_0^{-1}] 
= \phi(x) [I,\widetilde G_0 \circ \widetilde F_0^{-1}]
=\varepsilon_F([X, \phi \otimes x])1.\\
\end{align*}
Since the elements $[X, \phi \otimes x]$ linearly span 
$\mathsf{End}^{\vee}\!(F)$, we have the commutativity of the first
diagram. The commutativity of the second diagram is proved
similarly.
\end{proof}

\begin{definition}
 Let $H$ be a Hopf algebra. The \textbf{fibre functor cogroupoid of $H$}, denoted
${\rm Fib}(H)$, is the full subcogroupoid of ${\rm Mon}_k({\rm Comod}_f(H))$
whose objects are the fibre functors on ${\rm Comod}_f(H)$.
\end{definition}

We have now all the ingredients to prove Theorem \ref{galoiscog}.

\begin{proof}[Proof of Theorem \ref{galoiscog}]
Let $H$ be a Hopf algebra and let $A$ be a left $H$-Galois object. Consider the full subcogroupoid $\mathcal D$ of ${\rm Fib}(H)$ whose objects are $\Omega$ (the forgetful functor) and $\Omega^A$ (see Ulbrich's theorem). 
Let us check that $\mathcal D$ is connected.
Let $\phi : A \longrightarrow k$ be a linear map such that $\phi(1)=1$.
For any $X \in {\rm Comod}_f(H)$, the maps
$1 \otimes \phi : X\square_H A \rightarrow X$
define a non-zero element in ${\rm Nat}(\Omega^A , \Omega)$.
Hence by Lemma \ref{univcoend} we have $\mathsf{Hom}^{\vee}\!(\Omega,\Omega^A)\not=(0)$
and the cogroupoid $\mathcal D$ is connected by Proposition \ref{connec}.

For any $X \in {\rm Comod}_f(H)$, the inclusions $X\square_H A \subset X \otimes A$
define an element $$\beta_\bullet\in{\rm Nat}_{\otimes}(\Omega^A,\Omega\otimes A)$$
which corresponds by Lemma \ref{univcoendmono} to an algebra map
$g :  \mathsf{Hom}^{\vee}\!(\Omega,\Omega^A) \longrightarrow A$, and we leave it to the reader to check that the following diagram commutes
$$\xymatrix{
\mathsf{Hom}^{\vee}\!(\Omega,\Omega^A) \ar[r]^g 
\ar[d]^{\Delta_{\Omega, \Omega^A}^{\Omega}} & A \ar[d]^\rho \\
\mathsf{Hom}^{\vee}\!(\Omega,\Omega) \otimes \mathsf{Hom}^{\vee}\!(\Omega,\Omega^A)
\ar[r]^-{f \otimes g}  & H \otimes A
}$$
where $\rho$ stands for the coaction of $H$ on $A$ and $f$ is the bialgebra
isomorphism  of Theorem \ref{tanbial}.
This means that if we endow $\mathsf{Hom}^{\vee}\!(\Omega,\Omega^A)$
with the left $H$-comodule algebra structure transported from the Hopf algebra isomorphism $f$, then $g$ is an $H$-comodule algebra morphism. But $\mathsf{Hom}^{\vee}\!(\Omega,\Omega^A)$ is $H$-Galois since it is 
$\mathsf{Hom}^{\vee}\!(\Omega,\Omega)$-Galois (Proposition \ref{back}), and hence
by Proposition \ref{galgroupo} $g$ is an isomorphism. We can now consider
the connected cogroupoid $\mathcal D_0$  with two objects
$X,Y$ and
$$\mathcal D_0(X,X)=H, \ \mathcal D_0(X,Y)=A, \ \mathcal D_0(Y,X)=\mathsf{Hom}^{\vee}\!(\Omega^A,\Omega), \
 \mathcal D_0(Y,Y) =\mathsf{Hom}^{\vee}\!(\Omega^A,\Omega^A)$$
The structural maps of $\mathcal D_0$ are transported from the connected cogroupoid
$\mathcal D$ via the isomorphisms $f$ and $g$, and we are done. 
\end{proof}

We get the following result as a corollary of the proof.

\begin{proposition}
 Let $H$ be a Hopf algebra. The cogroupoid  ${\rm Fib}(H)$
is connected.
\end{proposition}

\begin{proof}
 If $A$ is a left $H$-Galois object, we have seen in the previous proof that $\mathsf{Hom}^{\vee}\!(\Omega,\Omega^A)\not=(0)$. Hence the result follows from Ulbrich's theorem (any fibre fonctor is isomorphic to some $\Omega^A$)
and Proposition \ref{connec}.
\end{proof}

\begin{proof}[Proof of Theorem \ref{schacog}]
 Let $H$, $L$ be some Hopf algebras and let 
$F : {\rm Comod}(H) \longrightarrow {\rm Comod}(L)$ be a $k$-linear monoidal equivalence. Then $F$ induces a monoidal equivalence 
${\rm Comod}_f(H) \longrightarrow {\rm Comod}_f(L)$, still denoted $F$
(see Proposition \ref{fini}). 
It is clear that $\Omega^L \circ F$ is a fibre functor on ${\rm Comod}_f(H)$.
It is not difficult to construct a linear map
$\mathsf{End}^{\vee}\!(\Omega^L)\longrightarrow \mathsf{End}^{\vee}\!(\Omega^L\circ F)$
which is a Hopf algebra map since $F$ is monoidal and is an isomorphism
since $F$ is an equivalence.
Let $\mathcal D$ be the (connected) subcogroupoid of ${\rm Fib}(H)$ whose objects
are $\Omega^H$ and $\Omega^L \circ F$. By Theorem \ref{tanbial}
we have Hopf algebra isomorphisms $\mathsf{End}^{\vee}\!(\Omega^H)\cong H$
and $\mathsf{End}^{\vee}\!(\Omega^L)\cong L$, and hence 
$\mathsf{End}^{\vee}\!(\Omega^L\circ F)\cong L$. 
Therefore we get, by transporting the appropriate structures from the connected cogroupoid $\mathcal D$, a cogroupoid $\mathcal D_0$ with two objects $X,Y$ such that
$H=\mathcal D_0(X,X)$ and $L=\mathcal D_0(Y,Y)$.
\end{proof}

\subsection{The weak Hopf algebra of a finite cogroupoid}
In this short subsection we connect the theory of cogroupoids (and hence of Hopf-Galois
objects) with the theory of weak Hopf algebras \cite{bonisz}
(we do not recall here the precise definition of a weak Hopf algebra).
Since cogroupoids are non commutative generalizations of groupoids, it is 
natural to wonder if they are linked with weak Hopf algebras, one of the most
well-known non commutative generalization of groupoids. Not surprisingly,
the following construction shows that
this is indeed the case. The (unpublished)
result was obtained in collaboration with Grunspan in 2003-2004.

\begin{proposition}\label{coweak}
 Let $\coc$ be a cogroupoid and let $X_1, \ldots , X_n \in {\rm ob}(\coc)$.
For  $i,j \in \{1, \ldots ,n\}$, Put $\coc(i,j) = \coc(X_i, X_j)$, $\Delta_{i,j}^k = \Delta_{X_i,X_j}^{X_k}$,
$S_{i,j}=S_{X_i,X_j}$, $\varepsilon_i = \varepsilon_{X_i}$. Consider the 
direct sum of algebras 
$$H = \bigoplus_{i,j=1}^n \coc(i,j)$$
Then $H$ has a weak Hopf algebra structure defined as follows.
\begin{enumerate}
 \item The comultiplication $\Delta : H \longrightarrow H \otimes H$
is defined by 
$$\Delta(a^{i,j}) = \sum_{k=1}^n \Delta_{i,j}^k(a^{i,j})=\sum_{k=1}^n a^{i,k}_{(1)} \otimes a^{k,j}_{(2)}$$
\item The counit $\varepsilon : H \longrightarrow k$ is defined
by $\varepsilon_{|\coc(i,j)} = \delta_{ij}\varepsilon_{i}$
\item The antipode $S: H \longrightarrow H$ is defined by $S_{|\coc(i,j)} = S_{i,j}$.
\end{enumerate}
\end{proposition}

The proof is done by a straightforward verification.

\begin{remark}{\rm 
\begin{enumerate}
\item The use of quantum groupoids (in a different framework)
in the study of monoidal equivalences also arose (much earlier)
 in the work of Brugui\`eres \cite{bru}.
 \item It is also worth to note that De Commer \cite{dec1,dec2} obtained independently
more or less the same result in his investigation
of Galois objects for multiplier Hopf algebras and operator algebraic
quantum groups.
\end{enumerate}}
\end{remark}

\section{Examples of connected cogroupoids}

This section is devoted to the presentation of some examples of connected
cogroupoids. We hope that these examples will convince
that the reader that very often in concrete situations 
it is easy and natural to work with cogroupoids.

\subsection{The bilinear cogroupoid $\mathcal B$}
Let $E \in \GL_n(k)$.
The algebra
$\mathcal B(E)$ is the algebra presented by generators
$(a_{ij})_{1 \leq i,j \leq n}$ submitted to the relations
$$E^{-\!1} a^t E a = I_n = a E^{-\!1} a^t E,$$
where $a$ is the matrix $(a_{ij})_{1 \leq i,j \leq n}$, $a^t$ is the transpose matrix  and $I_n$ is the $n \times n$ identity matrix.

It admits a Hopf algebra structure, with
comultiplication $\Delta$ defined by 
$$\Delta(a_{ij})
= \sum_{k=1}^n a_{ik} \otimes a_{kj}, \ 
\varepsilon(a_{ij}) = \delta_{ij}, \ 
S(a) = E^{-1}a^t E$$
This Hopf algebra was introduced by Dubois-Violette and Launer \cite{dvl}
as the \textsl{quantum automorphism group of the non degenerate bilinear form
associated to $E$}. This terminology comes from the following result.

\begin{proposition}[Universal property of $\B(E)$]

\ 
\begin{enumerate}
\item Consider the vector space $V = k^n$ with its canonical basis 
$(e_i)_{1 \leq i \leq n}$. Endow $V$ with the $\B(E)$-comodule structure 
defined by $\alpha(e_i) = \sum_{j=1}^n e_j \otimes a_{ji}$, 
$1 \leq i \leq n$. Then the linear map $\beta : V \otimes V \longrightarrow k$
defined by $\beta (e_i \otimes e_j) = \lambda_{ij}$, ${1 \leq i,j \leq n}$,
where $E = (\lambda_{ij})$, is a $\B(E)$-comodule morphism.

\item Let $H$ be a Hopf algebra and let $V$ be a finite-dimensional 
$H$-comodule of dimension $n$. 
Let $\beta : V \otimes V \longrightarrow k$ be an $H$-comodule morphism 
such that the associate bilinear form is non-degenerate.
Then there exists $E \in \GL_n(k)$ such that $V$ is a $\B(E)$-comodule,
that $\beta$ is a $\B(E)$-comodule morphism, and that there exists 
a unique Hopf algebra morphism $\phi : \B(E) \longrightarrow H$
with $({\rm id}_V \otimes \phi) \circ \alpha = \alpha'$, where $\alpha$
and $\alpha'$ denote the coactions on $V$ of $\B(E)$ and $H$ respectively. 
\end{enumerate}
\end{proposition}

\begin{proof}
 The proof is left as an exercise.
\end{proof}

It is not difficult to check that
$\mathcal O_q(\SL_2(k)) = \B(E_q)$, where 
$$E_q =  \left(\begin{array}{cc} 0 & 1 \\
                          -q^{-1} & 0\\
       \end{array} \right)$$ and hence the Hopf algebras
$\B(E)$ are generalizations of $\mathcal O_q(\SL_2(k))$.

\medskip

We now will describe $\B(E)$ as part of a cogroupoid.
We first need a version involving two matrices, as follows.

Let $E \in \GL_m(k)$ and let $F \in \GL_n(k)$.
The algebra $\mathcal B(E,F)$ is the universal algebra with generators
$a_{ij}$, $1\leq i \leq m,1\leq j \leq n$, 
satisfying the relations
$$F^{-\!1} a^t E a = I_n \ ; \   a F^{-\!1}  a^t E = I_m.$$
Of course the generator $a_{ij}$ in $\B(E,F)$ 
should be denoted $a_{ij}^{E,F}$ to express the dependence on $E$ and $F$, but 
when there is  no confusion and we simply denote it by $a_{ij}$.
 It is clear that $\B(E,E)=\B(E)$

In the following lemma we construct the structural maps
that will put the algebras $\B(E,F)$ in a cogroupoid framework. 

\begin{lemma}\label{strucB}
\begin{enumerate}
 \item For any $E\in \GL_m(k)$, $F\in \GL_n(k)$, $G \in \GL_p(k)$, 
there exists an algebra map
\begin{align*}
\Delta_{E,F}^G : \B(E,F) &\longrightarrow \B(E,G) \otimes \B(G,F)\\
a_{ij} &\longmapsto \sum_{k=1}^p a_{ik} \otimes a_{kj}
\end{align*}
and for any $M \in \GL_r(k)$, the following diagrams commute
\begin{equation*}
\begin{CD}
\B(E,F) @>\Delta_{E,F}^G>> \B(E,G) \otimes \B(G,F) \\
@V\Delta_{E,F}^MVV
@V\Delta_{E,G}^M \otimes 1VV\\
\B(E,M) \otimes \B(M,F)  @>1 \otimes \Delta_{M,F}^G>> 
\B(E,M)\otimes \B(M,G)  \otimes \B(G,F)
\end{CD}
\end{equation*}
$$\xymatrix{
\B(E,F) \ar[d]^{\Delta_{E,F}^F} \ar@{=}[rrd]^{\rm id} \\
 \B(E,F) \otimes \B(F,F) \ar[rr]^{1 \otimes \varepsilon_F} && \B(E,F)}
\quad
\xymatrix{
\B(E,F) \ar[d]^{\Delta_{E,F}^E} \ar@{=}[rrd]^{{\rm id}} \\
 \B(E,E) \otimes \B(E,F) \ar[rr]^{\varepsilon_E \otimes 1} && \B(E,F)}
$$ where $\varepsilon_E$ is the counit of $\B(E)$.

\item For any $E\in \GL_m(k)$, $F\in \GL_n(k)$, there exists an algebra map 
\begin{align*}
S_{E,F} :\B(E,F) & \longrightarrow \B(F,E)^{\rm op} \\
a & \longmapsto E^{-1}a^tF
\end{align*}
such that the following diagrams commute
{\footnotesize
$$\xymatrix{
\B(E,E) \ar[r]^{\varepsilon_E} \ar[d]^{\Delta_{E,E}^F} &
k \ar[r]^{u}  & \B(E,F) \\
\B(E,F) \otimes \B(F,E) \ar[rr]^{1 \otimes S_{F,E}} && \B(E,F) \otimes \B(E,F) \ar[u]^{m}} \quad \quad
\xymatrix{
\B(E,E) \ar[r]^{\varepsilon_E} \ar[d]^{\Delta_{E,E}^F} &
k \ar[r]^{u}  & \B(F,E) \\
\B(E,F) \otimes \B(F,E) \ar[rr]^{S_{E,F} \otimes 1} && \B(F,E) \otimes \B(F,E) \ar[u]^{m}}
$$}
\end{enumerate}
\end{lemma}

\begin{proof}
 It is a straighforward verification to construct 
the announced algebra maps: this is left to the reader. 
The maps involved in the diagrams of part (1) all are algebra maps, and hence 
it is enough to check the commutativity on the generators of $\B(E,F)$, which
is obvious. 
For the diagrams of part (2), the commutativity
follows from the immediate verification on the generators of $\B(E,E)$
and the fact that
$\Delta_{\bullet,\bullet}^\bullet$
and $S_{\bullet,\bullet}$ are algebra maps.
\end{proof}

Hence the lemma ensures that we have a cogroupoid. 

\begin{definition}
 The cogroupoid $\B$ is the cogroupoid defined as follows:
\begin{enumerate}
 \item ${\rm ob}(\B) = \{ E \in \GL_n(k), \ n \geq 1\}$,
\item For $E,F \in {\rm ob}(\B)$, the algebra $\B(E,F)$ is the algebra defined above,
\item The structural maps $\Delta_{\bullet,\bullet}^\bullet$, $\varepsilon_\bullet$
and $S_{\bullet,\bullet}$ are defined in the previous lemma.
\end{enumerate}
\end{definition}

So we have a cogroupoid linking all the Hopf algebras $\B(E)$, and the natural next
question is to study the connectedness of $\B$.

\begin{lemma}
 Let $E \in \GL_m(k)$, $F \in \GL_n(k)$ with $m,n\geq 2$. Then $\B(E,F) \not =(0)$ if and only 
if ${\rm tr}(E^{-1}E^t)={\rm tr}(F^{-1}F^t)$.
\end{lemma}

\begin{proof}
 It is left as an exercise to check that if $\B(E,F) \not =(0)$ then 
 ${\rm tr}(E^{-1}E^t)={\rm tr}(F^{-1}F^t)$. Conversely, assume
that ${\rm tr}(E^{-1}E^t)={\rm tr}(F^{-1}F^t)$. To show that $\mathcal B(E,F) \not =(0)$
we can assume that $k$ is algebraically closed and hence that there
exists $q \in k^*$ such that ${\rm tr}(E^{-1}E^t)=-q-q^{-1}={\rm tr}(F^{-1}F^t)$.
It is shown in \cite{bi1}, by using the diamond lemma, that $\B(E_q,F) \not=(0)$
and hence by Proposition \ref{connec} we conclude that
$\B(E,F) \not = (0)$.
\end{proof}

\begin{corollary}\label{noncleft}
 Let $\lambda \in k$. Consider the full subcogroupoid $\B^\lambda$ of $\B$ with objects
$${\rm ob}(\B^\lambda)= \{ E \in \GL_n(k), \ n\geq 2,  {\rm tr}(E^{-1}E^t)=\lambda\}$$
Then $\B^\lambda$ is a connected cogroupoid.

In particular for $E \in \GL_m(k)$, $F \in \GL_n(k)$ with $n,m \geq 2$
and  ${\rm tr}(E^{-1}E^t)={\rm tr}(F^{-1}F^t)$, then
$\B(E,F)$ is a $\B(E)$-$\B(F)$-bi-Galois object and is not cleft
if $m \not = n$.
\end{corollary}

\begin{proof}
 It follows from the previous lemma and Proposition \ref{back}
that $\B(E,F)$ is a $\B(E)$-$\B(F)$-bi-Galois object if ${\rm tr}(E^{-1}E^t)={\rm tr}(F^{-1}F^t)$. Let us check that it is not cleft if $m \not =n$.
For $E \in \GL_m(k)$, let $V_E$ be the $m$-dimensional $\B(E)$-comodule
with basis $v_1^E, \ldots ,v_m^E$ and with $\B(E)$-coaction
$\alpha(v_i^E)= \sum_j v_j^E \otimes a_{ji}$. Let 
\begin{align*}
\Theta : {\rm Comod}(\B(E)) &\cong^\otimes {\rm Comod}(\B(F))
\quad \quad \quad \quad \quad \quad  \Theta' : &{\rm Comod}(\B(F)) &\cong^\otimes {\rm Comod}(\B(E))
 \\
V & \longmapsto V \square_{\B(E)} \B(E,F) \quad 
& V & \longmapsto V \square_{\B(F)} \B(F,E)
\end{align*}
be the inverse monoidal equivalences induced by the bi-Galois objects $\B(E,F)$ and $\mathcal B(F,E)$ (see Theorem \ref{explicit}).
We shall show that $\Theta(V_E) \cong V_F$, which by Theorem \ref{cleftfibre}, will prove that $\B(E,F)$ is not a cleft $\B(E)$-Galois object. 
It is easy to check that the linear map 
\begin{align*}
 \nu_F : V_F &\longrightarrow \Theta(V_E) = V_E \square_{\B(E)} \B(E,F) \\
v_j^F &\longmapsto \sum_{i=1}^m v_i^E \otimes a_{ij}^{E,F}
\end{align*}
is $\B(F)$-colinear. We get a sequence of colinear maps
\begin{equation*}
 \begin{CD}
  V_F @>\nu_F>> \Theta(V_E) @>\Theta(\nu_E)>> \Theta\Theta'(V_F) @>\simeq>> V_F
 \end{CD}
\end{equation*}
whose composition is the identity map, as shown by the following concrete computation
\begin{align*}
 V_F & \longrightarrow  V_E \square_{\B(E)} \B(E,F)
 \longrightarrow \left(V_F \square_{\B(F)} \B(F,E)\right) \square_{\B(E)} \B(E,F)
 \cong  V_F \\
v_j^F &\longmapsto \sum_{i=1}^m v_i^E \otimes a_{ij}^{E,F} \quad \quad
 \longmapsto \sum_{i=1}^m \sum_{k=1}^n v_k^E \otimes a_{ki}^{F,E} \otimes a_{ij}^{E,F} \quad \quad   \longmapsfrom v_j^F
\end{align*}
Hence $\nu_F$ is injective, $\Theta(\nu_E)$ is surjective, and $\nu_F$ is surjective
by the symmetric argument.
\end{proof}

\begin{corollary}\label{sl}
Let $E \in \GL_n(k)$ and let $q \in k^*$ be such that
${\rm tr}(E^{-1}E^t)= -q-q^{-1}$. Then we have
a $k$-linear equivalence of monoidal categories
$${\rm Comod}(\B(E)) \cong^\otimes {\rm Comod}(\mathcal O_q(\SL_2(k))$$
\end{corollary}

\begin{proof}
 This follows from the previous corollary
and Corollary \ref{cogroumonoidal}.
\end{proof}

This result has a number of interesting consequences in characteristic zero:
\begin{enumerate}
 \item (\cite{bi1}) Any cosemisimple Hopf algebra having a (co-)representation semi-ring isomorphic to the one of $\mathcal O(\SL_2)$ is isomorphic to $\B(E)$ for some matrix
$E \in \GL_n(k)$ ($n\geq 2$)
such that the solutions of ${\rm tr}(E^{-1}E^t) =-q-q^{-1}$ are
generic (i.e. $q = \pm 1$ or $q$ is not a root of unity).
\item (\cite{bi1}) For $E \in \GL_m(k)$, $F \in \GL_n(k)$, the Hopf algebras
$\B(E)$ and $\B(F)$ are isomorphic if and only if $m=n$ and there exists
$P \in \GL_n(k)$ such that $F = PEP^t$ (i.e. the bilinear forms associated to $E$ and $F$ are equivalent, by \cite{rie} this is equivalent so say that
the matrices $E^{-1}E^t$ and $F^{-1}F^t$ are conjugate).
\item (\cite{binyjm}) For any $m \geq 1$, there exists cosemisimple Hopf algebras having
an antipode of order $2m$ .
\end{enumerate}

We shall also see in Section 4 that one can deduce easily   the classification
of Galois objects over $\B(E)$ from these results.

\subsection{The universal cosovereign cogroupoid $\mathcal H$}
We present in this subsection another example of cogroupoid involving
non cleft Galois objects. This is also an occasion to advertize on 
an interesting but not very well known class of Hopf algebras.

let $F \in \GL_n(k)$.
The algebra $\mathcal H(F)$ \cite{bi01} is defined to be the universal algebra
with generators
 $(u_{ij})_{1 \leq i,j \leq n}$,
 $(v_{ij})_{1 \leq i,j \leq n}$ and relations:
$$ u v^t = v^t u = I_n \quad ; \quad vF u^t F^{-1} = 
{F} u^t F^{-1}v = I_n,$$
where $u= (u_{ij})$, $v = (v_{ij})$ and $I_n$ is
the identity $n \times n$ matrix. The algebra $\mathcal H(F)$ has a Hopf algebra
structure defined by 
$$\Delta(u_{ij}) = \sum_k u_{ik} \otimes u_{kj}, \ 
\Delta(v_{ij}) = \sum_k v_{ik} \otimes v_{kj}, \
\varepsilon (u_{ij}) = \varepsilon (v_{ij}) = \delta_{ij}, \
S(u) = {^t \! v}, \ S(v) = F {^t \! u} F^{-1}$$
Furthermore, $H(F)$ is a cosovereign Hopf algebra 
(see \cite{bi01} for the precise meaning): in particular this means
that any finite-dimensional $\mathcal H(F)$-comodule is isomorphic to its bidual
(let us say that a Hopf algebra having this property is \textbf{coreflexive}). The Hopf algebras 
$\mathcal H(F)$ are called the universal cosovereign Hopf algebras in \cite{bi01} because they have the following universal property.

\begin{proposition}
Let $H$ be a Hopf algebra and let $V$ be a finite dimensional
$H$-comodule isomorphic to its bidual comodule $V^{**}$. 
Then there exists a matrix $F \in \GL_n(k)$ ($n = \dim V$) such that
$V$ is an $\mathcal H(F)$-comodule and such that there exists a unique Hopf algebra
morphism $\pi : \mathcal H(F) \longrightarrow H$ satisfying
$(1_V \otimes \pi) \circ \beta_V = \alpha_V$,
where $\alpha_V : V \longrightarrow V \otimes H$ and
$\beta_V : V \longrightarrow V \otimes \mathcal H(F)$ denote the coactions 
of $H$ and $H(F)$ on $V$ respectively. 
$$\xymatrix{
V \ar[rr]^{\beta_V} \ar[dr]^{\alpha_V} & & V \otimes \mathcal H(F)  \ar@{-->}[dl]_{1_V \otimes f} \\
& V \otimes H &
}$$
In particular, 
every  finitely generated coreflexive Hopf algebra is a homomorphic quotient
of some Hopf algebra $\mathcal H(F)$.
\end{proposition}

\begin{proof}
Let $e_1, \ldots ,e_n$ be a basis of $V$ and let $x_{ij}$, $1 \leq i,j \leq n$ be elements of $H$ such that
$\alpha_V(e_i) = \sum_k e_k \otimes x_{ki}$, $\forall i$.
Put $y_{ij} = S(x_{ji})$. Then we have $xy^t=I_n=y^tx$ for the matrices
$x= (x_{ij})$ and $y =(y_{ij})$ since $H$ is a Hopf algebra.
Let $f : V \longrightarrow V^{**}$ be an $H$-colinear isomorphism, with
$f(e_i) = \sum_k \lambda_{ki}e_k$, for $M=(\lambda_{ij}) \in M_n(k)$.
The $H$-colinearity of $f$ means that $S^2(x)M = Mx$, i.e. 
$S(y)^t = S^2(x) = MxM^{-1}$. We also have $yS(y) =I_n=S(y)y$ since $H$ is a Hopf algebra, so we get $yM^{-1^t}x^tM^t = I_n=M^{-1^t}x^tM^ty$. Hence we get a Hopf algebra map $\pi : H(M^{-1^t}) \longrightarrow H$ such $\pi(u_{ij})=x_{ij}$ and $\pi(v_{ij})=y_{ij}$. An $H(M^{-1^t})$-comodule structure on $V$ is defined by letting
$\beta_V(e_i) = \sum_k e_k \otimes u_{ki}$, and it is clear that
$\pi$ satisfies the property in the statement, while uniqueness is clear. 
The last assertion follows from the previous one and the fact that
a Hopf algebra that is finitely generated as an algebra is generated (as a Hopf algebra) by the coefficients of one finite-dimensional comodule.
\end{proof}

The above universal property indicates that the Hopf algebras
$\mathcal H(F)$ are the quantum analogues of $\mathcal O(\GL_n(k))$
(as soon as we believe that a finite-dimensional representation
of a quantum group should be isomorphic with its bidual).

Similarly to the previous subsection, we describe $\mathcal H(F)$ as part of a cogroupoid, and we begin with a generalization involving two matrices.

Let $E \in \GL_m(k)$ and let $F \in \GL_n(k)$. The
 algebra $\mathcal H(E,F)$ is the algebra presented by generators
$u_{ij}$, $ v_{ij}$, $1\leq i \leq m, 1\leq j \leq n$,
and submitted to the relations
$$u v^t = I_m = v F u^t E^{-1} \quad ; \quad
v^tu = I_n = F u^t E^{-1} v.$$
When $E = F$, we have $\mathcal H(F,F) = \mathcal H(F)$.

The structural morphisms of the corresponding cogroupoid are constructed in the following lemma.

\begin{lemma}
\begin{enumerate}
 \item For any $E\in \GL_m(k)$, $F\in \GL_n(k)$, $G \in \GL_p(k)$, 
there exists an algebra map
\begin{align*}
\Delta_{E,F}^G : \mathcal H(E,F) &
\longrightarrow \mathcal H(E,G) \otimes \mathcal H(G,F)\\
u_{ij},  \ v_{ij} &\longmapsto \sum_{k=1}^p u_{ik} \otimes u_{kj}, \ \sum_{k=1}^p v_{ik} \otimes v_{kj}
\end{align*}
and for any $M \in \GL_r(k)$, the following diagrams commute
\begin{equation*}
\begin{CD}
\mathcal H(E,F) @>\Delta_{E,F}^G>> \mathcal H(E,G) \otimes \mathcal H(G,F) \\
@V\Delta_{E,F}^MVV
@V\Delta_{E,G}^M \otimes 1VV\\
\mathcal H(E,M) \otimes \mathcal H(M,F)  @>1 \otimes \Delta_{M,F}^G>> 
\mathcal H(E,M)\otimes \mathcal H(M,G)  \otimes \mathcal H(G,F)
\end{CD}
\end{equation*}
$$\xymatrix{
\mathcal H(E,F) \ar[d]^{\Delta_{E,F}^F} \ar@{=}[rrd]^{\rm id} \\
 \mathcal H(E,F) \otimes \mathcal H(F,F) \ar[rr]^{1 \otimes \varepsilon_F} && \mathcal H(E,F)}
\quad
\xymatrix{
\mathcal H(E,F) \ar[d]^{\Delta_{E,F}^E} \ar@{=}[rrd]^{{\rm id}} \\
 \mathcal H(E,E) \otimes \mathcal H(E,F) \ar[rr]^{\varepsilon_E \otimes 1} && \mathcal H(E,F)}
$$

\item For any $E\in \GL_m(k)$, $F\in \GL_n(k)$, there exists an algebra map 
\begin{align*}
S_{E,F} :\mathcal H(E,F) & \longrightarrow \mathcal H(F,E)^{\rm op} \\
u, \ v & \longmapsto v^t, \ E u^t F^{-1}
\end{align*}
such that the following diagrams commute
{\footnotesize
$$\xymatrix{
\mathcal H(E,E) \ar[r]^{\varepsilon_E} \ar[d]^{\Delta_{E,E}^F} &
k \ar[r]^{u}  & \mathcal H(E,F) \\
\mathcal H(E,F) \otimes \mathcal H(F,E) \ar[rr]^{1 \otimes S_{F,E}} && \mathcal H(E,F) \otimes \mathcal H(E,F) \ar[u]^{m}} \quad \quad
\xymatrix{
\mathcal H(E,E) \ar[r]^{\varepsilon_E} \ar[d]^{\Delta_{E,E}^F} &
k \ar[r]^{u}  & \mathcal H(F,E) \\
\mathcal H(E,F) \otimes \mathcal H(F,E) \ar[rr]^{S_{E,F} \otimes 1} && \mathcal H(F,E) \otimes \mathcal H(F,E) \ar[u]^{m}}
$$}
\end{enumerate}
\end{lemma}

\begin{proof}
 The proof is similar to that of Lemma \ref{strucB}.
\end{proof}

Hence the lemma ensures that we have a cogroupoid. 

\begin{definition}
 The cogroupoid $\mathcal H$ is the cogroupoid defined as follows:
\begin{enumerate}
 \item ${\rm ob}(\mathcal H) = \{ F \in \GL_n(k), \ n \geq 1\}$,
\item For $E,F \in {\rm ob}(\mathcal H)$, the algebra $\mathcal H(E,F)$ is the algebra defined above,
\item The structural maps $\Delta_{\bullet,\bullet}^\bullet$, $\varepsilon_\bullet$
and $S_{\bullet,\bullet}$ are defined in the previous lemmma.
\end{enumerate}
\end{definition}

So we have cogroupoid linking all the Hopf algebras $\mathcal H(F)$, and the natural next
question is to study the connectedness of $\mathcal H$.

\begin{lemma}
 Let $E \in \GL_m(k)$, $F \in \GL_n(k)$ with $m,n\geq 2$. Then $H(E,F) \not =(0)$ if and only 
if ${\rm tr}(E)={\rm tr}(F)$ and ${\rm tr}(E^{-1})={\rm tr}(F^{-1})$.
\end{lemma}

\begin{proof}
 It is left as an exercise to check that if $H(E,F) \not =(0)$ then 
${\rm tr}(E)={\rm tr}(F)$ and ${\rm tr}(E^{-1})={\rm tr}(F^{-1})$. The proof of the converse
 uses the diamond lemma, it might be found
in \cite{bilms}.
\end{proof}

\begin{corollary}\label{noncleft2}
 Let $\lambda , \mu \in k$. Consider the full subcogroupoid $\mathcal H^{\lambda,\mu}$ of $\mathcal H$ with objects
$${\rm ob}(\mathcal H^{\lambda,\mu})= \{ F \in \GL_n(k), \ n\geq 2,  \ {\rm tr}(F) =\lambda, \
{\rm tr}(F^{-1})=\mu\}$$
Then $\mathcal H^{\lambda, \mu}$ is a connected cogroupoid.

In particular for $E \in GL_m(k)$, $F \in \GL_n(k)$ with $n,m \geq 2$, 
${\rm tr}(E) ={\rm tr}(F)$ and
${\rm tr}(E^{-1})={\rm tr}(F^{-1})$, then
$\mathcal H(E,F)$ is a $\mathcal H(E)$-$\mathcal H(F)$-bi-Galois object and is not cleft
if $n \not = m$.
\end{corollary}

\begin{proof}
 The proof is a copy and paste of the proof of Corollary \ref{noncleft}.
\end{proof}

We can use this result to deduce the
corepresentation theory of $\mathcal H(F)$ when $F$ is a generic matrix.  
Let us introduce some notation
and terminology.  

\noindent
$\bullet$
Let $F \in \GL_n(k)$. We say that $F$ is normalized if ${\rm tr}(F) =
{\rm tr} (F^{-1})$. We say that $F$ is normalizable if there exists 
$\lambda \in k^*$ such that 
${\rm tr}(\lambda F) = {\rm tr} ((\lambda F)^{-1})$. Over an algebraically 
closed field, any matrix is normalizable unless 
${\rm tr}(F) = 0 \not = {\rm tr} (F^{-1})$  or
${\rm tr}(F)\not = 0 = {\rm tr} (F^{-1})$.
The study of the Hopf algebra $H(F)$ for a normalizable matrix
reduces to the case when $F$ is normalized, since $\mathcal H(\lambda F) = \mathcal H(F)$.

\noindent
$\bullet$ Let $q \in k^*$. As usual, we say that $q$ is generic
if $q$ is not a root of unity of order $N \geq 3$.
We say that a matrix $F \in \GL_n(k)$ is generic if $F$ is normalized
and if the solutions of $q^2 -{\rm tr}(F)q +1 = 0$ are generic. 

\noindent 
$\bullet$ 
Let $q\in k^*$. We put $F_q = \left(\begin{array}{cc} q^{-1} & 0 \\
                          0 & q \\
       \end{array} \right) \in \GL_2(k)$.
The Hopf algebra $\mathcal H(F_q)$ is denoted by $\mathcal H(q)$.

\noindent
$\bullet$ Let $F \in \GL_n(k)$. The natural $n$-dimensional 
$H(F)$-comodules associated to the multiplicative matrices $u = (u_{ij})$
and $v = (v_{ij})$ are denoted by $U$ and $V$, with $V = U^*$.

\noindent
$\bullet$ We will consider the coproduct monoid $\mathbb N * \mathbb N$.
Equivalently $\mathbb N * \mathbb N$ is the free monoid on two generators,
which we denote, by $\alpha$ and $\beta$.
There is a unique antimultiplicative morphism
$^- : \mathbb N * \mathbb N \longrightarrow \mathbb N * \mathbb N$
such that $\bar{e} = e$, $\bar{\alpha} = \beta$ and $\bar{\beta} = \alpha$
($e$ denotes the unit element of $\mathbb N * \mathbb N$).

\medskip

The corepresensation theory of $\mathcal H(F)$ is described
in the following result from \cite{bilms}.
Here $k$ denotes an algebraically
closed field.

\begin{theorem}
Let $F \in \GL_n(k)$ ($n \geq 2$) be a normalized matrix.

\noindent
(a) Let  $q \in k^*$ be such that $q^2 -{\rm tr}(F)q +1 = 0$. Then we have
a $k$-linear equivalence of monoidal categories
$${\rm Comod}(H(F)) \cong^\otimes {\rm Comod}(H(q))$$
We assume now that $k$ is a characteristic zero field.

\noindent
(b) The Hopf algebra $H(F)$ is cosemisimple if and only if $F$ is 
a generic matrix.

\noindent
(c) Assume that $F$ is generic. 
To any element $x \in \mathbb N * \mathbb N$ corresponds a simple
$\mathcal H(F)$-comodule $U_x$, with $U_e = k$, $U_\alpha = U$ and $U_\beta = V$.
Any simple $\mathcal H(F)$-comodule is isomorphic to one of the $U_x$, and
$U_x \cong U_y$ if and only if $x=y$. For $x,y \in \mathbb N * \mathbb N$, 
we have $U_x^* \cong U_{\bar{x}}$ and
$$U_x \otimes U_y  \cong \bigoplus_{\{a,b,g \in 
\mathbb N * \mathbb N|x=ag,y={\bar g}b\}} U_{ab} \ .$$ 
\end{theorem}

Part (a) follows from Corollary \ref{noncleft2} and
Corollary \ref{cogroumonoidal} and reduces the study to the case of $\mathcal H(q)$.
One then constructs a Hopf algebra embedding
$\mathcal H(q) \subset k[z,z^{-1}]* \mathcal O_q(\SL_2(k))$
and uses properties of $\mathcal O_q(\SL_2(k))$ and results
on free products of cosemisimple Hopf algebras by Wang \cite{wa95}.
 See \cite{bilms} for details.  
The fusion rule formulas arose first in \cite{ba}, where another proof
of (c), for positive matrices $F \in \GL_n(\mathbb C)$, might be found.

\subsection{The 2-cocycle cogroupoid of a Hopf algebra}
We now come back to the familiar Hopf-Galois objects obtained from $2$-cocycles and see how they behave in the cogroupoid framework. 

Let $H$ be a Hopf algebra and let $\sigma, \tau \in Z^2(H)$.
The algebra $H(\sigma , \tau)$ is the algebra having $H$ as underlying vector space
and product defined by
$$x.y = \sigma(x_{(1)},y_{(1)}) \tau^{-1}(x_{(3)},y_{(3)}) x_{(2)}y_{(2)}$$
It is straighforward to check that $H(\sigma, \tau)$ is an associative algebra
(with the same unit as $H$). 

In the following lemma we define all the necessary structural maps for the 
$2$-cocycle cogroupoid of $H$.

\begin{lemma}
\begin{enumerate}
 \item Let $\sigma, \tau, \omega \in Z^2(H)$. The map
\begin{align*}
 \Delta_{\sigma, \tau}^\omega = \Delta : H(\sigma, \tau) &\longrightarrow
H(\sigma, \omega) \otimes H(\omega, \tau) \\
x &\longmapsto x_{(1)} \otimes x_{(2)}
\end{align*}
is an algebra map, and for any $\alpha \in Z^2(H)$, the following diagram commutes:
\begin{equation*}
\begin{CD}
H(\sigma,\tau) @>\Delta_{\sigma,\tau}^\omega>> 
H(\sigma,\omega) \otimes H(\omega,\tau) \\
@V\Delta_{\sigma,\tau}^{\alpha}VV
@V\Delta_{\sigma,\omega}^{\alpha}\otimes 1VV\\
H(\sigma,\alpha) \otimes H(\alpha,\tau)  @>1 \otimes \Delta_{\alpha,\tau}^\omega>> 
H(\sigma,\alpha)\otimes H(\alpha,\omega)  \otimes H(\omega,\tau)
\end{CD}
\end{equation*}
\item Let $\sigma \in Z^2(H)$. The linear map $\varepsilon_\sigma = \varepsilon : H(\sigma,\sigma) \longrightarrow k$ is an algebra map, and for any $\tau \in Z^2(H)$, the following
diagrams commute:
$$\xymatrix{
H(\sigma,\tau) \ar[d]^{\Delta_{\sigma,\tau}^\tau} \ar@{=}[rrd]^{\rm id} \\
 H(\sigma,\tau) \otimes H(\tau,\tau) \ar[rr]^{1 \otimes \varepsilon_\tau} && H(\sigma,\tau)}
\quad
\xymatrix{
H(\sigma,\tau) \ar[d]^{\Delta_{\sigma,\tau}^\tau} \ar@{=}[rrd]^{{\rm id}} \\
 H(\sigma,\sigma) \otimes H(\sigma,\tau) \ar[rr]^{\varepsilon_\sigma \otimes 1} && H(\sigma,\tau)}
$$
\item Let $\sigma, \tau \in Z^2(H)$. Consider the linear map
\begin{align*}
 S_{\sigma, \tau} : H(\sigma , \tau) &\longrightarrow H(\tau , \sigma) \\
x &\longmapsto \sigma(x_{(1)}, S(x_{(2)})) \tau^{-1}(S(x_{(4)}), x_{(5)})S(x_{(3)})
\end{align*}
The following diagrams commute:
{\footnotesize
$$\xymatrix{
H(\sigma,\sigma) \ar[r]^{\varepsilon_\sigma} \ar[d]^{\Delta_{\sigma,\sigma}^\tau} &
k \ar[r]^{u}  & H(\sigma,\tau) \\
H(\sigma,\tau) \otimes H(\tau,\sigma) \ar[rr]^{1 \otimes S_{\tau,\sigma}} && H(\sigma,\tau) \otimes H(\sigma,\tau) \ar[u]^{m}} \quad \quad
\xymatrix{
H(\sigma,\sigma) \ar[r]^{\varepsilon_\sigma} \ar[d]^{\Delta_{\sigma,\sigma}^\tau} &
k \ar[r]^{u}  & H(\tau,\sigma) \\
H(\sigma,\tau) \otimes H(\tau,\sigma) \ar[rr]^{S_{
\sigma,\tau} \otimes 1} && H(\tau,\sigma) \otimes H(\tau,\sigma) \ar[u]^{m}}
$$}
\end{enumerate}
\end{lemma}  

\begin{proof}
 Assertions (1) and (2) are obvious. Let $x \in H$. In $H(\sigma, \tau)$, we have
\begin{align*}
x_{(1)} \cdot S_{\tau,\sigma}(x_{(2)}) &= 
x_{(1)} \cdot \left(\tau(x_{(2)},S(x_{(3)}))\sigma^{-1}(S(x_{(5)}),x_{(6)})S(x_{(4)})\right) \\
&=\tau(x_{(4)},S(x_{(5)})) \sigma^{-1}(S(x_{(9)}),x_{(10)})
\sigma(x_{(1)}, S(x_{(8)})) \tau^{-1}(x_{(3)}, S(x_{(6)}))  x_{(2)} S(x_{(7)})\\
& =\tau^{-1}*\tau(x_{(3)},S(x_{(4)}))  \sigma^{-1}(S(x_{(7)}), x_{(8)}) \sigma(x_{(1)},S(x_{(6)}))
x_{(2)} S(x_{(5)}) \\
& = \sigma^{-1}(S(x_{(5)}), x_{(6)}) \sigma(x_{(1)},S(x_{(4)}))
x_{(2)} S(x_{(3)}) \\
& = \sigma(x_{(1)}, S(x_{(2)})) \sigma^{-1}(S(x_{(3)}),x_{(4)}) =
\varepsilon(x)1
\end{align*}
where the last identity is (a5) from Theorem 1.6 in \cite{doi}.
This proves the commutativity of the first diagram in (3). For the second diagram, we have in $H(\tau, \sigma)$
\begin{align*}
 S_{\sigma, \tau}(x_{(1)}) \cdot x_{(2)} &= 
 \left(\sigma(x_{(1)},S(x_{(2)}))\tau^{-1}(S(x_{(4)}),x_{(5)})S(x_{(3)})\right)
\cdot x_{(6)} \\
& = \sigma(x_{(1)},S(x_{(2)}))\tau^{-1}(S(x_{(6)}),x_{(7)})
\tau(S(x_{(5)}),x_{(8)})\sigma^{-1}(S(x_{(3)}),x_{(10)})S(x_{(4)})x_{(9)} \\
& = \sigma(x_{(1)},S(x_{(2)}))\tau^{-1}*\tau(S(x_{(5)}),x_{(6)})
\sigma^{-1}(S(x_{(3)}),x_{(8)})S(x_{(4)})x_{(7)} \\
& = \sigma(x_{(1)},S(x_{(2)}))\sigma^{-1}(S(x_{(3)}),x_{(6)})S(x_{(4)})x_{(5)} \\
& = \sigma(x_{(1)},S(x_{(2)})) \sigma^{-1}(S(x_{(3)}),x_{(4)}) = 
\varepsilon(x)1
\end{align*}
where again the last identity is (a5) from Theorem 1.6 in \cite{doi}.
\end{proof}

\begin{definition}
 Let $H$ be a Hopf algebra. The \textbf{$2$-cocycle cogroupoid of $H$}, denoted 
$\underline{H}$, is the cogroupoid defined as follows:
\begin{enumerate}
 \item ${\rm ob}(\underline{H}) = Z^2(H)$.
\item For $\sigma, \tau \in Z^2(H)$, the algebra $\underline{H}(\sigma,\tau)$
is the algebra $H(\sigma,\tau)$ defined above.
\item The structural maps $\Delta_{\bullet,\bullet}^\bullet$, $\varepsilon_\bullet$
and $S_{\bullet,\bullet}$ are the ones defined in the previous lemmma.
\end{enumerate}
\end{definition}
 
The Hopf algebra $H(\sigma,\sigma)$ is the Hopf algebra $H^\sigma$ defined by Doi in \cite{doi}. The algebras $_\sigma H$ and $H_{\sigma^{-1}}$ considered in Example 1.3
are the algebras $H(\sigma,1)$ and $H(1,\sigma)$ respectively (where $1$ stand for $\varepsilon \otimes \varepsilon$). It follows that $H(\sigma,1)$ is an 
$H^\sigma$-$H$-bi-Galois object and that $H(1,\sigma)$ is an $H$-$H^\sigma$-bi-Galois object. The associated monoidal equivalence
$${\rm Comod}(H) \cong^\otimes {\rm Comod}(H^\sigma)$$
is isomorphic (as a monoidal functor) to the monoidal
equivalence whose underlying functor is the identity (recall that $H=H^\sigma$ as coalgebras) and
whose monoidal constraints is given by the following
isomorphisms
\begin{align*}
 V \otimes W &\cong V \otimes W \\
v \otimes w &\mapsto \sigma^{-1}(v_{(1)},w_{(1)}) v_{(0)} \otimes w_{(0)}
\end{align*}

\subsection{The multiparametric $\GL_n$-cogroupoid}
It is in general unpleasant and difficult to work with
explicit cocycles, and
the 2-cocycle cogroupoid in the previous subsection
is a theoretical tool. 
In concrete examples, it is much easier to work with explicit algebras
(although it is useful to know that a cocycle lies behind 
the constructions, e.g. to prove that the algebras are non-zero),
and in this subsection we present
an example of a  full subcogroupoid of the 2-cocycle cogroupoid
of $\mathcal O(\GL_n(k))$. 

We say that a matrix ${\bf p}=(p_{ij}) \in M_n(k)$ is an AST-matrix
(after Artin-Schelter-Tate \cite{ast}) if 
$p_{ii} = 1$ and $p_{ij}p_{ji} = 1$ for all $i$ and $j$.
The trivial AST-matrix (i.e. $p_{ij} =1$ for all $i$ and $j$)
is denoted by ${\bf 1}$. 
We denote by ${\rm AST}(n)$ the set of AST matrices of size $n$ .

For ${\bf p} \in {\rm AST}(n)$, the algebra $\mathcal O_{\mathbf p}({\rm GL}_n(k))$ (\cite{ast}) is the algebra presented by generators
$x_{ij}$, $y_{ij}$, $ 1\leq i,j\leq n$, submitted to the following relations ($1 \leq i,j,k,l \leq n$):
$$x_{kl}x_{ij} = p_{ki}p_{jl}x_{ij}x_{kl}, \quad y_{kl}y_{ij} = p_{ki}p_{jl}y_{ij}y_{kl}, \quad
y_{kl}x_{ij} = p_{ik}p_{lj}x_{ij}y_{kl}$$
$$\sum_{k=1}^n x_{ik}y_{jk} = \delta_{ij} = \sum_{k=1}^n x_{ki}y_{kj}$$ 
The presentation we have given avoids the use of the quantum determinant.
The algebra $ \mathcal O_{\mathbf p}({\rm GL}_n(k))$ has a standard Hopf algebra structure, described
as follows:
$$\Delta(x_{ij})=\sum_k x_{ik} \otimes x_{kj}, \ \Delta(y_{ij})=\sum_k y_{ik} \otimes y_{kj}, \
\varepsilon(x_{ij}) =\delta_{ij} = \varepsilon(y_{ij}), \ S(x_{ij}) = y_{ji}, \ S(y_{ij})= x_{ji}$$
When ${\bf p}={\bf 1}$, one gets the usual Hopf algebra
 $\mathcal O({\rm GL}_n(k))$.

To put $\mathcal O_{\mathbf p}({\rm GL}_n(k))$ in a cogroupoid framework, we use 
the following algebras. For ${\bf p}, {\bf q} \in {\rm AST}(n)$, the algebra
$\mathcal O_{\mathbf p, \mathbf q }({\rm GL}_n(k))$ is the algebra presented by generators
$x_{ij}$, $y_{ij}$, $ 1\leq i,j\leq n$, submitted to the following relations ($1 \leq i,j,k,l \leq n$):
$$x_{kl}x_{ij} = p_{ki}q_{jl}x_{ij}x_{kl}, \quad y_{kl}y_{ij} = p_{ki}q_{jl}y_{ij}y_{kl}, \quad
y_{kl}x_{ij} = p_{ik}q_{lj}x_{ij}y_{kl}$$
$$\sum_{k=1}^n x_{ik}y_{jk} = \delta_{ij} = \sum_{k=1}^n x_{ki}y_{kj}$$ 

\begin{lemma} In this lemma we note $\mathcal O_{\mathbf p, \bf q}=\mathcal O_{\mathbf p, \bf q}({\rm GL}_n(k))$.

\begin{enumerate}
 \item For any ${\bf p}, {\bf q}, {\bf r} \in {\rm AST}(n)$, 
there exists an algebra map
\begin{align*}
\Delta_{{\bf p},{\bf q}}^{\bf r} : \mathcal O_{\mathbf p, \bf q} &
\longrightarrow  \mathcal O_{\mathbf p, \bf r} \otimes \mathcal O_{\mathbf r, \bf p}\\
x_{ij},  \ y_{ij} &\longmapsto \sum_{k=1}^n x_{ik} \otimes x_{kj}, \ \sum_{k=1}^n y_{ik} \otimes y_{kj}
\end{align*}
and for any $\bf s \in {\rm AST}(n)$, the following diagrams commute
\begin{equation*}
\begin{CD}
\mathcal O_{\mathbf p,\bf q} @>\Delta_{\bf p,\bf q}^{\bf r}>> 
\mathcal O_{\mathbf p, \bf r} \otimes \mathcal O_{\bf r, \mathbf p} \\
@V\Delta_{\bf p,\bf q}^{\bf s}VV
@V\Delta_{\bf p, \bf r}^{\bf s} \otimes 1VV\\
\mathcal O_{\mathbf p, \bf s} \otimes  \mathcal O_{\bf s ,\mathbf q} @>1 \otimes \Delta_{\bf s,\bf q}^{\bf r}>> 
\mathcal O_{\mathbf p, \bf s}\otimes \mathcal O_{\mathbf s, \bf r}  \otimes \mathcal O_{\mathbf r, \bf q}
\end{CD}
\end{equation*}
$$\xymatrix{
\mathcal O_{\mathbf p, \bf q} \ar[d]^{\Delta_{\bf p, \bf q}^{\bf q}} \ar@{=}[rrd]^{\rm id} \\
 \mathcal O_{\mathbf p, \bf q} \otimes \mathcal O_{\mathbf q, \bf q} \ar[rr]^{1 \otimes \varepsilon_{\bf q}} && \mathcal O_{\mathbf p, \bf q}}
\quad
\xymatrix{
\mathcal O_{\mathbf p, \bf q} \ar[d]^{\Delta_{\bf p,\bf q}^{\bf p}} \ar@{=}[rrd]^{{\rm id}} \\
 \mathcal O_{\mathbf p, \bf p} \otimes \mathcal O_{\mathbf p, \bf q} \ar[rr]^{\varepsilon_{\bf p} \otimes 1}&& \mathcal O_{\mathbf p, \bf q}}
$$

\item For any $ \bf p, \bf q \in {\rm AST}(n)$, there exists an algebra map 
\begin{align*}
S_{\bf p, \bf q} :\mathcal O_{\mathbf p, \bf q}({\rm GL}_n(k)) & \longrightarrow \mathcal O_{\bf q , \mathbf p}({\rm GL}_n(k))^{\rm op} \\
x, \ y & \longmapsto y^t, \  x^t 
\end{align*}
such that the following diagrams commute
{\footnotesize
$$\xymatrix{
\mathcal O_{\mathbf p,\bf p} \ar[r]^{\varepsilon_{\bf p}} \ar[d]^{\Delta_{\bf p,\bf p}^{\bf q}} &
k \ar[r]^{u}  & \mathcal O_{\mathbf p, \bf q} \\
\mathcal O_{\mathbf p, \bf q} \otimes \mathcal O_{\bf q, \mathbf p} \ar[rr]^{1 \otimes S_{\bf q, \bf p}} && \mathcal O_{\mathbf p, \bf q} \otimes \mathcal O_{\mathbf p, \bf q} \ar[u]^{m}} \quad \quad
\xymatrix{
\mathcal O_{\mathbf p,\bf p} \ar[r]^{\varepsilon_{\bf p}} \ar[d]^{\Delta_{\bf p,\bf p}^{\bf q}} &
k \ar[r]^{u}  & \mathcal O_{\mathbf q,\bf p} \\
\mathcal O_{\mathbf p,\bf q} \otimes \mathcal O_{\mathbf q,\bf p} \ar[rr]^{S_{\bf p,\bf q} \otimes 1} && \mathcal O_{\mathbf q,\bf p} \otimes \mathcal O_{\mathbf q,\bf p} \ar[u]^{m}}
$$}
\end{enumerate}
\end{lemma}

\begin{proof}
 Exercise.
\end{proof}

The lemma ensures that we have a cogroupoid. 

\begin{definition}
 The multiparametric $\GL_n$-cogroupoid, denoted $\underline{\GL_n}$, is defined 
as follows:
\begin{enumerate}
 \item ${\rm ob}(\underline{\GL_n}) = {\rm AST}(n)$,
\item For $\bf p, \bf q \in {\rm AST}(n)$, the algebra  $\underline{\GL_n}(\bf p, \bf q)$ is the algebra $\mathcal O_{\mathbf p, \mathbf q }({\rm GL}_n(k))$
 defined above.
\item The structural maps $\Delta_{\bullet,\bullet}^\bullet$, $\varepsilon_\bullet$
and $S_{\bullet,\bullet}$ are defined in the previous lemmma.
\end{enumerate}
\end{definition}

\begin{proposition}
 The multiparametric $\GL_n$-cogroupoid is connected.
\end{proposition}

\begin{proof}
 We know from Proposition \ref{connec} that it is enough to show that
$\mathcal O_{\mathbf p, \mathbf 1 }({\rm GL}_n(k)) \not = (0)$ 
for any ${\bf p} \in {\rm AST}(n)$. Consider
the algebra $k_{\bf p}[t_1^{\pm 1}, \cdots , t_n^{\pm 1}]$
presented by generators $t_1, \cdots , t_n, t_1^{- 1}, \cdots , t_n^{-1}$
submitted to the relations ($1 \leq i,k \leq n$)
$$t_it_i^{-1} = 1 = t_{i}^{-1} t_i, \ t_it_k = p_{ik}t_k t_i, \ t_i^{-1}t_k^{-1} = p_{ik}t_k^{-1} t_i^{-1}, \ t_it_k^{-1} = p_{ki}t_k^{-1}t_i$$
It is not difficult to check that $k_{\bf p}[t_1^{\pm 1}, \cdots , t_n^{\pm 1}]$
is a non-zero algebra either by using the diamond lemma or by showing
that it is isomorphic to the twisted group algebra $k_\sigma[\mathbb Z^n]$ for
the $2$-cocycle considered in \cite{ast}. It is straighforward to check
that we have a surjective algebra map
\begin{align*}
 \mathcal O_{\mathbf p, \mathbf 1 }({\rm GL}_n(k))&\longrightarrow
k_{\bf p}[t_1^{\pm 1}, \cdots , t_n^{\pm 1}] \\
x_{ij}, y_{ij} &\longmapsto \delta_{ij} t_i, \delta_{ij}t_i^{-1}
\end{align*}
and hence $\mathcal O_{\mathbf p, \mathbf 1 }({\rm GL}_n(k)) \not = (0)$ .
\end{proof}

 \begin{corollary}\label{multimonoidal}
For any $\bf p \in {\rm AST}(n)$, we have
a $k$-linear equivalence of monoidal categories
$${\rm Comod}(\mathcal O_{\mathbf p}({\rm GL}_n(k))) \cong^\otimes {\rm Comod}(\mathcal O(\GL_n(k))$$
\end{corollary}

\section{Classifying Hopf-Galois objects: the fibre functor method}

In this section we explain how to classify Hopf-Galois
objects by using fibre functors via Ulbrich's Theorem.
The idea is to study how the fibre functor associated to the Galois
object will transform the ``fundamental'' morphisms of the category
of comodules.

\subsection{Hopf-Galois objects over $\B(E)$}
Let $E \in \GL_m(k)$ with $m \geq 2$. We first state the classification result
for (left) $\B(E)$-Galois objects. We already know from
the previous section that $\B(E,F)$ is a left $\B(E)$-Galois object
if $F \in \GL_n(k)$ satisfies ${\rm tr}(E^{-1}E^t) = {\rm tr}(F^{-1}F^t)$.
Let 
$$X_0(E) = \{F \in \GL_n(k),  \ n \geq 2 \ | \ {\rm tr}(E^{-1}E^t) = {\rm tr}(F^{-1}F^t)\}$$
and let $\sim$ be the equivalence relation on 
$X_0(E)$ defined for $F,G \in X_0(E)$ with $F \in \GL_n(k)$ and
$G\in \GL_p(k)$ by 
$$F \sim G \iff [n=p \ {\rm and} \ \exists P \in \GL_n(k) \ {\rm with} 
\ F = PGP^t]$$
We put $X(E)=X_0(E)/\sim$.
The following result is stated in \cite{au}.
\begin{theorem}\label{classBE}
 Let $E \in \GL_m(k)$ with $m \geq 2$. The map
\begin{align*}
 X_0(E) &\longrightarrow {\rm Gal}(\B(E)) \\
F &\longmapsto [\B(E,F)]
\end{align*}
induces a bijection $X(E)\simeq {\rm Gal}(\B(E))$.
\end{theorem}

Of course the meaning of bracket symbol $[,]$ is that 
we consider the isomorphism class of the Galois object.

let $V_E$ be the $m$-dimensional $\B(E)$-comodule
with basis $v_1^E, \ldots , v_m^E$ and with $\B(E)$-coaction
$\alpha(v_i^E)= \sum_j v_j^E \otimes a_{ji}$. 
We know that the nondegenerate bilinear form
$$\beta_E :V_E \otimes V_E \longmapsto k, \
v_i^E \otimes v_j^E \longmapsto \lambda_{ij}, \ E= (\lambda_{ij})$$
is $\B(E)$-colinear and generates ${\rm Comod}(\B(E))$ as a tensor
category (by the universal property of $\B(E)$).
So to study the fibre functors on ${\rm Comod}(\B(E))$
 we have to study how they transform $\beta_E$. We then deduce
informations on the corresponding Galois object.

\begin{proposition}\label{surj}
 Let $A$ be a left $\B(E)$-Galois object. Then there exists 
$F \in X_0(E)$ such that $A \cong \B(E,F)$ as Galois objects.
\end{proposition}
 
\begin{proof}
Let $\Omega^A : {\rm Comod}_f(\B(E))\longmapsto k$ be the fibre
functor corresponding to $A$ and consider the bilinear map defined by the composition
\begin{equation*}
 \begin{CD}
  \beta' : \Omega^A(V_E) \otimes \Omega^A(V_E) @>\sim>> \Omega^A(V_E\otimes V_E)
@>\Omega^A(\beta_E)>> \Omega^A(k) @>\sim>> k \\
@| @| @| @| \\
V_E\square_{\B(E)} A \otimes V_E\square_{\B(E)} A @>\sim>>
(V_E \otimes V_E)\square_{\B(E)} A @>\beta_E\otimes 1>>
k\square_{\B(E)} A @>\sim>>k
 \end{CD}
\end{equation*}
Let $w_1, \ldots ,w_n$ be a basis of $\Omega^A(V_E)$ (we have $n \geq 2$ by Proposition \ref{fini}).
For $i \in \{1, \ldots , n\}$, write 
$$w_i = \sum_{k=1}^m v_k^E \otimes z_{ki}$$
for elements $z_{ki}\in A$, $1 \leq k \leq m$, $1 \leq i \leq n$.
Let $F=(\mu_{ij}) \in M_n(k)$ be such that 
$\beta'(w_i \otimes w_j) = \mu_{ij}$. We have
 $$\beta'(w_i \otimes w_j)1 = \beta'(\sum_{k,l=1}^m v_k^E \otimes z_{ki} \otimes v_l^E \otimes z_{lj})1 = \sum_{k,l=1}^m z_{ki}\lambda_{kl}z_{lj} =
\mu_{ij}1$$
which means that if $z =(z_{ki}) \in M_{m,n}(A)$, we have
$$z^t E z = F$$
We have to check that the matrix $F$ is invertible.
We use the fact that $\beta$ gives a dual for $V_E$ in the category
of $\B(E)$-comdules. For this
consider the linear map
$$\delta : k \longrightarrow V_E \otimes V_E, 1 \longmapsto \sum_{k,l}\lambda_{kl}^{-1} v_k^E \otimes v_l^E, \ ({\rm with} \ E^{-1}=(\lambda_{kl}^{-1}))$$
This map is $\B(E)$-colinear (direct verification) and hence we can consider the linear
map $\delta'$ defined by the composition
\begin{equation*}
 \begin{CD}
  \delta' :  k @>\sim>> \Omega^A(k) @>\Omega^A(\delta)>>
\Omega^A(V_E\otimes V_E) @>(\widetilde{\Omega}^A_{V_E,V_E})^{-1}>>\Omega^A(V_E) \otimes \Omega^A(V_E)\\
@| @| @| @| \\
 k @>\sim>> k\square_{\B(E)}A @>\delta\otimes 1>>
(V_E \otimes V_E)\square_{\B(E)} A 
  @>\sim>> V_E\square_{\B(E)} A \otimes V_E\square_{\B(E)} A
 \end{CD}
\end{equation*}
We have 
$$\delta'(1) = \sum_{i,j}\gamma_{ij} w_i \otimes w_j = \sum_{k,l} \sum_{i,j}
\gamma_{ij}
v_k^E \otimes z_{ki} \otimes v_l^E \otimes v_{lj} \otimes z_{lj}$$ for
$G=(\gamma_{ij}) \in M_n(k)$ and hence
$$\widetilde{\Omega}^A_{V_E,V_E}\circ\delta'(1) =
\sum_{k,l}v_k^E \otimes v_l^E \otimes (\sum_{i,j} z_{ki}\gamma_{ij}z_{lj})
$$
On the other hand
$$\widetilde{\Omega}^A_{V_E,V_E}\circ\delta'(1) =
\sum_{k,l}v_k^E \otimes v_l^E \otimes \lambda_{kl}^{-1}1$$
which shows that $zGz^t=E^{-1}$.
We have
$$(1\otimes \beta_E) \circ (\delta \otimes 1) =1_{V_E}$$
in the monoidal category of $\B(E)$-comodules, so the monoidal functor
$\Omega^A$ transforms this equation into
$$(1\otimes \beta') \circ (\delta' \otimes 1) =1_{\Omega^A(V_E)}$$
This means that $FG=I_n$, so that $F$ is invertible with $G=F^{-1}$.
We have $z^t E z = F$ and $zGz^t=E^{-1}$, hence
$$F^{-1}z^t E z= I_n \ {\rm and} \ zF^{-1}z^tE=I_m
$$
Thus there exists a unique algebra morphism $f : \B(E,F) \longrightarrow A$
defined by $f(a_{ki})=  z_{ki}$. Moreover it is easy to see that $f$
is left $\B(E)$-colinear since
$$w_i = \sum_k v_{k}^E \otimes z_{ki} \in V_E \square_{\B(E)} A$$ 
We conclude by Proposition \ref{galgroupo} that $f$ is an isomorphism.  
\end{proof}

To classify the Galois objects $\B(E,F)$ up to isomorphism, 
it is also useful to examine the associated fibre functors.

\begin{proposition}\label{inj}
 Let $F,G \in X_0(E)$ with $F \in \GL_n(k)$ and $G \in\GL_p(k)$. The following
assertions are equivalent.
\begin{enumerate}
 \item The left $\B(E)$-comodule algebras $\B(E,F)$ and $\B(E,G)$ are isomorphic.
\item $F \sim G$, i.e. $n=p$  and $\exists P \in \GL_n(k)$ with 
 $F = PGP^t$.
\end{enumerate}
\end{proposition}

\begin{proof}
We use the notation in the previous proof.
 Let $f : \B(E,F) \longrightarrow \B(E,G)$ be a left $\B(E)$-comodule algebra isomorphism.
Then $f$ induces a linear isomorphim
$$1 \otimes f : V_E \square_{\B(E)} \B(E,F) \longrightarrow V_E \square_{\B(E)} \B(E,G)$$  
By corollary \ref{noncleft} and its proof we have $n=p$ and the elements
$$w_i =\sum_k v_k^E \otimes a_{ki}^{E,F}, \quad
w_i' =\sum_k v_k^E \otimes a_{ki}^{E,G}, \quad 1\leq i \leq n$$
form respective bases of $V_E \square_{\B(E)} \B(E,F)$ and $V_E \square_{\B(E)} \B(E,G)$. Then there exists $M = (m_{ij}) \in \GL_n(k)$ such that
$1 \otimes f(w_i) =\sum_j m_{ji} w_j'$, $\forall i$, and hence
$$\sum_k v_k^E \otimes f(a_{ki}^{E,F}) = \sum_k v_k^E \otimes
(\sum_j a_{kj}^{E,G} m_{ji})$$
which in matrix form gives $f(a^{E,F}) = a^{E,G}M$ or dropping
the exponents $f(a) = aM$. We have $F^{-1}a^tEa=I_n$ in $\B(E,F)$,
hence
$$I_n=f(F^{-1}a^tEa) = F^{-1}M^ta^tEaM\Rightarrow
MF^{-1}M^ta^tEa = I_n \ {\rm in} \ \B(E,G)$$
But we have $aG^{-1}a^tE=I_m$ so $MF^{-1}M^ta^tE = G^{-1}a^tE\Rightarrow
MF^{-1}M^ta^t = G^{-1}a^t$ and we conclude
by multiplying on the right by $EaG^{-1}$ that $MF^{-1}M^t=G^{-1}$. This proves
that $F \sim G$.

Conversely, if $n=p$ and $\exists P \in \GL_n(k)$ such that
$F=PGP^t$, it is not difficult so show that there exists a unique
algebra map $f : \B(E,F) \longrightarrow \B(E,G)$ such that 
$f(a) = aP^t$, that $f$ is an isomorphism and is $\B(E)$-colinear.
\end{proof}

Theorem \ref{classBE} follows from these two propositions.

\begin{remark}{\rm 
 There is an interesting equivalence relation on Galois objects, called
homotopy, introduced by Kassel and further developped by Kassel and Schneider
in \cite{kassch}.
The classification of Galois objects over $\B(E)$ up to homotopy is studied
in \cite{au}.
}
\end{remark}

\subsection{Hopf-Galois objects over $\mathcal H(F)$}
We provide now the classification of the Galois objects over the universal cosovereign Hopf algebras $\mathcal H(F)$, following the same line as the one of the previous subsection.

Let $E \in \GL_m(k)$ with $m \geq 2$.  We  know from
the previous section that $\mathcal H(E,F)$ is a left $\mathcal H(E)$-Galois object
if $F \in \GL_n(k)$ satisfies 
${\rm tr}(E) = {\rm tr}(F)$ and ${\rm tr}(E^{-1}) = {\rm tr}(F^{-1})$.
Let 
$$\overline{X}_0(E) = \{F \in \GL_n(k),  \ n \geq 2 \ | \ 
{\rm tr}(E) = {\rm tr}(F), \ {\rm tr}(E^{-1}) = {\rm tr}(F^{-1})\} $$
and let $\approx$ be the equivalence relation on 
$\overline{X}_0(E)$ defined for $F,G \in \overline{X}_0(E)$ with $F \in \GL_n(k)$ and
$G\in \GL_p(k)$ by 
$$F \approx G \iff [n=p \ {\rm and} \ \exists P \in \GL_n(k) \ {\rm with} 
\ F = PGP^{-1}]$$
We put $\overline{X}(E)=\overline{X}_0(E)/\approx$.

\begin{theorem}\label{classHE}
 Let $E \in \GL_m(k)$ with $m \geq 2$. The map
\begin{align*}
 \overline{X}_0(E) &\longrightarrow {\rm Gal}(\mathcal H(E)) \\
F &\longmapsto [\mathcal H (E,F)]
\end{align*}
induces a bijection $\overline{X}(E)\simeq {\rm Gal}(\mathcal H(E))$.
\end{theorem}

\begin{proposition}\label{surj2}
 Let $A$ be a left $\mathcal H(E)$-Galois object. Then there exists 
$F \in \overline{X}_0(E)$ such that $A \cong \mathcal H(E,F)$ as Galois objects.
\end{proposition}

\begin{proof}
The reasoning is similar to the one of Proposition \ref{surj} so we give less detail. Let $U_E$ be the $m$-dimensional $\mathcal H(E)$-comodule
with basis $u_1^E, \ldots , u_m^E$ and with $\mathcal H(E)$-coaction given by
$\alpha(u_i^E)= \sum_j u_j^E \otimes u_{ji}$. The coaction of the dual comodule
$U_E^*$ is given by $\alpha^*(u_i^{E*})= \sum_j u_j^{E*} \otimes v_{ji}$.
The following linear maps 
$$ e : U_E^* \otimes U_E \longrightarrow k, \ u_i^* \otimes u_j \longmapsto \delta_{ij}$$
$$d : k \longrightarrow U_E \otimes U_E^*, \ 1 \longmapsto \sum_{i=1}^m u_i \otimes u_i^* $$
$$ \epsilon :U_E \otimes U_E^* \longrightarrow k, \ u_i \otimes u_j^* \longmapsto \lambda_{ij}^{-1}$$
$$ \delta : k \longrightarrow U_E^* \otimes U_E, \ 1 \longmapsto \sum_{i,j=1}^m \lambda_{ij} u_i^* \otimes u_j$$
where $E= (\lambda_{ij})$ and $E^{-1}=(\lambda_{ij}^{-1})$, are $\mathcal H(E)$-colinear.
 Let $\Omega^A : {\rm Comod}_f(\mathcal H(E))\longmapsto k$ be the fibre
functor corresponding to $A$.
Consider the linear maps defined by the compositions
\begin{equation*}
 \begin{CD}
  e' : \Omega^A(U_E^*) \otimes \Omega^A(U_E) @>\sim>> \Omega^A(U_E^*\otimes U_E)
@>\Omega^A(e)>> \Omega^A(k) @>\sim>> k \\
@| @| @| @| \\
U_E^*\square_{\mathcal H(E)} A \otimes U_E\square_{\mathcal H(E)} A @>\sim>>
(U_E^* \otimes U_E)\square_{\mathcal H(E)} A @>e\otimes 1>>
k\square_{\mathcal H(E)} A @>\sim>>k
 \end{CD}
\end{equation*}
\begin{equation*}
 \begin{CD}
  d' :  k @>\sim>> \Omega^A(k) @>\Omega^A(d)>>
\Omega^A(U_E\otimes U_E^*) @>(\widetilde{\Omega}^A_{U_E,U_E^*})^{-1}>>\Omega^A(U_E) \otimes \Omega^A(U_E^*)\\
@| @| @| @| \\
 k @>\sim>> k\square_{\mathcal H(E)}A @>d\otimes 1>>
(U_E \otimes U_E^*)\square_{\mathcal H(E)} A 
  @>\sim>> U_E\square_{\mathcal H(E)} A \otimes U_E^*\square_{\mathcal H(E)} A
 \end{CD}
\end{equation*}
Then $(\Omega^A(U_E^*),e',d')$ is a left dual for $\Omega^A(U_E)$ since $\Omega^A$ is a monoidal functor. Hence there exist bases $w_1, \ldots , w_n$ and $w'_1, \ldots , w'_n$ of $\Omega^A(U_E) $ and $\Omega^A(U_E^*)$ respectively ($n=\dim\Omega^A(U_E)\geq 2$) such that $e'(w'_i\otimes w_j) = \delta_{ij}$ and $d'(1)=\sum_i w_i \otimes w_i'$.
 Now write $w_i= \sum_k u_k^E \otimes z_{ki}$ and $w'_i= \sum_k u_k^{E*} \otimes c_{ki}$ for elements $z_{ki}, c_{ki} \in A$. Then we see that $c^tz= I_n$ and
$zc^t=I_m$.

Consider now the linear maps defined by the compositions
\begin{equation*}
 \begin{CD}
  \epsilon' : \Omega^A(U_E) \otimes \Omega^A(U_E^*) @>\sim>> \Omega^A(U_E\otimes U_E^*)
@>\Omega^A(\epsilon)>> \Omega^A(k) @>\sim>> k \\
@| @| @| @| \\
U_E\square_{\mathcal H(E)} A \otimes U_E^*\square_{\mathcal H(E)} A @>\sim>>
(U_E \otimes U_E^*)\square_{\mathcal H(E)} A @>\epsilon\otimes 1>>
k\square_{\mathcal H(E)} A @>\sim>>k
 \end{CD}
\end{equation*}
\begin{equation*}
 \begin{CD}
  \delta' :  k @>\sim>> \Omega^A(k) @>\Omega^A(\delta)>>
\Omega^A(U_E^*\otimes U_E) @>(\widetilde{\Omega}^A_{U_E^*,U_E})^{-1}>>\Omega^A(U_E^*) \otimes \Omega^A(U_E)\\
@| @| @| @| \\
 k @>\sim>> k\square_{\mathcal H(E)}A @>\delta\otimes 1>>
(U_E^* \otimes U_E)\square_{\mathcal H(E)} A 
  @>\sim>> U_E^*\square_{\mathcal H(E)} A \otimes U_E\square_{\mathcal H(E)} A
 \end{CD}
\end{equation*}
Since $(\Omega^A(U_E),\epsilon',\delta')$ is a left dual for $\Omega^A(U_E^*)$, there
exists $M =(m_{ij}) \in \GL_n(k)$ such that $\epsilon'(w_i \otimes w_j')=m_{ij}$
and $\delta'(1) =\sum_{i,j=1}^n m_{ij}^{-1}w'_i \otimes w_j$. From this we see that
$M=z^tE^{-1}c$ and $E=cM^{-1}z^t$. Hence there exists an algebra map
$f : \mathcal H(E,M^{-1}) \longrightarrow A$ such that $f(u_{ki})=z_{ki}$ and
$f(v_{ki})=c_{ki}$, which is $\mathcal H(E)$-colinear, and is an isomorphism by Proposition \ref{galgroupo}.
\end{proof}
 
\begin{proposition}\label{inj2}
 Let $F,G \in \overline{X}_0(E)$ with $F \in \GL_n(k)$ and $G \in\GL_p(k)$. The following assertions are equivalent.
\begin{enumerate}
 \item The left $\mathcal H(E)$-comodule algebras $\mathcal H(E,F)$ and $\mathcal H(E,G)$ are isomorphic.
\item $F \approx G$, i.e. $n=p$  and $\exists P \in \GL_n(k)$ with 
 $F = PGP^{-1}$.
\end{enumerate}
\end{proposition}

\begin{proof}
The proof is similar to that of Proposition \ref{inj}:
 if $f : \mathcal H(E,F) \longrightarrow \mathcal H(E,G)$ is a left $\mathcal H (E)$-comodule algebra isomorphism, then $n=p$ and there exists $M \in \GL_n(k)$ such that $f(u) = uM$ and $f(v)=uM^{t^{-1}}$, and then we have $G = M^{t^{-1}}FM^t$, so $F \approx G$. Conversely, if there exists $ P \in \GL_n(k)$ with 
 $F = PGP^{-1}$, one defines a left $\mathcal H(E)$-comodule algebra isomorphism  $f : \mathcal H(E,F) \longrightarrow \mathcal H(E,G)$ by $f(u) = u P^t$ and 
$f(v) = v P^{t^{-1}}$.
\end{proof}

Theorem \ref{classHE} follows from these two propositions.

\subsection{Some remarks}
\begin{enumerate}
 \item The problem of classifying the fibre functors
(or Galois objects) on a comodule category is a particular
case of the problem of classifying the module categories over the comodule category
(see \cite{os03}), considered for $\mathcal O_q(\SL_2(k))$
by Etingof and Ostrik (\cite{etos} and \cite{os08}).
Recent nice contributions on the problem of classifying module
categories (and hence Hopf-Galois objects) were done by Mombelli
\cite{momb1,momb2}

\item As we have seen the fibre functor method provides a powerful
method for classifying Hopf-Galois objects, as least if the comodule
category is generated by a few equationnaly well defined morphisms.
It does not seem to work so well for pointed Hopf algebras, or 
at least I could not make it work well.
The best technique for pointed Hopf algebras has been for a long time
 the one by Masuoka \cite{ma}, but recently new techniques have emerged,
by Masuoka again \cite{ma2} or more interestingly yet by the already mentionned work of Mombelli \cite{momb1,momb2}.
\end{enumerate}

\section{Applications to comodule algebras}

In this section we present applications of bi-Galois theory
to the study of comodule algebras. The basic result is the following one.

\begin{theorem}\label{comodalg}
 Let $H$, $L$ be Hopf algebras with a $k$-linear monoidal equivalence
$$
F : {\rm Comod}(H) \cong^\otimes {\rm Comod}(L) 
$$
Then $F$ induces an equivalence of categories between right $H$-comodule algebras
and right $L$-comodule algebras. If $B$ is a right $H$-comodule algebra, the above
monoidal equivalence induces an algebra isomorphism between the respective coinvariant
subalgebras
$B^{{\rm co} H} \cong F(B)^{{\rm co} L}$. 
\end{theorem}

\begin{proof}
 A right $H$-comodule algebra is exactly an algebra in the monoidal category
of right $H$-comodules, thus the above functor induces the announced
category equivalence. 
Let $B$ be a right $H$-comodule algebra.
We have an algebra isomorphism
\begin{align*}
 {\rm Hom}_H(k,B) &\longrightarrow B^{{\rm co} H}\\
f &\longmapsto f(1) 
\end{align*}
where the algebra structure on ${\rm Hom}_H(k,B)$ is given by
$f \cdot g = m_B \circ (f \otimes g)$.
 The above monoidal functor
 induces an algebra isomorphism
$${\rm Hom}_H(k,B) \cong {\rm Hom}_L(F(k),F(B)) \cong{\rm Hom}_L(k,F(B))$$
and hence induces an algebra isomorphism
$B^{{\rm co} H} \cong F(B)^{{\rm co} L}$.
\end{proof}

The cogroupoid picture might be useful to determine exactly
the algebras $F(B)$ and
$F(B)^{{\rm co} L}$ (with $F(B) =B\square_HA$
and $F(B)^{{\rm co} L} = (B\square_HA)^{{\rm co} L}$ for an $H$-$L$-bi-Galois object $A$), for example to get
a presentation by generators and relation from the ones of $B$ and $B^{{\rm co} H}$.
 We present some illustrations.

\subsection{A model comodule algebra for $\B(E)$}
We continue the study of $\B(E)$-comodules by providing explicit models
for the simple $\B(E)$-comodules in the cosemisimple case.
We assume that $k$ has characteristic zero in this subsection.

\begin{definition}
 Let $H$ be a cosemisimple Hopf algebra. A \textbf{model
$H$-comodule algebra} is a right $H$-comodule algebra in which every simple $H$-comodule appears with multiplicity one. 
\end{definition}

It is well known that if $q$ is not a root of unity of order $\geq 3$,
then the quantum plane algebra $k_q[x,y]$ is a model
$\mathcal O_q(\SL_2(k))$-comodule algebra, since the degree grading provides a 
decomposition
$$k_q[x,y] =\bigoplus_{n \in \mathbb N}k_q[x,y]_n$$
and the comodules $k_q[x,y]_n$, $n \in \mathbb N$, provide a complete list
of simple $\mathcal O_q(\SL_2(k))$-comodules (see e.g. \cite{ks}).
It is clear from Theorem \ref{comodalg} that the property of having a model comodule algebra
is monoidally invariant and hence $\B(E)$ will have a model comodule
algebra when it is cosemisimple.

\begin{definition}
 Let $M = (\alpha_{ij}) \in M_m(k)$ and $t \in k$. The algebra $\mathcal A_{M,t}$
is the algebra presented by generators 
$x_1, \ldots , x_m$ submitted to the relation
$\sum_{i,j=1}^m \alpha_{ij}x_ix_j =t$.
\end{definition}

When $$M = E_q^{-1} =  \left(\begin{array}{cc} 0 & -q \\
                          1 & 0\\
       \end{array} \right)$$
the algebra $\mathcal A_{M,0}$ is the quantum plane $k_q[x,y]$, while
that algebra $\mathcal A_{M,1}$ is the quantum Weyl algebra $\mathcal A_1^q(k)$.

\begin{proposition}\label{AEcomodalg}
 Let $E \in \GL_m(k)$. Then $\mathcal A_{E^{-1},t}$ has a right $\B(E)$-comodule
algebra structure defined by
the formula
$$\alpha(x_i) =\sum_{k=1}^m x_k \otimes a_{ki}$$
Let $F \in \GL_n(k)$ be such that 
${\rm tr}(E^{-1}E^t) = {\rm tr}(F^{-1}F^t)$ ($m,n \geq 2$) and consider
the monoidal equivalence $\Theta : {\rm Comod}(\B(E)) \longrightarrow {\rm Comod}(\B(F))$
induced by the $\B(E)$-$\B(F)$-bi-Galois object $\B(E,F)$.
We have an isomorphism of $\B(F)$-comodule algebras $\Theta(\mathcal A_{E^{-1},t}) \cong \mathcal A_{F^{-1},t}$.
\end{proposition}

\begin{proof}
 It is straighforward to check that $\mathcal A_{E^{-1},t}$ has a $\B(E)$-comodule algebra structure given by the above formula. We give two sketches of proof that 
$\Theta(\mathcal A_{E^{-1},t}) \cong \mathcal A_{F^{-1},t}$. In the first proof
we proceed very similarly to the proof of Corollary \ref{noncleft} 
when we showed that $\Theta(V_E) \cong V_F$. It is straightforward to check
that one has an algebra map
\begin{align*}
 \iota_F :\mathcal A_{F^{-1},t} &\longrightarrow \mathcal A_{E^{-1},t} \square_{\B(E)}
\B(E,F) =\Theta(\mathcal A_{E^{-1},t})\\
x_i^F &\longmapsto \sum_{k=1}^m x_k^E\otimes a_{ki}^{E,F}
\end{align*}
which is $\B(F)$-colinear
and the argument of the proof of Corollary \ref{noncleft}
to show that $\iota_F$ is an isomorphism works without any change.

A second sketch of proof is as follows \cite{mov}.
Consider the tensor algebra $T(V_E)= \bigoplus_{l \in \mathbb N} V_E^{\otimes l}$. This
is a $\B(E)$-comodule algebra and we have $\Theta(T(V_E)) 
\cong T(\Theta(V_E)) \cong T(V_F)$ as $\B(F)$-comodule algebras.
As in the proof of Proposition \ref{surj}, consider the $\B(E)$-colinear map
$$\delta_E : k \longrightarrow V_E \otimes V_E, 1 \longmapsto \sum_{k,l}\lambda_{kl}^{-1} v_k^E \otimes v_l^E, \ ({\rm with} \ E^{-1}=(\lambda_{kl}^{-1}))$$
By construction $\mathcal A_{E^{-1},t}=T(V_E)/\langle \delta_E-t1\rangle$
an with an obvious abuse of notation we have 
$\Theta(\delta_E)= \delta_F$ (check this)
and hence 
$$\Theta(\mathcal A_{E^{-1},t})=
\Theta(T(V_E)/\langle \delta_E-t1\rangle) \cong 
T(V_F)/\langle \Theta(\delta_E-t1)\rangle)=
T(V_F)/\langle \delta_F-t1\rangle= \mathcal A_{F^{-1},t}
$$
This concludes the proof.
\end{proof}

\begin{corollary}
 Let $E \in \GL_n(k)$ be such that there exists
$q \in k^*$ with ${\rm tr}(E^{-1}E^t) = -q-q^{-1}$ and $q$ is not a root of unity of order $\geq 3$. Then for $t \in \{0,1\}$, $\mathcal A_{E^{-1},t}$ is a model $\B(E)$-comodule algebra.
\end{corollary}

\begin{proof}
 The proof follows from Corollary \ref{sl}, the previous proposition and
the fact that $\mathcal A_{E_q^{-1},0}=k_q[x,y]$ and $\mathcal A_{E_q^{-1},1}= \mathcal A_{1}^q(k)$ are model $\mathcal O_q(\SL_2(k))$-comodule algebras
(that $\mathcal A_1^q(k)$ is a model $\mathcal O_q(\SL_2(k))$-comodule algebra
can be deduced from the fact that $k_q[x,y]$ is the graded algebra associated to the natural degree filtration on $\mathcal A_1^q(k)$).
\end{proof}

\begin{remark}{\rm
 It follows from the previous results that if $\B(E)$ is cosemisimple,
the simple $\B(E)$-comodule corresponding to $n \in \mathbb N$
is the subspace of degree $n$ elements in $\mathcal A_{E^{-1},0}$.
}
\end{remark}

\subsection{An application to invariant theory}
The considerations of the beginning of the section
allows us to get nice theorems in invariant theory
almost for free. We provide an illustration taken from \cite{bi00}.

The first and most well-known theorem in invariant theory
is the structure theorem for the algebra of symmetric polynomials.
So each time we have a $2$-cocycle deformation of $\mathcal O(S_n)$
(the algebra of functions on the symmetric group $S_n$),
we can derive a twisted version of the classical theorem
on symmetric polynomials. 

We begin by introducing some notation.
We assume that ${\rm char}(k) \not=2$ . Let $n\geq 2$
For $i \in \{1,\ldots,2n\}$, we put 
$$\begin{cases}
i'= i-1 \ {\rm and} \  i^*= i/2 & {\rm if}\  i \ {\rm is} \ {\rm even} \\
i' = i + 1 \ {\rm and} \  i^* = i'/2 & {\rm if}\  i \ {\rm is} \ {\rm odd}
\end{cases}$$
 We have $i'^*= i^*$. Let ${\rm AST}_2(n)$ be the set of AST matrices ${\bf p}=(p_{ij}) \in
  M_n(k)$ with $p_{ii} = 1$ and $p_{ij} = p_{ji} = \pm 1$ for any
    $i$ and $j$.

For ${\bf p} \in {\rm AST}_2(n)$, the algebra $\mathcal O_{\bf p}(S_{2n})$ is the algebra presented by generators
$(x_{ij})_{1\leq i,j \leq 2n}$ and submitted to the relations
$(1 \leq i,j,k,l \leq 2n)$:
\begin{eqnarray}
x_{ij} x_{ik} = \delta_{jk} x_{ij} \quad ; \quad
x_{ji} x_{ki} = \delta_{jk} x_{ji} \quad ; \quad
\sum_{l=1}^{2n} x_{il} = 1 = \sum_{l=1}^{2n} x_{li} \ . \\
(3 + p_{i^*j^*})x_{kj}x_{li} + (1 - p_{i^*j^*})x_{kj}x_{li'}
+ (1 - p_{i^*j^*})x_{kj'}x_{li} +(p_{i^*j^*} -1)x_{kj'}x_{li'} = \nonumber \\
(3 + p_{l^*k^*})x_{li}x_{kj} + (1 - p_{l^*k^*})x_{l'i}x_{kj}
+ (1 - p_{l^*k^*})x_{li}x_{k'j} + (p_{l^*k^*} - 1)x_{l'i}x_{k'j} \ . 
\end{eqnarray}
If ${\bf p}={\bf 1}$ then $\mathcal O_{\bf
  1}(S_{2n})$ is just the algebra of functions on $S_{2n}$
(see \cite{wa98}). It is a direct verification to check that 
$\mathcal O_{\bf p}(S_{2n})$ is a Hopf algebra with
$$\Delta(x_{ij}) = \sum_k x_{ik} \otimes x_{kj}, \quad 
\varepsilon(x_{ij}) = \delta_{ij}, \quad \ S(x_{ij})= x_{ji}$$
We consider the algebra
$k_{\bf p}[x_1, \ldots, x_{2n}]$
presented by generators $x_1, \ldots , x_{2n}$
submitted to the relations
\begin{eqnarray}
4x_ix_j = 
(3 + p_{i^*j^*})x_jx_i + (1 - p_{i^*j^*})x_{j'}x_i
+ (1 - p_{i^*j^*})x_jx_{i'} + (p_{i^*j^*} - 1)x_{j'}x_{i'} \ .
\end{eqnarray} 
It is clear that $k_{\bf 1}[x_1, \ldots, x_{2n}]= k[x_1, \ldots, x_{2n}]$.

Our aim is to prove the following result.

\begin{theorem}\label{invs2n}
 The algebra
$k_{\bf p}[x_1, \ldots, x_{2n}]$ has a right
$\mathcal O_{\bf p}(S_{2n})$-comodule algebra structure
defined by
$\delta(x_i) = \sum_k x_k \otimes x_{ki}$, and
the algebra of coinvariants $k_{\bf p}[x_1, \ldots, x_{2n}]^{{\rm co}  \mathcal O_{\bf p}(S_{2n})}$
is a (commutative) algebra of polynomials on $2n$ variables.
\end{theorem}

The strategy is clear: we have construct a $k$-linear monoidal
equivalence 
$${\rm Comod}(\mathcal O(S_{2n})) \cong^\otimes {\rm Comod}(\mathcal O_{\bf p}(S_{2n}))$$ sending $k[x_1, \ldots, x_{2n}]$
to $k_{\bf p}[x_1, \ldots, x_{2n}]$, and combine the classical theorem on symmetric polynomials with Theorem \ref{comodalg}.
In order to do this we construct an appropriate connected cogroupoid.

Let ${\bf p},{\bf q} \in {\rm AST}_2(n)$. The algebra
$\mathcal O_{{\bf p}, {\bf q}}(S_{2n})$ is the 
algebra presented by generators $(x_{ij})_{1 \leq i,j \leq 2n}$ submitted to the relation:  
\begin{eqnarray}
x_{ij} x_{ik} = \delta_{jk} x_{ij} \quad ; \quad
x_{ji} x_{ki} = \delta_{jk} x_{ji} \quad ; \quad
\sum_{l=1}^{2n} x_{il} = 1 = \sum_{l=1}^{2n} x_{li} \ . \\
(3 + q_{i^*j^*})x_{kj}x_{li} + (1 - q_{i^*j^*})x_{kj}x_{li'}
+ (1 - q_{i^*j^*})x_{kj'}x_{li} +(q_{i^*j^*} -1)x_{kj'}x_{li'} = \nonumber \\
(3 + p_{l^*k^*})x_{li}x_{kj} + (1 - p_{l^*k^*})x_{l'i}x_{kj}
+ (1 - p_{l^*k^*})x_{li}x_{k'j} + (p_{l^*k^*} - 1)x_{l'i}x_{k'j} \ . 
\end{eqnarray}

\begin{lemma} In this lemma we note $\mathcal O_{\mathbf p, \bf q}=\mathcal O_{\mathbf p, \bf q}(S_{2n})$.

\begin{enumerate}
 \item For any ${\bf p}, {\bf q}, {\bf r} \in {\rm AST}_2(n)$, 
there exists an algebra map
\begin{align*}
\Delta_{{\bf p},{\bf q}}^{\bf r} : \mathcal O_{\mathbf p, \bf q} &
\longrightarrow  \mathcal O_{\mathbf p, \bf r} \otimes \mathcal O_{\mathbf r, \bf p}\\
x_{ij},  &\longmapsto \sum_{k=1}^n x_{ik} \otimes x_{kj}
\end{align*}
and for any ${\bf s} \in {\rm AST}_2(n)$, the following diagrams commute
\begin{equation*}
\begin{CD}
\mathcal O_{\mathbf p,\bf q} @>\Delta_{\bf p,\bf q}^{\bf r}>> 
\mathcal O_{\mathbf p, \bf r} \otimes \mathcal O_{\bf r, \mathbf p} \\
@V\Delta_{\bf p,\bf q}^{\bf s}VV
@V\Delta_{\bf p, \bf r}^{\bf s} \otimes 1VV\\
\mathcal O_{\mathbf p, \bf s} \otimes  \mathcal O_{\bf s ,\mathbf q} @>1 \otimes \Delta_{\bf s,\bf q}^{\bf r}>> 
\mathcal O_{\mathbf p, \bf s}\otimes \mathcal O_{\mathbf s, \bf r}  \otimes \mathcal O_{\mathbf r, \bf q}
\end{CD}
\end{equation*}
$$\xymatrix{
\mathcal O_{\mathbf p, \bf q} \ar[d]^{\Delta_{\bf p, \bf q}^{\bf q}} \ar@{=}[rrd]^{\rm id} \\
 \mathcal O_{\mathbf p, \bf q} \otimes \mathcal O_{\mathbf q, \bf q} \ar[rr]^{1 \otimes \varepsilon_{\bf q}} && \mathcal O_{\mathbf p, \bf q}}
\quad
\xymatrix{
\mathcal O_{\mathbf p, \bf q} \ar[d]^{\Delta_{\bf p,\bf q}^{\bf p}} \ar@{=}[rrd]^{{\rm id}} \\
 \mathcal O_{\mathbf p, \bf p} \otimes \mathcal O_{\mathbf p, \bf q} \ar[rr]^{\varepsilon_{\bf p} \otimes 1}&& \mathcal O_{\mathbf p, \bf q}}
$$

\item For any $ \bf p, \bf q \in {\rm AST}(n)$, there exists an algebra map 
\begin{align*}
S_{\bf p, \bf q} :\mathcal O_{\mathbf p, \bf q}(S_{2n}) & \longrightarrow \mathcal O_{\bf q , \mathbf p}(S_{2n})^{\rm op} \\
x  & \longmapsto x^t 
\end{align*}
such that the following diagrams commute
{\footnotesize
$$\xymatrix{
\mathcal O_{\mathbf p,\bf p} \ar[r]^{\varepsilon_{\bf p}} \ar[d]^{\Delta_{\bf p,\bf p}^{\bf q}} &
k \ar[r]^{u}  & \mathcal O_{\mathbf p, \bf q} \\
\mathcal O_{\mathbf p, \bf q} \otimes \mathcal O_{\bf q, \mathbf p} \ar[rr]^{1 \otimes S_{\bf q, \bf p}} && \mathcal O_{\mathbf p, \bf q} \otimes \mathcal O_{\mathbf p, \bf q} \ar[u]^{m}} \quad \quad
\xymatrix{
\mathcal O_{\mathbf p,\bf p} \ar[r]^{\varepsilon_{\bf p}} \ar[d]^{\Delta_{\bf p,\bf p}^{\bf q}} &
k \ar[r]^{u}  & \mathcal O_{\mathbf q,\bf p} \\
\mathcal O_{\mathbf p,\bf q} \otimes \mathcal O_{\mathbf q,\bf p} \ar[rr]^{S_{\bf p,\bf q} \otimes 1} && \mathcal O_{\mathbf q,\bf p} \otimes \mathcal O_{\mathbf q,\bf p} \ar[u]^{m}}
$$}
\end{enumerate}
\end{lemma}

\begin{proof}
This is a straightforward verification. One way to simplify the  computations is as follows. For $i,j,l,k \in \{1, \ldots , 2n\}$ and ${\bf p} \in {\rm ASR}_2(n)$, put
$$R_{ij}^{kl}({\bf p}) = \delta_{i^*l^*}\delta_{j^*k^*}
\left(1 + (-1)^{i-l} + (-1)^{j-k} + (-1)^{i-l+j-k}p_{j^*i^*}\right)$$
With this notation, the second relations defining $\mathcal O_{\bf p, \bf q}$ are
$$\sum_{\alpha,\beta} R_{\alpha\beta}^{kl}({\bf p}) x_{\alpha i} x_{\beta j}
= \sum_{\alpha,\beta} R^{\alpha\beta}_{ij}({\bf q}) x_{k\alpha } x_{l\beta}$$
With this compact notation the verification that $\Delta_{\bf p,\bf q}^{\bf r}$
is well-defined is easy.
\end{proof}

The lemma ensures that we have a cogroupoid whose objects are the elements
of ${\rm AST}_2(n)$ and whose algebras are the above algebras $\mathcal O_{\mathbf p, \bf q}(S_{2n})$. The following lemma and Proposition \ref{connec} ensure that the cogroupoid is connected, and hence that each $\mathcal O_{{\bf p}, {\bf q}}(S_{2n})$
is an $\mathcal O_{{\bf p}}(S_{2n})$-$\mathcal O_{{\bf q}}(S_{2n})$-bi-Galois object.

\begin{lemma}
We have $\mathcal O_{{\bf p}, {\bf 1}}(S_{2n})\not=(0)$ for any ${\bf p} \in {\rm AST}_2(n)$.
\end{lemma}

\begin{proof}
Consider the algebra $A$ presented by generators $t_1, \ldots , t_n$ submitted to the relations $t_i^2=1$, $\forall i$, $t_it_j= p_{ij}t_jt_i$. It is not difficult to check that $A$ is a non-zero algebra (it is isomorphic with the twisted group algebra
$k_\sigma[(\mathbb Z / 2 \mathbb Z)^n]$ for an appropriate 2-cocycle). One then 
checks that there exists an algebra map
\begin{align*}
 \mathcal O_{{\bf p}, {\bf 1}}(S_{mn}) &\longrightarrow A \\
x_{ij} &\longmapsto \frac{\delta_{i^* , j^*}}{2}(1+(-1)^{j-i}t_{i^*}) 
\end{align*}
and hence $\mathcal O_{{\bf p}, {\bf 1}}(S_{mn})$ is a non-zero
algebra.
\end{proof}

\begin{proof}[Proof of Theorem \ref{invs2n}]
 The existence of the announced $\mathcal O_{\bf p}(S_{2n})$-comodule algebra structure on $k_{\bf p}[x_1, \ldots , x_{2n}]$ is left to the reader. Also the reader
will check the existence of an algebra map
\begin{align*}
 k_{\bf p}[x_1, \ldots , x_{2n}] &\longrightarrow k_{\bf q}[x_1, \ldots , x_n]
\square_{\mathcal O_{\bf q}(S_{2n})} \mathcal O_{\bf q, \bf p}(S_{2n}) \\
x_i &\longmapsto \sum_{k=1}^n x_k \otimes x_{ki}
\end{align*}
for any ${\bf p}, {\bf q} \in {\rm AST}_2(n)$. From this he will check, using the technique already used several times (e.g. in the proof of Corollary \ref{noncleft}), 
that the monoidal equivalence associated to the $\mathcal O_{\bf q}(S_{2n})$-$\mathcal O_{\bf p}(S_{2n})$-bi-Galois object $\mathcal O_{\bf q, \bf p}(S_{2n})$ transforms the algebra  $k_{\bf q}[x_1, \ldots , x_n]$  into the algebra $k_{\bf p}[x_1, \ldots , x_n]$. Theorem \ref{comodalg} and the classical theorem on the structure of symmetric polynomials conclude the proof. 
\end{proof}

\section{Yetter-Drinfeld modules}

In this section the previous results and constructions are applied 
to Yetter-Drinfeld modules. We give an application to Brauer groups of Hopf
algebras, while an application to bialgebra cohomology
will be given at the end of the next Section.

Let $H$ be a Hopf algebra. Recall that a (right-right) Yetter-Drinfeld
module over $H$ is a right $H$-comodule and right $H$-module $V$
satisfying the condition
$$(v \leftarrow x)_{(0)} \otimes  (v \leftarrow x)_{(1)} =
v_{(0)} \leftarrow x_{(2)} \otimes S(x_{(1)}) v_{(1)} x_{(3)}$$
The category of Yetter-Drinfeld modules over $H$ is denoted $\yd_H^H$:
the morphisms are the $H$-linear $H$-colinear maps.
Endowed with the usual tensor product of 
modules and comodules, it is a monoidal category.

\subsection{Monoidal equivalence between categories of Yetter-Drinfeld modules}
Recently some authors (\cite{chzh,benperwit}) have used the fact that 
the Yetter-Drinfeld categories of Hopf algebras
that are cocycle deformation of each other are monoidally equivalent.
This result is  a particular case of the following result.

\begin{theorem}\label{coyd}
 Let $H$ and $L$ be two Hopf algebras such that there exists a $k$-linear
monoidal equivalence between ${\rm Comod}(H)$ and ${\rm Comod}(L)$.
Then there exists a a $k$-linear
monoidal equivalence between $\yd_H^H$ and $\yd_L^L$.
\end{theorem}

A sketch of  a quick proof of the theorem is as follows:
the weak center (see the appendix in \cite{scsurv}) of 
the category ${\rm Comod}(H)$ is monoidally
equivalent with the category $\yd_H^H$, and hence a monoidal
equivalence between  ${\rm Comod}(H)$ and ${\rm Comod}(L)$
will induce a monoidal equivalence between $\yd_H^H$ and $\yd_L^L$.

The result that the weak center of ${\rm Comod}(H)$ is $\yd_H^H$
is stated in \cite{scsurv}. It is probably
well-known to many specialists, but to the best of my knowledge
 no complete proof
is avalaible in the litterarure (although it is similar
to the more familiar  module case).
Hence it might be useful to have a direct and more constructive
proof using bi-Galois objects and cogoupoids.
It is the aim of the section to provide such a proof.

\begin{proposition}
 Let $\coc$ be cogroupoid, let $X,Y \in {\rm ob}(\coc)$ and let $V$
be a right $\coc(X,X)$-module.
\begin{enumerate}
 \item $V \otimes \coc(X,Y)$ has a right $\coc(Y,Y)$-module structure
defined by
$$\left(v \otimes a^{X,Y}\right) \leftarrow b^{Y,Y} = v\leftarrow b_{(2)}^{X,X}
\otimes S_{Y,X}(b_{(1)}^{Y,X})a^{X,Y}b_{(3)}^{X,Y}$$ 
Endowed with the right $\coc(Y,Y)$-comodule defined by $1 \otimes \Delta_{X,Y}^Y$,
$V \otimes \coc(X,Y)$ is a Yetter-Drinfeld module over $\coc(Y,Y)$.
\item  If moreover $V$ is Yetter-Drinfeld module, then $V \square_{\coc(X,X)} \coc(X,Y)$ is Yetter-Drinfeld submodule
of $V \otimes \coc(X,Y)$.  
\end{enumerate}
\end{proposition}

Note that when $V=k$ is the trivial comodule, then the action 
of $\coc(Y,Y)$ on $\coc(X,Y)$ is the Myashita-Ulbrich action (see \cite{scsurv}, Definition 2.1.8).

\begin{proof} We freely use Proposition \ref{anti}.
 Let us first check that the above formula defines a right $\coc(Y,Y)$-module structure on $V \otimes \coc(X,Y)$. It is clear that $(v \otimes a^{X,Y}) \leftarrow
1^{Y,Y} = v \otimes a^{X,Y}$. We have
\begin{align*}
 (v \otimes a^{X,Y}) \leftarrow (b^{Y,Y}c^{Y,Y})
& = v\leftarrow(b_{(2)}^{X,X}c_{(2)}^{X,X}) \otimes S_{Y,X}(b_{(1)}^{Y,X}c_{(1)}^{Y,X})a^{X,Y}b_{(3)}^{X,Y}c_{(3)}^{X,Y} \\
& = (v\leftarrow b_{(2)}^{X,X})\leftarrow c_{(2)}^{X,X} \otimes
S_{Y,X}(c_{(1)}^{Y,X})S_{Y,X}(b_{(1)}^{Y,X})a^{X,Y}b_{(3)}^{X,Y}c_{(3)}^{X,Y} \\
& = \left(v\leftarrow b_{(2)}^{X,X} \otimes
S_{Y,X}(b_{(1)}^{Y,X})a^{X,Y}b_{(3)}^{X,Y} \right)
\leftarrow c^{X,X} \\
& = \left(v \otimes a^{X,Y} \leftarrow b^{Y,Y}\right) \leftarrow c^{Y,Y}
\end{align*}
and hence $V \otimes \coc(X,Y)$ is a right $\coc(Y,Y)$-module.
We also have
\begin{align*}
 ((v \otimes a^{X,Y}) \leftarrow b^{Y,Y})_{(0)} & 
\otimes ((v \otimes a^{X,Y}) \leftarrow b^{Y,Y})_{(1)}  \\
& = 
v \leftarrow b_{(2)}^{X,X} \otimes \Delta_{X,Y}^Y(S_{Y,X}(b_{(1)}^{Y,X})a^{X,Y}b_{(3)}^{X,Y}) \\
& =  v \leftarrow b_{(3)}^{X,X} \otimes
S_{Y,X}(b_{(2)}^{Y,X})a_{(1)}^{X,Y}b_{(4)}^{X,Y} \otimes
S_{Y,Y}(b_{(1)}^{Y,Y})a_{(2)}^{Y,Y}b_{(5)}^{Y,Y} \\
& = ((v \otimes a_{(1)}^{X,Y}) \leftarrow b_{(2)}^{Y,Y}) \otimes
 S_{Y,Y}(b_{(1)}^{Y,Y})a_{(2)}^{Y,Y}b_{(3)}^{Y,Y} \\
& = (v \otimes a^{X,Y})_{(0)} \leftarrow b_{(2)}^{Y,Y} \otimes
S_{Y,Y}(b_{(1)}^{Y,Y})(v \otimes a^{X,Y})_{(1)}b_{(3)}^{Y,Y}
\end{align*}
and this shows that  $V \otimes \coc(X,Y)$ is a Yetter-Drindeld module
over $\coc(Y,Y)$.

Assume now that $V$ is a Yetter-Drinfeld module and let us check  that $V \square_{\coc(X,X)} \coc(X,Y)$ is Yetter-Drinfeld submodule
of $V \otimes \coc(X,Y)$. We already know that
it is a subcomodule. Let $\sum_i v_i \otimes a_i^{X,Y} \in V \square_{\coc(X,X)} \coc(X,Y)$.
We must check that
$$(\sum_i v_i \otimes a_i^{X,Y}) \leftarrow b^{Y,Y} = \sum_i
v_i\leftarrow b_{(2)}^{X,X}
\otimes S_{Y,X}(b_{(1)}^{Y,X})a_i^{X,Y}b_{(3)}^{X,Y} \in V \square_{\coc(X,X)} \coc(X,Y)$$ 
We have
\begin{align*}
 \sum_i &
(v_i\leftarrow b_{(2)}^{X,X})_{(0)} \otimes (v_i\leftarrow b_{(2)}^{X,X})_{(1)}
\otimes S_{Y,X}(b_{(1)}^{Y,X})a_i^{X,Y}b_{(3)}^{X,Y}  \\
& =\sum_i
v_{i(0)}\leftarrow b_{(3)}^{X,X} \otimes S_{X,X}(b_{(2)}^{X,X})v_{i(1)}b_{(4)}^{X,X} 
\otimes S_{Y,X}(b_{(1)}^{Y,X})a_i^{X,Y}b_{(5)}^{X,Y}  \\
& = \sum_i
v_{i}\leftarrow b_{(3)}^{X,X} \otimes S_{X,X}(b_{(2)}^{X,X})a_{i(1)}^{X,X}b_{(4)}^{X,X} 
\otimes S_{Y,X}(b_{(1)}^{Y,X})a_{i(2)}^{X,Y}b_{(5)}^{X,Y}  \\
&= \sum_i
v_{i}\leftarrow b_{(2)}^{X,X} \otimes
\Delta_{X,Y}^X(S_{Y,X}(b_{(1)}^{Y,X})a_{i}^{X,Y}b_{(3)}^{X,Y})  
\end{align*}
which gives the result.
\end{proof}

\begin{theorem}\label{expliyetter}
 Let $\coc$ be connected cogroupoid. Then for any $X,Y \in {\rm ob}(\coc)$,
the functor
\begin{align*}
 \yd_{\coc(X,X)}^{\coc(X,X)} & \longrightarrow \yd_{\coc(Y,Y)}^{\coc(Y,Y)} \\
V &\longmapsto V \square_{\coc(X,X)} \coc(X,Y)
\end{align*}
is an equivalence of $k$-linear monoidal categories.
\end{theorem}

\begin{proof}
First we have to check that we indeed have a functor, i.e. that
if $f : V\longrightarrow W$ is a morphism of Yetter-Drinfeld modules,
then $f \otimes 1 :  V \square_{\coc(X,X)} \coc(X,Y) \longrightarrow
W \square_{\coc(X,X)} \coc(X,Y)$ is a morphism of Yetter-Drindeld modules.
We already know it is colinear, and it is not difficult
to check that it is also $\coc(Y,Y)$-linear. Hence we have our functor
\begin{align*}
\mathbb F : \yd_{\coc(X,X)}^{\coc(X,X)} & \longrightarrow \yd_{\coc(Y,Y)}^{\coc(Y,Y)} \\
V &\longmapsto V \square_{\coc(X,X)} \coc(X,Y)
\end{align*}
Consider the symmetrically defined functor
\begin{align*}
\mathbb G : \yd_{\coc(Y,Y)}^{\coc(Y,Y)} & \longrightarrow \yd_{\coc(X,X)}^{\coc(X,X)} \\
V &\longmapsto V \square_{\coc(Y,Y)} \coc(Y,X)
\end{align*}
We know from Theorem \ref{explicit} and its proof that these are inverse equivalences
on comodules and hence we have to check that the comodule isomophisms
$\mathbb G\circ \mathbb F \cong {\rm id}$ and $\mathbb F\circ \mathbb G \cong {\rm id}$
are linear. The first isomorphism is given for $V \in {\rm Comod}(\coc(X,X))$ by
\begin{align*}
\theta_V : V &\longrightarrow (V \square_{\coc(X,X)} \coc(X,Y))\square_{\coc (Y,Y)} \coc(Y,X) \\
v &\longmapsto v_{(0)} \otimes \Delta_{X,X}^Y(v_{(1)})= v_{(0)} \otimes v_{(1)}^{X,Y} \otimes v_{(2)}^{Y,X}
\end{align*}
For $a^{X,X} \in \coc(X,X)$, we have
\begin{align*}
 \theta_V(v\leftarrow a^{X,X}) &= (v\leftarrow a^{X,X})_{(0)} \otimes \Delta_{X,X}^Y((v\leftarrow a^{X,X})_{(1)}) \\
&= v_{(0)}\leftarrow a_{(2)}^{X,X} \otimes 
\Delta_{X,X}^Y(S_{X,X}(a_{(1)}^{X,X})v_{(1)} a_{(3)}^{X,X}) \\
&= v_{(0)}\leftarrow a_{(3)}^{X,X} \otimes 
S_{Y,X}(a_{(2)}^{Y,X})v_{(1)}^{X,Y} a_{(4)}^{X,Y} \otimes
S_{X,Y}(a_{(1)}^{X,Y})v_{(2)}^{Y,X} a_{(5)}^{Y,X}
\\
& = ((v_{(0)} \otimes v_{(1)}^{X,Y})\leftarrow a_{(2)}^{Y,Y}) \otimes 
S_{X,Y}(a_{(1)}^{X,Y})v_{(2)}^{Y,X} a_{(3)}^{Y,X} \\
& = (v_{(0)} \otimes v_{(1)}^{X,Y} \otimes v_{(2)}^{Y,X}) \leftarrow a^{X,X}
= \theta_V(v)\leftarrow a^{X,X}
\end{align*}
So $\theta_V$ is $\coc(X,X)$-linear, and we have
inverse equivalences of Yetter-Drinfeld module categories.

We have to check that $\mathbb F$ is monoidal. It is clear that
$k \cong \mathbb F(k)$ is $\coc(Y,Y)$-linear and hence
it remains to check that for $V,W \in \yd_{\coc(X,X)}^{\coc(X,X))}$ the 
$\coc(Y,Y)$-colinear isomorphism
\begin{align*}
 \widetilde{F}_{V,W} :(V \square_{\coc(X,X)} \coc(X,Y)) \otimes
(W \square_{\coc(X,X)} \coc(X,Y)) 
&\longrightarrow (V\otimes W) \square_{\coc(X,X)} \coc(X,Y) \\
(\sum_i v_i \otimes a_i^{X,Y}) \otimes
(\sum_j w_j \otimes b_j^{X,Y}) &\longmapsto
\sum_{i,j} v_i \otimes w_j\otimes a_i^{X,Y}b_j^{X,Y}
\end{align*}
is $\coc(Y,Y)$-linear. For $c^{Y,Y} \in \coc(Y,Y)$ we have
\begin{align*}
 \widetilde{F}_{V,W}&\left(\left((\sum_{i} v_i \otimes a_i^{X,Y}) \otimes
(\sum_j w_j \otimes b_j^{X,Y})\right)\leftarrow c^{Y,Y} \right)\\
 &= \widetilde{F}_{V,W}\left(\sum_{i,j} v_i \leftarrow c_{(2)}^{X,X} \otimes S_{Y,X}(c_{(1)}^{Y,X})a_i^{X,Y}c_{(3)}^{X,Y} \otimes
w_j \leftarrow c_{(5)}^{X,X}\otimes S_{Y,X}(c_{(4)}^{Y,X})b_j^{X,Y}c_{(6)}^{X,Y}\right)
\\
 &= \sum_{i,j} v_i \leftarrow c_{(2)}^{X,X} \otimes
w_j \leftarrow c_{(5)}^{X,X}\otimes S_{Y,X}(c_{(1)}^{Y,X})a_i^{X,Y}c_{(3)}^{X,Y}
 S_{Y,X}(c_{(4)}^{Y,X})b_j^{X,Y}c_{(6)}^{X,Y}
\\
&= \sum_{i,j} v_i \leftarrow c_{(2)}^{X,X} \otimes
w_j \leftarrow c_{(3)}^{X,X}\otimes S_{Y,X}(c_{(1)}^{Y,X})a_i^{X,Y}
 b_j^{X,Y}c_{(4)}^{X,Y}
\\
& = \left(\sum_{i,j} v_i \otimes w_j\otimes a_i^{X,Y}b_j^{X,Y}\right) \leftarrow c^{Y,Y} \\
& = \widetilde{F}_{V,W}\left(\left((\sum_{i} v_i \otimes a_i^{X,Y}) \otimes
(\sum_j w_j \otimes b_j^{X,Y})\right)\right)\leftarrow c^{Y,Y} 
\end{align*}
This finishes the proof.
\end{proof}

Theorem \ref{coyd} now follows from Theorem \ref{schacog} and Theorem \ref{expliyetter}.

\subsection{Application to Brauer groups of Hopf algebras} 
In this section we briefly indicate applications
of the previous considerations to Brauer groups of Hopf algebras.
 
Let $H$ be a Hopf algebra. Recall that for 
$V, W \in \yd_H^H$, the linear map
\begin{align*}
 c_{V,W} : V \otimes W &\longrightarrow W \otimes V \\
v \otimes w &\longmapsto w_{(0)} \otimes v \leftarrow w_{(1)}
\end{align*}
defines a prebraiding on $\yd_H^H$ (see e.g. the appendix in \cite{scsurv})
which is a braiding if the antipode of $H$ is bijective.

The Brauer group ${\rm Br}(H)$ of a Hopf algebra $H$ with bijective antipode is defined as the Brauer group of the braided category of finite-dimensional Yetter-Drinfeld modules over $H$ \cite{voz}.

Chen and Zhang \cite{chzh} observed that for $\sigma \in Z^2(H)$ there is an equivalence
of braided monoidal categories between $\yd^H_H$ and $\yd_{H^\sigma}^{H^{\sigma}}$
and the Brauer group of a braided monoidal category clearly being
an invariant of the braided category, this yields a group isomorphism
${\rm Br}(H)\cong {\rm Br}(H^\sigma)$.
To generalize the result to general monoidal equivalences
${\rm Comod}(H) \cong^\otimes  {\rm Comod}(L)$, one may use the center
argument or it is enough to check that
the monoidal equivalence of Theorem \ref{expliyetter} is a braided equivalence.

\begin{proposition}
 Let $\coc$ be connected cogroupoid. Then for any $X,Y \in {\rm ob}(\coc)$,
the monoidal equivalence
$$
 \yd_{\coc(X,X)}^{\coc(X,X)} \cong^\otimes \yd_{\coc(Y,Y)}^{\coc(Y,Y)}$$
in Theorem \ref{expliyetter} is an equivalence of prebraided categories.
\end{proposition}

\begin{proof} We have to check that the diagram
\begin{equation*}
\begin{CD}
  (V \square_{\coc(X,X)} \coc(X,Y)) \otimes
(W \square_{\coc(X,X)} \coc(X,Y)) 
@>\widetilde{F}_{V,W}>> (V\otimes W) \square_{\coc(X,X)} \coc(X,Y) \\
@Vc_{-,-}VV  @Vc_{V,W}\otimes 1VV \\
(W \square_{\coc(X,X)} \coc(X,Y)) \otimes
(V \square_{\coc(X,X)} \coc(X,Y)) 
@>\widetilde{F}_{W,V}>> (W\otimes V) \square_{\coc(X,X)} \coc(X,Y)
\end{CD}
\end{equation*}
is commutative.
We have
\begin{align*}
 & \widetilde{F}_{V,W}\circ c_{-,-}\left((\sum_{i} v_i \otimes a_i^{X,Y}) \otimes
(\sum_j w_j \otimes b_j^{X,Y})\right) \\
& = \widetilde{F}_{V,W}\left(\sum_{i,j} 
(w_j \otimes b_j^{X,Y})_{(0)} \otimes (v_i \otimes a_i^{X,Y}) \leftarrow
(w_j\otimes b_j^{X,Y})_{(2)} \right) \\
& = \widetilde{F}_{V,W}\left(\sum_{i,j} 
w_j \otimes b_{j(1)}^{X,Y} \otimes (v_i \otimes a_i^{X,Y}) \leftarrow
b_{j(2)}^{Y,Y} \right) \\
& = \widetilde{F}_{V,W}\left(\sum_{i,j} 
w_j \otimes b_{j(1)}^{X,Y} \otimes v_{i}\leftarrow b_{j(3)}^{X,X} \otimes S_{Y,X}(b_{j(2)}^{Y,X})a_i^{X,Y}b_{j(4)}^{X,Y}\right) \\
& = \sum_{i,j} 
w_j \otimes v_{i}\leftarrow b_{j(3)}^{X,X} \otimes b_{j(1)}^{X,Y}S_{Y,X}(b_{j(2)}^{Y,X})a_i^{X,Y}b_{j(4)}^{X,Y} \\
& = \sum_{i,j} 
w_j \otimes v_{i}\leftarrow b_{j(1)}^{X,X} \otimes a_i^{X,Y}b_{j(2)}^{X,Y}
= \sum_{i,j} 
w_{j(0)} \otimes v_{i}\leftarrow w_{j(1)}^{X,X} \otimes a_i^{X,Y}b_{j}^{X,Y} \\
& =(c_{V,W}\otimes 1) \circ \widetilde{F}_{V,W}\left((\sum_{i} v_i \otimes a_i^{X,Y}) \otimes
(\sum_j w_j \otimes b_j^{X,Y})\right)
\end{align*}
and we are done.
\end{proof}

\begin{corollary}
 Let $H$, $L$ be  Hopf algebras. If ${\rm Comod}(H) \cong^\otimes  {\rm Comod}(L)$ then ${\rm Br}(H) \cong {\rm Br}(L)$.
\end{corollary}

\section{Homological algebra applications}

This section is devoted to various applications in homological algebra.
The first two subsections deal with Hochschild (co)homology of Hopf-Galois objects and comodule algebras. We do not 
recall the definitions related to Hochschild homology
(see \cite{kascons} for a concise introduction and the references therein).
We will not try to give the more complete and strongest results, but we 
will concentrate on simple concrete applications. 
The last section deals with bialgebra cohomology introduced by Gerstenhaber and Schack
\cite{gs1,gs2}.

\subsection{Homology of Hopf-Galois objects} In this subsection we explain
how to relate the (co)homology of a Hopf-Galois object to the (co)homology of 
the underlying Hopf algebra. A deeper study is done in \cite{ste}.  

\begin{theorem}\label{hopfgaloishomology}
 Let $H$ be a Hopf algebra and let $A$ be left $H$-Galois object.
There exists two functors $','' : {\rm Bimod}(A) \longrightarrow {\rm Bimod}(H)$,  such that for any $A$-bimodule $M$ we have
$$\forall n\geq 0, \ H_n(A,M) \cong H_n(H,M') \ {\rm and} \
H^n(A,M) \cong H^n(H,M'')$$ 
\end{theorem}

We begin with a lemma.

\begin{lemma}
 Let $H$ be a Hopf algebra and let $A$ be left $H$-Galois object.
\begin{enumerate}
 \item Let $M$ be a left $A$-module. For any $n \geq 0$, the linear map
\begin{align*}
  A^{\otimes n} \otimes M&\longrightarrow H^{\otimes n} \otimes M \\
a_1 \otimes \cdots \otimes a_n \otimes m &\longmapsto a_{1(-1)} \otimes \cdots \otimes a_{n(-1)}
\otimes a_{1(0)} \cdots a_{n(0)}\cdot m
\end{align*}
is an isomorphism.
\item Let $M$ be a right $A$-module. For any $n \geq 0$, we have a linear isomorphism
$$ \Hom_k(A^{\otimes n}, M) \cong \Hom_k(H^{\otimes n} , M)$$
\end{enumerate}
\end{lemma}

\begin{proof}
By Theorem \ref{galoiscog}
we can assume that $H=\coc(X,X)$ and $A= \coc(X,Y)$ for objects $X,Y$
of a connected cogroupoid $\coc$. If $M$ is a left $A$-module,
the inverse isomorphism
$H^{\otimes n} \otimes M\rightarrow A^{\otimes n} \otimes M$
is given by
$$a_1^{X,X} \otimes \cdots \otimes a_n^{X,X} \otimes m
\longmapsto a_{1(1)}^{X,Y} \otimes \cdots \otimes a_{n(1)}^{X,Y} \otimes 
S_{Y,X}(a_{1(2)}^{Y,X}\cdots a_{n(2)}^{Y,X})\cdot m$$
Assume now that $M$ is a right $A$-module. Consider the linear map
\begin{align*}
 \Hom_k(A^{\otimes n}, M) &\longrightarrow \Hom_k(H^{\otimes n} , M) \\
f &\longmapsto \widehat{f}, \ \widehat{f}(a_1^{X,X} \otimes \cdots \otimes a_n^{X,X})
= f(a_{1(1)}^{X,Y} \otimes \cdots \otimes a_{n(1)}^{X,Y}) \cdot
S_{Y,X}(a_{1(2)}^{Y,X}\cdots a_{n(2)}^{Y,X})
\end{align*}
The inverse is given by $g \longmapsto \tilde{g}$, where
$$\tilde{g}(a_1^{X,Y} \otimes \cdots \otimes a_n^{X,Y})=
g(a_{1(1)}^{X,X} \otimes \cdots \otimes a_{n(1)}^{X,Y})\cdot
a_{1(2)}^{X,Y} \cdots  a_{n(2)}^{X,Y}
$$
and we have our result.
\end{proof}

\begin{proof}[Proof of Theorem \ref{hopfgaloishomology}]
Again we can assume that $H=\coc(X,X)$ and $A= \coc(X,Y)$ for objects $X,Y$
of a connected cogroupoid $\coc$. Let $M$ be an $A$-bimodule.
It is straighforward to check that one defines an $H$-bimodule structure
on $M$ by letting
$$a^{X,X}\cdot m = \varepsilon_X(a^{X,X})m, 
 \quad
m\cdot a^{X,X} = S_{Y,X}(a^{Y,X}_{(2)})\cdot m \cdot a_{(1)}^{X,Y}
$$
This $H$-bimodule is denoted $M'$, and we get a functor
${\rm Bimod}(A) \longrightarrow {\rm Bimod}(H)$, $M \longmapsto M'$.

Similarly there is an $H$-bimodule structure
on $M$ defined by
$$a^{X,X}\cdot m = a^{X,Y}_{(1)}\cdot m \cdot S_{Y,X}(a_{(2)}^{Y,X}), \quad
m\cdot a^{X,X} = \varepsilon_X(a^{X,X})m
$$
This $H$-bimodule is denoted $M''$, and we get a second functor
${\rm Bimod}(A) \longrightarrow {\rm Bimod}(H)$, $M \longmapsto M''$.
It is not difficult to check that 
the linear isomorphisms of the previous lemma induce
isomorphisms between the standard complexes
 $$C_*(A,M) \cong C_*(H,M'), \quad C^*(A,M) \cong C^*(H,M'')$$
defining Hochschild homology and cohomology respectively.
See \cite{kascons} (note that for $C_*(A,M)$ one has to switch
 $M \otimes A^{\otimes *}$ to $A^{\otimes *}\otimes M$).
\end{proof}

\begin{corollary}
Let $q \in k^*$ and let $F \in \GL_n(k)$ satisfying
${\rm tr}(F^{-1}F^t)=-q-q^{-1}$. Then
 for any $p>3$ and any $\mathcal B(E_q,F)$-bimodule $M$, we have $H_p(\mathcal B(E_q,F),M)=(0) = H^p(\mathcal B(E_q,F),M)$.
\end{corollary}

\begin{proof}
 It is known that for any $\mathcal O_q(\SL_2(k))$-bimodule $M$ one has $H_p(\mathcal O_q(\SL_2(k)),M)=(0) = H^p(\mathcal O_q(\SL_2(k)),M)$ if $p>3$ (see e.g. \cite{hk}).
Hence since $\mathcal B(E_q) = \mathcal O_q(\SL_2(k))$, the result follows from Corollary \ref{noncleft} and Theorem
\ref{hopfgaloishomology}.
\end{proof}

\subsection{Equivariant resolutions for comodules algebras}

In this subsection we show how to use monoidal equivalences
to get informations on the Hochschild (co)homology of a comodule algebra, by
transporting appropriate resolutions.
We apply this to the previous algebras $\mathcal A_{F,t}$ (with $F$ invertible) to get the 
 following result. The techniques are taken from \cite{mov,os08}, where more general algebras are studied.

\begin{theorem}\label{homology}
 Let $F \in \GL_n(k)$ and $t\in \{0,1\}$. Then for any $\mathcal A_{F,t}$-bimodule $M$ and any $m > 2$
we have $H_m(\mathcal A_{F,t},M)=(0)=H^m(\mathcal A_{F,t},M)$.
\end{theorem}

The idea is that we can transport well-known resolutions
for the quantum plane $k_q[x,y]$ and for the quantum Weyl algebra $\mathcal A_1^q(k)$ to $\mathcal A_{F,t}$ through a monoidal equivalence.

We begin with some generalities. If $A$ is any algebra, we denote by $A^e$ the algebra $A \otimes A^{\rm op}$. The categories of $A$-bimodules and of $A^e$-modules are identified in the usual manner.

\begin{definition}
 Let $H$ be a Hopf algebra and let $A$ be a right $H$-comodule algebra.
An \textbf{$H$-equivariant free resolution of $A$} is an exact sequence
{\tiny
\begin{equation*}
 \begin{CD}
\cdots  A \otimes V_{n+1} \otimes A @>d_{n+1}>>
 A \otimes V_n \otimes A @>d_n>> \cdots A \otimes V_2 \otimes A @>d_2>>  A \otimes V_1 \otimes A @>d_1>> A\otimes A @>d_0=m>> A @>>> 0
 \end{CD}
\end{equation*}}
where $\forall n\geq 1$, $V_n$ is an $H$-comodule and
$d_n$ is a left $A$-linear, right $A$-linear and right $H$-colinear map.
\end{definition}

An $H$-equivariant free resolution of $A$ is in particular a resolution of $A$ by free $A^e$-modules.



\begin{proposition}\label{transporeso}
 Let $H$, $L$ be Hopf algebras and assume that there exists a $k$-linear monoidal equivalence  $F : {\rm Comod}(H) \cong^\otimes {\rm Comod}(L)$.
Let $A$ be a right $H$-comodule algebra. Then any $H$-equivariant
free resolution of $A$ 
{\tiny
\begin{equation*}
 \begin{CD}
\cdots  A \otimes V_{n+1} \otimes A @>d_{n+1}>>
 A \otimes V_n \otimes A @>d_n>> \cdots A \otimes V_2 \otimes A @>d_2>>  A \otimes V_1 \otimes A @>d_1>> A\otimes A @>d_0=m>> A @>>> 0
 \end{CD}
\end{equation*}}
induces a $L$-equivariant free resolution of $F(A)$
{\footnotesize
\begin{equation*}
 \begin{CD}
\cdots F(A) \otimes F(V_{n+1})\otimes F(A)@>>>
 F(A) \otimes F(V_n) \otimes F(A)@>>> \cdots \\ \cdots  F(A) \otimes F(V_1) \otimes F(A)@>>> F(A) \otimes  F(A)@>>> F(A) @>>> 0
 \end{CD}
\end{equation*}}
\end{proposition}

\begin{proof}
 The $L$-equivariant free resolution of $F(A)$ is defined through the commutative diagram
{\footnotesize
\begin{equation*}
 \begin{CD}
 \cdots F(A) \otimes F(V_{n+1}) \otimes F(A)@>>>
 F(A) \otimes F(V_n) \otimes F(A)\cdots    @>>> F(A) \otimes F(A)@>>> F(A)@>>> 0 \\
@V\widetilde{F} VV @V\widetilde{F} VV @V\widetilde{F}_{A,A}VV @|  \\
\cdots F(A \otimes V_{n+1}\otimes A)@>F(d_{n+1})>>
 F(A \otimes V_n\otimes A) \cdots    @>F(d_1)>> F(A\otimes A) @>F(m)>> F(A) @>>> 0
 \end{CD}
\end{equation*}}
where the $\widetilde{F}$'s denote the monoidal constraints of $F$
(the maps on the top of the diagram are $F(A)$-linear by definition of the algebra structure on $F(A)$).
\end{proof}



We now apply these considerations to the algebras $\mathcal A_{F,t}$
by using the algebras $k_q[x,y]$ and  $\mathcal A_1^q(k)$. 
These algebras are
$\mathcal O_q(\SL_2(k))$-comodule algebras and the associated classical Kozsul type
resolution are 
$\mathcal O_q(\SL_2(k))$-equivariant.

\begin{proposition}\label{kqreso}
For $A_{q,t} = \mathcal A_{E_q^{-1},t}$, $t \in\{0,1\}$, 
$A_{q,t} = k\langle x,y\rangle /(yx-qxy=t1)$, we have an
$\mathcal O_q(\SL_2(k))$-equivariant free resolution of $A_{q,t}$
\begin{equation*}
\begin{CD}
0 @>>> A_{q,t} \otimes A_{q,t} @>>> A_{q,t} \otimes V \otimes A_{q,t}@>>> A_{q,t} \otimes A_{q,t} @>>> A_{q,t} @>>> 0
\end{CD}
\end{equation*}
where $V=kx\oplus ky$.
\end{proposition}

\begin{proof}
The resolution is
 \begin{equation*}
\begin{CD}
0 @>>> A_{q,t} \otimes A_{q,t} @>\gamma>> A_{q,t} \otimes V \otimes A_{q,t}@>m\otimes 1-1 \otimes m>> A_{q,t}\otimes A_{q,t} @>m>> A_{q,t} @>>> 0
\end{CD}
\end{equation*}
where $\gamma$ is the unique $A_{q,t}$-bimodule map such that 
$$\gamma(1\otimes 1)=
-qx\otimes y \otimes 1 -q 1 \otimes x \otimes y 
+ y \otimes x \otimes 1 + 1 \otimes y \otimes x$$
See \cite{kas92,wam} for the verification that this complex of $A_{q,t}$-bimodules is exact.
It is straightforward to check that $\gamma$ is $\mathcal O_q(\SL_2(k))$-colinear, and 
since the multiplication of a comodule algebra is colinear, we conclude that 
we indeed have an 
$\mathcal O_q(\SL_2(k))$-equivariant resolution of $A_{q,t}$.
\end{proof}

\begin{proof}[Proof of Theorem \ref{homology}]
 First assume that there exists $q \in k^*$ such that
${\rm tr}(F^{-1}F^t)=-q-q^{-1}$, so that $\mathcal B(E_q, F^{-1})$
is a $\mathcal O_q(\SL_2(k))$-$\B(F^{-1})$-bi-Galois object. The result follows by using
Propositions \ref{AEcomodalg}, \ref{kqreso}, \ref{transporeso} and the definition of Hochshild (co)homology by Tor and Ext.

Otherwise let $k\subset k'$ a field extension such that there exists
$q \in k'^*$ such that  ${\rm tr}(F^{-1}F^t)=-q-q^{-1}$.
Propositions  \ref{kqreso} and \ref{transporeso} yield
a $k' \otimes \B(F^{-1})$-equivariant free resolution of $k'\otimes \mathcal A_{F,t}$.
Writing explicitely this resolution
we see that it is obtained  by tensoring an exact sequence of the following type
with $k'$
\begin{equation*}
\begin{CD}
0 @>>>  \mathcal A_{F,t}\otimes \mathcal A_{F,t}@>>> \mathcal A_{F,t} \otimes V_F \otimes \mathcal A_{F,t}@>>>  \mathcal A_{F,t} \otimes \mathcal A_{F,t}@>>> \mathcal A_{F,t} @>>> 0
\end{CD}
\end{equation*}
which is therefore a resolution of $\mathcal A_{F,t}$ by free $\mathcal A_{F,t}^e$-modules, and we are done.
\end{proof}

\begin{remark}{\rm 
 The algebras $\mathcal A_{F,0}$ are studied directly by Dubois-Violette in \cite{dv} by using Koszul algebras techniques.}
\end{remark}

\begin{remark}{\rm 
Let us briefly explain how one can also use a monoidal equivalence to deduce the 
well-known resolution for $k_q[x,y]$ of Proposition \ref{kqreso} from the even
more well-known  resolution for  $k[x,y]$. For this consider the resolution
of Proposition \ref{kqreso} at $q=1$ slightly modified as follows:
{\footnotesize
\begin{equation*}
\begin{CD}
0 @>>> k[x,y] \otimes {\det} \otimes k[x,y]@>>> k[x,y] \otimes V \otimes k[x,y] @>>> k[x,y] \otimes k[x,y]@>>> k[x,y] @>>> 0
\end{CD}
\end{equation*}}
where $\det$ is the one dimensional $\mathcal O(\GL_2(k))$-comodule
associated to the determinant. This resolution is formed by
$\mathcal O(\GL_2(k))$-colinear maps. Consider now the
${\rm AST}$-matrix 
$$\bf q  =\left(\begin{array}{cc} 1 & q \\
                          q^{-1} & 1\\
       \end{array} \right)$$
Put $\mathcal O_{q,q^{-1}}(\GL_2(k))=\mathcal O_{\bf q}(\GL_2(k))$
(see Subsection 3.4). The reader will check that the monoidal equivalence
of Corollary \ref{multimonoidal} transforms the $\mathcal O(\GL_2(k))$-comodule algebra
$k[x,y]$ into the $\mathcal O_{q,q^{-1}}(\GL_2(k))$-comodule algebra $k_q[x,y]$.
Hence one gets the resolution of Proposition \ref{kqreso} by using Proposition \ref{transporeso}.

Of course we could not use $\mathcal O_q(\SL_2(k))$ to get the resolution 
because there
does not exists a monoidal equivalence bewteen $\mathcal O(\SL_2(k))$ and
$\mathcal O_q(\SL_2(k))$ if $q \not = 1$.
}
\end{remark}


\subsection{Application to bialgebra cohomology}
The cohomology of a bialgebra was introduced by Gerstenhaber and 
Schack \cite{gs2,gs1}: it is defined by means of an explicit bicomplex whose
arrows are modelled on the Hochschild complex of the underlying algebra
and columns are modelled on the Cartier complex of the underlying coalgebra.
If $H$ is a Hopf algebra, let us denote by $H^*_{b}(H)$
the resulting cohomology. 
Taillefer \cite{tai04} proved that
$$H_b^*(H) \cong {\rm Ext}^*_{\mathcal M(H)}(H,H)$$
where $\mathcal M(H)$ is the category of Hopf bimodules over $H$.
Combined with the monoidal equivalence between Hopf bimodules and Yetter-Drinfeld modules \cite{sc94}, this yields an isomorphism 
$$H_b^*(H) \cong {\rm Ext}^*_{\yd_H^H}(k,k)$$
Moreover it is proved in \cite{tai01} that if $H$ and $L$ are Hopf algebras
such that ${\rm Mod}(H) \cong^\otimes  {\rm Mod}(L)$, then
$H_b^*(H) \cong H_b^*(L)$. The proof is done 
by constructing a monoidal equivalence between 
$\mathcal M(H)$ and $\mathcal M(L)$ from the given monoidal
equivalence between ${\rm Mod}(H)$ and ${\rm Mod}(L)$
(the construction was done by using the dual notion of Hopf-Galois
object from \cite{sua}) and by using the ${\rm Ext}$-description.
Of course the argument works if one starts from a monoidal equivalence between
comodule categories. It is worth  summarizing this  in the following statement.

\begin{theorem}
 Let $H$ and $L$ be two Hopf algebras. Assume that  
$${\rm Mod}(H) \cong^\otimes  {\rm Mod}(L) \ {\rm or} \ {\rm Comod}(H) \cong^\otimes  {\rm Comod}(L)$$ Then $H_b^*(H) \cong H_b^*(L)$.
\end{theorem}

Note that explicit computations
of $H_b^*(H)$ are known only in very few cases: see \cite{tai07}
and the references therein.






\end{document}